\documentclass[11pt]{amsart}

\usepackage[english]{babel}

\usepackage{epigamath}
\usepackage{microtype}

\numberwithin{equation}{section}

\usepackage[utf8]{inputenc}
\usepackage[T1]{fontenc}
\usepackage{amsxtra}
\usepackage{amscd}
\usepackage{amsthm}
\usepackage{amsfonts}
\usepackage{eucal}
\usepackage[all]{xy}
\usepackage{graphicx}
\usepackage{comment}
\usepackage{epsfig}
\usepackage{psfrag}
\usepackage{mathrsfs}
\usepackage{amscd}
\usepackage{rotating}
\usepackage{lscape}
\usepackage{amsbsy}
\usepackage{verbatim}
\usepackage{moreverb}
\usepackage{url}
\usepackage{color}
\usepackage{xcolor}
\usepackage{stmaryrd}
\usepackage{enumitem}

\newcommand{\nc}{\newcommand}
\nc{\renc}{\renewcommand}

\nc\restr[2]{{ 
  \left.\kern-\nulldelimiterspace    #1  
  \vphantom{\big|}  
  \right|_{#2}  
  }}

\newtheorem{prop-intro}{Proposition}
\newtheorem{thm-intro}[prop-intro]{Theorem}
\newtheorem{cor}[subsubsection]{Corollary}
\newtheorem{lem}[subsubsection]{Lemma}
\newtheorem{prop}[subsubsection]{Proposition}

\newtheorem{construction}[subsubsection]{Construction}
\newtheorem{nota}[subsubsection]{Notation}

\newtheorem{thm}[subsubsection]{Theorem}

\newtheorem{notation}[subsubsection]{Notation}

\renc{\sec}{\section}
\nc{\ssec}{\subsection}
\nc{\sssec}{\subsubsection}

\theoremstyle{definition}
\newtheorem{defi}[subsubsection]{Definition}

\newtheorem{rem}[subsubsection]{Remark}

\nc{\on}{\operatorname}

\nc\wt{\widetilde}
\nc\wh{\widehat}
\nc\ol{\ov}
\nc{\oc}[1]{{\overset{\circ}{#1}}}
\nc{\ov}[1]{{\overline{#1}}}
\nc{\isor}[1]{{\xrightarrow[\raisebox{0.25 em}{\smash{\ensuremath{\sim}}}]{#1}}}
\nc{\modmod}{/ \! \! /}

\nc{\mc}{\mathcal}
\nc{\mf}{\mathfrak}
\nc{\mr}{\mathrm}
\nc{\mb}{\mathbb}
\nc{\mbf}{\mathbf}
\nc{\ms}{\mathscr}

\nc{\R}{{\mathbb R}}
\nc{\Z}{{\mathbb Z}}
\nc{\N}{{\mathbb N}}
\nc{\C}{{\mathbb C}}
\nc{\Q}{{\mathbb Q}}

\nc{\Fq}{{\mathbb F}_q}
\nc{\Fl}{{\mathbb F}_\ell}
\nc{\Fqbar}{\ol{{\mathbb F}_q}}
\nc{\Flbar}{\ol{{\mathbb F}_\ell}}
\nc{\Zl}{{\mathbb Z}_\ell}
\nc{\Zlbar}{\ol{{\mathbb Z}_\ell}}
\nc{\Ql}{E}
\nc{\Qlbar}{\ol{{\mathbb Q}_\ell}}
\nc{\hl}{\overset{\leftarrow}h{}}
\nc{\hr}{\overset{\rightarrow}h{}}
\nc{\Gr}{{\on{Gr}}}
\nc{\Hecke}{\on{Hecke}}
 \nc{\Hom}{\on{Hom}}
 \nc{\Coker}{\on{Coker}}
 \nc{\Ker}{\on{Ker}}
 \nc{\Lie}{\on{Lie}}
\nc{\Loc}{\on{Loc}}
\nc{\Pic}{\on{Pic}}
\nc{\Bun}{\on{Bun}}
\nc{\IC}{\on{IC}}
\nc{\Aut}{\on{Aut}}
\nc{\Perv}{\on{Perv}}
\nc{\pos}{{\on{pos}}}
\nc{\Sym}{\on{Sym}}

\nc{\ta} {{}^\tau}
\nc {\tu}[1]{{}^{\tau^{#1}}\!}

\nc{\Id}{\on{Id}}
\nc{\Fil}{\on{Fil}}
\nc{\pr}{\on{pr}}
\nc{\Res}{\on{Res}}
\nc{\cusp}{\on{cusp}}
\nc{\Frob}{\on{Frob}}
\nc{\diag}{\Delta}
\nc{\gr}{\on{gr}}
\nc{\Inj}{\on{Inj}}
\nc{\Bl}{\on{Bl}}
\nc{\dem}{\noindent {\bf Proof. }}
\nc{\cqfd}{{\ }\hfill $\square$ \vskip 1mm}
\nc{\s}[1]{\langle #1 \rangle}
\nc{\Cht}{\on{Cht}}
\nc{\isom}{\overset {\thicksim}{\to}}
\nc{\sm}{\smallsetminus}

\EpigaVolumeYear{4}{2020} \EpigaArticleNr{6} \ReceivedOn{June 6,
2019}
\InFinalFormOn{}
\AcceptedOn{April 1, 2020}

\title{Finiteness of cohomology groups of stacks of shtukas as modules over Hecke algebras, and applications}
\titlemark{Cohomologies of stacks of shtukas}

\author{Cong Xue}
\address{Cong Xue: DPMMS, Centre for Mathematical Sciences, Wilberforce Road, Cambridge, CB3 0WB, UK}
\email{cx233@cam.ac.uk}

\authormark{C. Xue}

\AbstractInEnglish{In this paper we prove that the cohomology groups with compact support of stacks of shtukas are modules of finite type over a Hecke algebra. As an application, we extend the construction of excursion operators, defined in \cite{vincent} on the space of cuspidal automorphic forms, to the space of automorphic forms with compact support. This gives the Langlands parametrization for some quotient spaces of the latter, which is compatible with the constant term morphism.}

\MSCclass{14G35; 11F70}

\KeyWords{Stacks of shtukas; $\ell$-adic cohomology; Hecke algebra}

\TitleInFrench{Finitude des groupes de cohomologie des champs de shtukas comme modules sur les alg\`ebres de Hecke et applications}

\AbstractInFrench{Dans cet article nous démontrons que les groupes de cohomologie à support compact des champs de chtoucas sont des modules de type fini sur une alg\`ebre de Hecke. En guise d'application, nous étendons la construction des op\'erateurs d'excursion, d\'efinis dans \cite{vincent} sur l'espace des formes automorphes paraboliques, \`a l'espace des formes automorphes \`a support compact. Cela fournit la param\'etrisation de Langlands pour certains quotients de ce dernier espace, param\'etrisation compatible avec le morphisme «~terme constant~».}

\acknowledgement{This work has received funding from the European Research Council (ERC) under the European Union's Horizon 2020 research and innovation programme (grant agreement No 714405).}

\begin{document}

%%%%%%%%%%%%%%%%%%%%%%%%%%%%%%%
% Title page
%%%%%%%%%%%%%%%%%%%%%%%%%%%%%%%

%\removeabove{}
%\removebetween{}
%\removebelow{}

\maketitle

\begin{prelims}

\DisplayAbstractInEnglish

\bigskip

\DisplayKeyWords

\medskip

\DisplayMSCclass

\bigskip

\languagesection{Fran\c{c}ais}

\bigskip

\DisplayTitleInFrench

\medskip

\DisplayAbstractInFrench

\end{prelims}

%%%%%%%%%%%%%%%%%%%%%
% Table of Contents
%%%%%%%%%%%%%%%%%%%%%

\newpage

\setcounter{tocdepth}{2}

\tableofcontents

%%%%%%%%%%%%%%%%%%%%%
% Content begins here
%%%%%%%%%%%%%%%%%%%%%

\section*{Introduction}

Let $X$ be a smooth projective geometrically connected curve over a finite field $\Fq$. We denote by $F$ its function field, $\mb A$ the ring of adèles of $F$ and $\mb O$ the ring of integral adèles. 

Let $G$ be a connected split reductive group over $\Fq$.

\quad

Let $\Xi$ be a cocompact subgroup in $Z_G(F) \backslash Z_G(\mb A) $, where $Z_G$ is the center of $G$. Then the quotient $Z_G(F) \backslash Z_G(\mb A) / Z_G(\mb O) \Xi $ is finite. Let $N \subset X$ be a finite subscheme. We denote by $\mc O_N$ the ring of functions on $N$ and $K_N: = \Ker(G(\mb O) \rightarrow G(\mc O_N))$.

Let $\ell$ be a prime number not dividing $q$. Let $E$ be a finite extension of $\mathbb{Q}_{\ell}$ containing a square root of $q$. Let $\mc O_E$ be the ring of integers of $E$.
We denote by $C_c( G(F) \backslash G(\mb A) / K_N \Xi , E)$ the vector space of automophic forms with compact support.

Let $u$ be a place in $X \sm N$. Let $\mc O_u$ be the complete local ring at $u$ and let $F_u$ be its field of fractions. Let $\ms H_{G, u}:=C_{c}(G(\mc O_u)\backslash G(F_u)/G(\mc O_u), E)$ be the Hecke algebra of $G$ at the place $u$. 
The algebra $\ms H_{G, u}$ acts on $C_c( G(F) \backslash G(\mb A) / K_N \Xi , E)$ by convolution on the right. 
The following proposition was already known by some experts:
\begin{prop-intro}    \label{prop-space-of-auto-form-is-Hecke-mod-type-fini}
For any place $u$ of $X \sm N$, the vector space $C_c( G(F) \backslash G(\mb A) / K_N \Xi , \Ql)$ is a $\ms H_{G, u}$-module of finite type.
\end{prop-intro}

Due to the lack of reference, in Section 1 we recall the proof of Proposition \ref{prop-space-of-auto-form-is-Hecke-mod-type-fini} for $G=SL_2$ to illustrate the idea of the proof in the general case. 

\quad

Let $\wh G$ be the Langlands dual group of $G$ over $\Ql$ (i.e. the reductive group whose roots and weights are the coroots and coweights of $G$, and vice-versa). Let $I$ be a finite set and $W$ be a finite dimensional $\Ql$-linear representation of $\wh G^I$. Associated to these data, in \cite{var} Varshavsky defined the stack classifying $G$-shtukas (denoted by $\Cht_{G, N, I, W}$) over $(X \sm N)^I$ and
its $\ell$-adic cohomology groups with compact support in any degree $j \in \Z$ (denoted by $H_{G, N, I, W, \ov{x}}^{j}$, 
where $\ov{x}$ is a geometric point of $(X \sm N)^I$). 
In particular, when $I=\emptyset$ and $W = \bf 1$ (the one-dimensional trivial representation of the trivial group $\wh G^{\emptyset}$), the cohomology group $H_{G, N, \emptyset, {\bf 1}, \ov{x}}^0$ coincides with $C_c(G(F) \backslash G(\mathbb{A}) / K_N \Xi, \Ql)$. 

The cohomology group $H_{G, N, I, W, \ov{x}}^j$ is equipped with an action of $\ms H_{G, u}$ by Hecke correspondences on the stacks of shtukas. Suppose that $\ov{x}$ satisfies the condition in \ref{subsection-defi-geo-point-x-bar} (for example $\ov{x}$ can be a geometric generic point of $X^I$). 
One main result in this paper is
\begin{thm-intro}   \label{thm-coho-cht-is-Hecke-mod-type-fini-Hecke-en-v}
For any place $u$ of $X \sm N$, the cohomology group $H_{G, N, I, W, \ov{x}}^j$ is of finite type as $\ms H_{G, u}$-module.
\end{thm-intro}

Since $H_{G, N, \emptyset, {\bf 1}, \ov{x}}^0=C_c(G(F) \backslash G(\mathbb{A}) / K_N \Xi, \Ql)$, Theorem \ref{thm-coho-cht-is-Hecke-mod-type-fini-Hecke-en-v} is a generalization of Proposition \ref{prop-space-of-auto-form-is-Hecke-mod-type-fini}.
The proof of Theorem \ref{thm-coho-cht-is-Hecke-mod-type-fini-Hecke-en-v} is given in Section 2. The idea is similar to the case of automorphic forms explained in Section 1. In addition, we use the constant term morphisms of the cohomology groups of stacks of shtukas and the contractibility of deep enough Harder-Narasimhan strata established in \cite{cusp-coho}.

\quad

As an application, in Section 3, we extend the excursion operators, defined in \cite{vincent} on the vector space of cuspidal automorphic forms $C_c^{\on{cusp}}(G(F) \backslash G(\mathbb{A}) / K_N \Xi, \Ql)$, to $C_c(G(F) \backslash G(\mathbb{A}) / K_N \Xi, \Ql)$.

Concretely, fix an algebraic closure $\ov F$ of $F$. Let $\eta^I$ be the generic point of $X^I$ and fix a geometric point $\ov{\eta^I}$ over $\eta^I$. The cohomology group $H_{G, N, I, W, \ov{\eta^I}}^j$ is equipped with an action of $\pi_1(\eta^I, \ov{\eta^I})$ and an action of the partial Frobenius morphisms. 
In \cite{vincent}, V. Lafforgue 
\begin{enumerate}
\item defined a Hecke-finite part of $H_{G, N, I, W, \ov{\eta^I}}^j$ (denoted by $H_{G, N, I, W, \ov{\eta^I}}^{j, \, \on{Hf}}$). Then he used Drinfeld's lemma (for $\mc O_E$-modules of finite type) to construct an action of $\on{Gal}(\ov F / F)^I$ on $H_{G, N, I, W, \ov{\eta^I}}^{j, \, \on{Hf}}$.

\item with the help of (1), he constructed the excursion operators acting on $C_c^{\on{cusp}}(G(F) \backslash G(\mathbb{A}) / K_N \Xi, \Ql)$. 

\item with the help of (2), he obtained a canonical decomposition of $C_c^{\on{cusp}}(G(F) \backslash G(\mathbb{A}) / K_N \Xi, \Qlbar)$ indexed by Langlands parameters.
\end{enumerate}

\noindent In this paper, 

\begin{enumerate}[label=(\arabic*')]
\item thanks to Theorem \ref{thm-coho-cht-is-Hecke-mod-type-fini-Hecke-en-v}, we apply a variant of Drinfeld's lemma (for $\ms H_{G, u}$-modules of finite type) to $H_{G, N, I, W, \ov{\eta^I}}^j$ and obtain an action of $\on{Weil}(\ov F / F)^I$ on $H_{G, N, I, W, \ov{\eta^I}}^j$ (Proposition \ref{prop-FWeil-on-H-G-factors-through-Weil}).

\item with the help of (1'), 
we construct the excursion operators acting on $C_c(G(F) \backslash G(\mathbb{A}) / K_N \Xi, E)$ in Section \ref{subsection-excursion-operator}. When restricted to $C_c^{\on{cusp}}(G(F) \backslash G(\mathbb{A}) / K_N \Xi, \Ql)$, these excursion operators coincide with those defined in (2). 

\item Let $\ms I$ be a finite codimensional ideal of $\ms H_{G, u}$. 
With the help of (2'), 
we obtain a canonical decomposition of the quotient vector space $$C_c( G(F) \backslash G(\mathbb{A}) / K_N \Xi, \Qlbar) / \ms I \cdot C_c(G(F) \backslash G(\mathbb{A}) / K_N \Xi, \Qlbar)$$ indexed by 
Langlands parameters (Theorem \ref{thm-decomposition-of-C-c-quotient-I-by-param-Langlands}).
\end{enumerate}

More generally, in Section \ref{subsection-excursion-operator-on-cohomology}, we construct the excursion operators acting on the cohomology groups of stacks of shtukas, which are compatible with the action of the Hecke algebras and the action of the Galois group. Moreover, there is a canonical decomposition similar to (3') for some quotient vector spaces of the cohomology groups.

\quad

In Section 4, we show that the action of the excursion operators commutes with the constant term morphisms. As a consequence, the parametrization in (3') is compatible with the constant term morphisms.

\subsection*{Acknowledgments}
I would like to thank Dennis Gaitsgory, Vincent Lafforgue, G\'erard Laumon, Jack Thorne and Zhiyou Wu for stimulating discussions.

\section{Example of automorphic forms for $SL_2$}

In this section, we consider the cohomology group of stacks of shtukas without paws, i.e. the space of automorphic forms.
In this section, let $N=\emptyset$ and $G=SL_2$.

\sssec{}
Let $\Bun_{G}$ be the classifying stack of $G$-bundles over $X$. For $G=SL_2$, it is the classifying stack of rank $2$ vector bundles on $X$ with trivial determinant. It is well-known that (see \cite[rem.~8.21]{vincent}  for more details) $$G(F) \backslash G(\mb A) / G(\mb O) = \Bun_{G}(\Fq) .$$ Thus $$C_c( G(F) \backslash G(\mb A) / G(\mb O) , \Ql) = C_c( \Bun_{G}(\Fq)  , \Ql).$$ It is also well-known that the action of $\ms H_{G, u}$ by convolution on $C_c( G(F) \backslash G(\mb A) / G(\mb O) , \Ql) $ coincides with the action of $\ms H_{G, u}$ by Hecke correspondences on $C_c( \Bun_{G}(\Fq)  , \Ql)$.

\sssec{}
The stack $\Bun_{G}$ has a Harder-Narasimhan stratification. For $G=SL_2$, we have
$$\Bun_{G} = \varinjlim _{d \in \Z_{\geq 0}} \Bun_{G}^{\leq (d, -d)} ,$$ 
where $\Bun_G^{\leq (d, -d)}$ is the open substack of $\Bun_{G}$ defined by the condition that the degree of any line subbundle of the rank $2$ vector bundle is $\leq d$. 
For any $d_1, d_2 \in \Z_{\geq 0}$ such that $d_1 \leq d_2$, we have an open immersion $\Bun_G^{\leq (d_1, -d_1)}  \hookrightarrow \Bun_G^{\leq (d_2, -d_2)}$. It induces an inclusion $C_c(\Bun_G^{\leq (d_1, -d_1)}(\Fq) , \Ql) \subset C_c(\Bun_G^{\leq (d_2, -d_2)}(\Fq) , \Ql)$ by extension by zero.
Thus the vector space $C_c(\Bun_G(\Fq) , \Ql)$ has a filtration $F_G^\bullet$ indexed by $d \in \Z_{\geq 0}$:
$$F_G^{d} := C_c(\Bun_G^{\leq (d, -d)}(\Fq) , \Ql) ,$$
with associated graded quotients (for $d \geq 1$)
$$gr_G^d := F_G^d / F_G^{d-1} = C_c(\Bun_G^{= (d, -d)}(\Fq) , \Ql),$$
where $\Bun_G^{= (d, -d)}$ is the closed substack of $\Bun_G^{\leq (d, -d)}$ defined by the condition that for the canonical Harder-Narasimhan filtration $0 \subset \mc L \subset \mc E$ of the rank $2$ vector bundle $\mc E$, we have $\on{deg} \mc L = d$. (When $d=0$, we have $\Bun_G^{= (0, 0)} = \Bun_G^{\leq (0, 0)}$ which classifies the semistable vector bundles.)

\quad

Since every $F_G^{d}$ has finite dimension, Proposition \ref{prop-space-of-auto-form-is-Hecke-mod-type-fini} is a direct consequence of Lemma \ref{lem-function-case-C-c-generated-by-d-0} below.

\sssec{}
Let $l := \on{deg}(u)$. Let $g$ be the genus of $X$. Let $d_0 = \on{max} \{0, g-1 \} + l$.

\begin{lem}    \label{lem-function-case-C-c-generated-by-d-0}
We have
$$C_c(\Bun_G(\Fq)  , \Ql) = \ms H_{G, u} \cdot F_G^{d_0}.$$ 
\end{lem}
\dem
Applying successively Lemma \ref{lem-function-case-F-d+1-generated-by-F-d} below to $d_0+1$, $d_0+2$, $d_0+3$, $\cdots$, we deduce that for any $d > d_0$, we have $F_G^d \subset \ms H_{G, u} \cdot F_G^{d_0}.$ This implies Lemma \ref{lem-function-case-C-c-generated-by-d-0}.
\cqfd

\quad

\sssec{}
Fix a Borel subgroup $B$ of $G=SL_2$. The group $G$ has only one simple coroot $\check{\alpha}=(1, -1)$. 
A generator of the coweight lattice is $\omega=\check{\alpha}$. 
Let $\varpi$ be a uniformizer of $\mc O_u$. Let $h^G_{\omega}$ be the Hecke operator in $\ms H_{G, u}$ associated to the characteristic function of $G(\mc O_u) \varpi^{\omega} G(\mc O_u)$, where
$$\varpi^{\omega} = \left( \begin{array}{cc} \varpi &  \\
           & \varpi^{-1} \end{array} \right) .$$

\sssec{}
The action of $h^G_{\omega}$ on $F_G^\bullet$ is defined in the following way:
let $\Gamma( G(\mc O_u) \varpi^{\omega} G(\mc O_u))$ be the groupoid classifying pairs $(\mc G \dashrightarrow \mc G')$, where $\mc G$, $\mc G'$ are rank $2$ vector bundles of trivial determinant, $\mc G \dashrightarrow \mc G'$ is an isomorphism outside $u$ such that when restricted to the formal disc on $u$, the relative position of $\mc G$ and $\mc G'$ is equal to $\varpi^{\omega}$. 
We have a Hecke correspondence (where the arrows are morphisms of groupoids):
$$ \Bun_{G}(\Fq)  \xleftarrow{\pr_1}    \Gamma( G(\mc O_u) \varpi^{\omega} G(\mc O_u))   \xrightarrow{\pr_2}   \Bun_{G}(\Fq).$$
$$\mc G \mapsfrom (\mc G \dashrightarrow \mc G') \mapsto \mc G' $$

For any $d \in \Z_{\geq 0}$ such that $d - l > 0$, the projection $\pr_1$ sends $\pr_2^{-1}( \Bun_{G}^{\leq (d-l, -(d-l))} ) $ to $\Bun_{G}^{\leq (d, -d)}  $.
We deduce a morphism
$$h^G_{\omega}: F_{G}^{d-l} \rightarrow F_{G}^{d}, \quad f \mapsto (\pr_1)_!(\pr_2)^* f .$$ 
Similarly, we have a morphism $h^G_{\omega}: F_{G}^{d-l-1} \rightarrow F_{G}^{d-1}$.
We deduce a morphism 
\begin{equation}   \label{equation-h-G-omega-gr-d-l-to-gr-d}
h^G_{\omega}: gr_G^{d-l} \rightarrow gr_G^{d}.
\end{equation}

\begin{lem} \label{lem-operatpr-on-filtration-GL-2}
If $d > d_0$, then the morphism (\ref{equation-h-G-omega-gr-d-l-to-gr-d}) is an isomorphism.
\end{lem}

The proof of Lemma \ref{lem-operatpr-on-filtration-GL-2} will be given in the remaining part of this section.

\begin{lem} \label{lem-function-case-F-d+1-generated-by-F-d}
For any $d \in \Z_{\geq 0}$ such that $d > d_0$, we have
$$F_G^{d} = h^G_{\omega} \cdot F_G^{d-l} + F_G^{d-1} \subset \ms H_{G, u} \cdot F_G^{d-1}.$$
\end{lem}
\dem
This is a direct consequence of Lemma \ref{lem-operatpr-on-filtration-GL-2}.
\cqfd

\sssec{}   
To prove Lemma \ref{lem-operatpr-on-filtration-GL-2}, we need some preparation. Recall that $B$ is a Borel subgroup of $G$.
Let $\Bun_B$ be the classifying stack of $B$-bundles on $X$. For $G=SL_2$, it is the classifying stack of pairs $(\mc L \subset \mc E)$, where $\mc E$ is in $\Bun_G$ and $\mc L$ is a line subbundle.
Let $T$ be the torus of $G$. Let $\Bun_T$ be the classifying stack of $T$-bundles on $X$. For $G=SL_2$, it is the classifying stack of pairs of line bundles $(\mc L, \mc L^{-1})$ on $X$. 
We have a correspondence 
\begin{equation}    \label{equation-Bun-G-Bun-B-Bun-T}
\Bun_G \xleftarrow{i} \Bun_B \xrightarrow{\pi} \Bun_T, \quad \mc E \mapsfrom (\mc L \subset \mc E) \mapsto (\mc L, \mc E / \mc L).
\end{equation}

\sssec{}
For any $d' \in \Z$, let $\Bun_T^{=(d', -d')}$ be the open and closed substack of $\Bun_T$ defined by the condition that $\on{deg}(\mc L) = d'$. (Note that $d'$ can be negative.) Then $\Bun_T = \bigsqcup_{d' \in \Z} \Bun_T^{=(d', -d')}$.
Let $H_T^{=d'}:=C_c(\Bun_T^{=(d', -d')}(\Fq) , \Ql)$. 
Consider the vector space $\varinjlim _{d \in \Z_{\geq 0}}  \prod_{d' \in \Z, \, d' \leq d } H_T^{=d'} $. It has a filtration $F_T^\bullet$ indexed by $d \in \Z_{\geq 0}$:
$$F_T^{d}:= \prod_{d' \in \Z, \, d' \leq d} H_T^{=d'} , $$
with associated graded quotients
$$ gr_T^{d}:= F_T^{d} / F_T^{d-1} = H_T^{=d} .$$ 

\sssec{}
For any $d\in \Z_{\geq 0}$, we define $\Bun_B^{\leq (d, -d)} := i^{-1}(\Bun_G^{\leq (d, -d)})$ and $\Bun_T^{\leq (d, -d)} := \sqcup_{d' \leq d} \Bun_T^{={(d', -d')}}.$
Then (\ref{equation-Bun-G-Bun-B-Bun-T}) induces a morphism of groupoids 
$$\Bun_G^{\leq (d, -d)}(\Fq) \xleftarrow{i} \Bun_B^{\leq (d, -d)} (\Fq)  \xrightarrow{\pi} \Bun_T^{\leq (d, -d)}(\Fq).$$
We have a constant term morphism 
$$C_G^B: F_G^d \rightarrow F_T^d, \quad f \mapsto \pi_! i^* f.$$ 
When $d \geq 1$, we have a commutative diagram:
$$
\xymatrix{
F_G^{d-1} \ar[r]   \ar[d]^{C_G^B}
& F_G^d   \ar[r]   \ar[d]^{C_G^B}
& gr_G^d  \ar[d]^{C_G^B} \\
F_T^{d-1}  \ar[r]   
& F_T^{d}  \ar[r] 
& gr_T^d
}
$$

\sssec{}
Let $\ms H_{T, u}:=C_{c}(T(\mc O_u)\backslash T(F_u)/T(\mc O_u), \Ql)$ be the Hecke algebra of $T$ at the place $u$. It acts on $F_T^\bullet$ by Hecke correspondences. 
Let $h^T_{\omega} \in \ms H_{T, u}$ be the Hecke operator associated to the characteristic function of $T(\mc O_u) \varpi^{\omega} T(\mc O_u)$.
By the Satake isomorphism, we have $h^G_{\omega} \in \ms H_{G, u} \isom \ms H_{T, u}^{\mf S_2} \subset \ms H_{T, u}$, where $ \mf S_2$ is the symmetric group. 
In $\ms H_{T, u}$, we have $h^G_{\omega} = q ( (h^{T}_{\omega})^{ -1} + h^T_{\omega} ) $.

\begin{lem}   \label{lem-function-case-h-G-act-on-F-T}
\leavevmode
\begin{enumerate}[label=$(\roman*)$]
\item For any $d \in \Z$, the action of $h^G_{\omega}$ induces a morphism $h^G_{\omega}: F_T^d \rightarrow F_T^{d+l}$;
\item the induced morphism on the associated graded quotients $ h^G_{\omega}: gr_T^d \rightarrow gr_T^{d+l}$ is an isomorphism.
\end{enumerate}
\end{lem}
\dem
(i) For any $d \in \Z$, the action of $h^T_{\omega}$
induces an isomorphism:
$$\Bun_T^{=(d, -d)} \isom \Bun_T^{=(d + l, -(d+l))} , \quad (\mc L_1, \mc L_2) \mapsto (\mc L_1(u), \mc L_2(-u))$$ 
(the inverse is induced by $(h^{T}_{\omega})^{ -1}$).
Thus the Hecke operator $h^T_{\omega}$ induces an isomorphism
\begin{equation}   \label{equation-SL2-case-H-T-d-to-H-T-d+1}
h^T_{\omega}: H_T^{=d}  \isom H_T^{=d+l} .
\end{equation}
The Hecke operator $(h^{T}_{\omega})^{ -1}$ induces 
\begin{equation}   \label{equation-SL2-case-H-T-d-to-H-T-d-1}
(h^{T}_{\omega})^{ -1}: H_T^{=d}  \isom H_T^{=d-l}. 
\end{equation}
As a consequence, the Hecke operator $h^G_{\omega} = q ( (h^{T}_{\omega})^{ -1} + h^T_{\omega})$ induces a morphism
\begin{equation}   \label{equation-SL2-case-H-T-d-to-H-T-d-1-d+1}
h^G_{\omega} : H_T^{=d}  \rightarrow H_T^{=d-l} \times H_T^{=d+l} .
\end{equation}
Applying (\ref{equation-SL2-case-H-T-d-to-H-T-d-1-d+1}) to all $d' \leq d$ and taking the product, we get a morphism: $h^G_{\omega}: \prod_{d' \leq d} H_T^{=d'} \rightarrow \prod_{d'' \leq d+l} H_T^{=d''}$, i.e. $h^G_{\omega}: F_T^d \rightarrow F_T^{d+l}$.

\quad

(ii) The morphisms $h^G_{\omega}: F_T^d \rightarrow F_T^{d+l}$ and $h^G_{\omega}: F_T^{d-1} \rightarrow F_T^{d-1+l}$ induce a morphism
$$h^G_{\omega}: F_T^d  / F_T^{d-1} \rightarrow F_T^{d+l} / F_T^{d+l-1}.$$

Recall that $gr_T^d = F_T^d  / F_T^{d-1} = H_T^{=d}$. 
The action of $h^G_{\omega} = q ( h^T_{\omega} + (h^{T}_{\omega})^{ -1})$ on $F_T^d  / F_T^{d-1}$ is equal to the action of $q h^T_{\omega}$ on $F_T^d  / F_T^{d-1}$, because by (\ref{equation-SL2-case-H-T-d-to-H-T-d-1}) we have $$(h^T_{\omega})^{-1}(  H_T^{=d}  ) = H_T^{=d-l}  \subset F_T^{d+l-1},$$ i.e. the morphism $(h^{T}_{\omega})^{ -1} : F_T^d  / F_T^{d-1}  \rightarrow F_T^{d+l} / F_T^{d+l-1}$ is the zero morphism. 
\cqfd          
           
\quad         

\noindent {\bf Proof of Lemma \ref{lem-operatpr-on-filtration-GL-2}.} 
Let $d>l$. Since the constant term morphism commutes with the action of the Hecke algebra $\ms H_{G, u}$, we have a commutative diagram
\begin{equation}   \label{diagram-example-SL-2-gr-G-gr-T-commute}
\xymatrix{
gr_G^{d-l} \ar[d]^{C_G^B}  \ar[r]^{h^G_{\omega}}
& gr_G^{d}  \ar[d]^{C_G^B}  \\
gr_T^{d-l} \ar[r]^{h^G_{\omega}}
& gr_T^{d}
}
\end{equation}

It is well-known that if $d>\on{max} \{0, g-1 \}$, then $C_G^B: gr_G^d \rightarrow gr_T^d$ is an isomorphism. (Indeed, for any $\mc E \in \Bun_G^{=(d, -d)} (\Fq)$, let $0 \subset \mc L \subset \mc E$ be its canonical Harder-Narasimhan filtration, then $\on{deg} \mc L = d$. We have $\mc E / \mc L = \mc L^{-1}$. The condition $d>\on{max} \{0, g-1 \}$ implies that $\on{Ext}^1(\mc L^{-1}, \mc L) \simeq H^1(X, \mc L^2)=0$. Thus $\mc E \simeq \mc L \oplus \mc L^{-1}$. Let $\Bun_B^{=(d, -d)} := i^{-1}(\Bun_G^{=(d, -d)} ) \cap \pi^{-1} (\Bun_T^{=(d, -d)} )$. 
We deduce that
\begin{equation}  \label{equation-Bun-G-Bun-B-Bun-T-strata-bij}
\Bun_G^{=(d, -d)} (\Fq)\leftarrow \Bun_B^{=(d, -d)}(\Fq) \rightarrow \Bun_T^{=(d, -d)}(\Fq)
\end{equation}
are bijective. Thus $C_c(\Bun_G^{= (d, -d)}(\Fq) , \Ql) \isom C_c(\Bun_T^{=(d, -d)}(\Fq) , \Ql)$.) 

As a consequence, if $d-l >\on{max} \{0, g-1 \}$, then the two vertical morphisms in (\ref{diagram-example-SL-2-gr-G-gr-T-commute}) are isomorphisms.
Moreover, by Lemma \ref{lem-function-case-h-G-act-on-F-T}, the lower line in (\ref{diagram-example-SL-2-gr-G-gr-T-commute}) is an isomorphism. So the upper line in (\ref{diagram-example-SL-2-gr-G-gr-T-commute}) is an isomorphism.
\cqfd

\begin{rem}
For a general group $G$, to prove Proposition \ref{prop-space-of-auto-form-is-Hecke-mod-type-fini}, we need \cite[Proposition~9.2.2]{DG15} to show that an analogue of (\ref{equation-Bun-G-Bun-B-Bun-T-strata-bij}) is bijective. We do not give details here because we will prove some more general statements in Section 2.
\end{rem}

\quad

\section{Proof of Theorem \ref{thm-coho-cht-is-Hecke-mod-type-fini-Hecke-en-v}}

In this section, we consider the general case.

\subsection{Reminder of cohomologies of stacks of shtukas}

The stacks of shtukas and their cohomologies are defined in \cite{var}, recalled in \cite[Section~2]{vincent} and in \cite[Section~1 and~2]{cusp-coho}.

As in the introduction, let $N \subset X$ be a finite subscheme. Let $I$ be a finite set and $W$ be a representation of $\wh G^I$ (in this paper, this always means a finite dimensional $\Ql$-linear representation of $\wh G^I$). 

\sssec{}   \label{subsection-Cht-G-N-I-W}
In \cite[Definition~2.3.1]{cusp-coho}, 
we defined a classifying stack of $G$-shtukas $\Cht_{G, N, I, W}$ over $(X \sm N)^I$ (this stack is denoted by $\Cht_{N, I, W}^{(I)}$ in \cite[notation~4.4]{vincent}). It is a Deligne-Mumford stack locally of finite type.  

\sssec{} \label{subsection-Lambda-G-ad-+-Q}
Let $G^{\mr{ad}}$ be the adjoint group of $G$. Let $\wh\Lambda_{G^{\mr{ad}}}$ denote the coweight lattice of $G^{\mr{ad}}$ and $\wh \Lambda_{G^{\mr{ad}}}^{\Q} : = \wh\Lambda_{G^{\mr{ad}}} \otimes_{\Z} \Q$. Fix a Borel subgroup $B \subset G$. Let $\wh\Lambda_{G^{\mr{ad}}}^{+} \subset \wh\Lambda_{G^{\mr{ad}}}$ denote the monoid of dominant coweights. Let $\wh\Lambda_{G^{\mr{ad}}}^{+, \Q} \subset \wh \Lambda_{G^{\mr{ad}}}^{\Q} $ be the corresponding rational cone. 

For any $\mu_1, \mu_2 \in \wh \Lambda_{G^{\mr{ad}}}^{+, \Q}$, we say $\mu_1 \leq \mu_2$ if and only if $\mu_2 - \mu_1$ is a linear combination with coefficients in $\Q_{\geq 0}$ of simple coroots of $G^{\mr{ad}}$, or equivalently, of simple coroots of $G$ modulo $\wh \Lambda_{Z_G}^{\Q} $, where $\wh \Lambda_{Z_G}$ is the coweight lattice of $Z_G$ and $\wh{\Lambda}_{Z_G}^{\Q} := \wh{\Lambda}_{Z_G} \otimes_{\Z} \Q$.

\sssec{}
In \cite[d\'efi.~2.1]{vincent} and \cite[\S~2.3.3]{cusp-coho}, for any $\mu \in \wh \Lambda_{G^{\mr{ad}}}^{+, \Q}$,  we defined an open substack $\Cht_{G, N, I, W}^{\leq \mu}$
of $\Cht_{G, N, I, W}$.
For any $\mu_1, \mu_2 \in \wh \Lambda_{G^{\mr{ad}}}^{+, \Q}$ and $\mu_1 \leq \mu_2$, we have an open immersion: 
\begin{equation}   \label{morphism-open-immersion-Cht}
\Cht_{G, N, I, W}^{\leq \mu_1}  \hookrightarrow \Cht_{G, N, I, W}^{\leq \mu_2}  .
\end{equation}
We have $\Cht_{G, N, I, W}  = \bigcup_{ \mu \in \wh \Lambda_{G^{\mr{ad}}}^{+, \Q} } \Cht_{G, N, I, W}^{\leq \mu} $. 

\begin{rem}
The above Harder-Narasimhan truncations $\Cht_{G, N, I, W}^{\leq \mu}$ on $\Cht_{G, N, I, W}$ come from the Harder-Narasimhan truncations defined in \cite{schieder} and \cite{DG15} for $\Bun_G$.
\end{rem}

\sssec{}
Recall that we have fixed $\Xi \subset Z_G(F) \backslash Z_G(\mb A) $ in the introduction. As recalled in \cite[Proposition~2.3.4]{cusp-coho}, the quotient $\Cht_{G, N, I, W}^{\leq \mu} / \Xi$ is a Deligne-Mumford stack of finite type. 

\sssec{}
In \cite[d\'efi.~4.5]{vincent} and \cite[Definition~2.4.7]{cusp-coho}, we defined a Satake perverse sheaf $\mc F_{G, N, I, W}^{\Xi}$ over $\Cht_{G, N, I, W} / \Xi$ (with the perverse normalization relative to $(X \sm N)^I$). When $W$ is irreducible, $\mc F_{G, N, I, W}^{\Xi}$ is (non canonically) isomorphic to the intersection complex of $\Cht_{G, N, I, W} / \Xi$ (with coefficients in $\Ql$ and with the perverse normalization relative to $(X \sm N)^I$).

\sssec{}
As in \cite[\S~1.1.7]{cusp-coho}, let $\mf p_G: \Cht_{G, N, I, W} / \Xi \rightarrow (X \sm N)^I$ be the projection of paws.
Let $\ov{x} \rightarrow (X \sm N)^I$ be a geometric point. 

\begin{defi} [\emph{cf.} d\'efinition~4.7 in \cite{vincent} and Definition~2.5.1 in \cite{cusp-coho}] \label{def-H-I-W-leq-mu-sheaf}
For any $\mu \in \wh \Lambda_{G^{\mr{ad}}}^{+, \Q}$ and any $j \in \Z$, we define
$$\mc H_{G, N, I, W}^{j, \, \leq\mu} := R^j (\mf p_G)_! ( \restr{ \mc F_{G, N, I, W}^{\Xi} }{ \Cht_{G, N, I, W}^{\leq \mu} / \Xi } );$$
$$H_{G, N, I, W, \ov{x}}^{j, \, \leq\mu} := \restr{ \mc H_{G, N, I, W}^{j, \, \leq\mu} }{ \ov{x} } .$$
\end{defi}

\begin{rem} 
(i) Since $\Cht_{G, N, I, W}^{\leq \mu} / \Xi$ is of finite type, $\mc H_{G, N, I, W}^{j, \, \leq\mu} $ is an $\Ql$-constructible sheaf on $(X \sm N)^I$ and $H_{G, N, I, W, \ov{x}}^{j, \, \leq\mu}$ is a finite dimensional $\Ql$-vector space.  

(ii) $\mc H_{G, N, I, W}^{j, \, \leq\mu} $ and $H_{G, N, I, W, \ov{x}}^{j, \, \leq\mu}$ depend on $\Xi$. We do not write $\Xi$ in the index to simplify the notations.

(iii) Let $d$ be the dimension of $\Cht_{G, N, I, W}$ relatively to $(X \sm N)^I$. By \cite[\S~4.2]{bbd}, the cohomology group $H_{G, N, I, W, \ov{x}}^{j, \, \leq\mu} $ is concentrated in degrees $j \in [-d, d]$. 
\end{rem}

\sssec{}
For any $\mu_1, \mu_2 \in \wh \Lambda_{G^{\mr{ad}}}^{+, \Q}$ and $\mu_1 \leq \mu_2$, the open immersion (\ref{morphism-open-immersion-Cht}) induces a morphism of sheaves:
$$\mc H_{G, N, I, W}^{j, \, \leq\mu_1}  \rightarrow \mc H_{G, N, I, W}^{j, \, \leq\mu_2} .$$

\begin{defi} [\emph{cf.} Definition~2.5.3 in \cite{cusp-coho}]   \label{def-H-G-I-W-sheaf}
We define
$$\mc H_{G, N, I, W}^j:= \varinjlim _{\mu}   \mc H_{G, N, I, W}^{j, \, \leq\mu} ,$$ in the category of inductive limits of constructible sheaves over $(X \sm N)^I$. And we define
$$H_{G, N, I, W, \ov{x}}^j:= \varinjlim _{\mu}   H_{G, N, I, W, \ov{x}}^{j, \, \leq\mu} .$$
We have $H_{G, N, I, W, \ov{x}}^j = \restr{ \mc H_{G, N, I, W}^j }{ \ov{x} }$.
\end{defi}

\sssec{}
By \cite[\S~2.1.9]{cusp-coho}, the cohomology $\mc H_{G, N, I, W}^j$ is functorial on $W$. In particular, for $W=W_1 \oplus W_2$, we have
$\Cht_{G, N, I, W} = \Cht_{G, N, I, W_1} \cup \Cht_{G, N, I, W_2}$ and $\mc H_{G, N, I, W}^j = \mc H_{G, N, I, W_1}^j \oplus \mc H_{G, N, I, W_2}^j$.

\sssec{}
When $I=\emptyset$ and $W = \bf 1$, we have $\Cht_{G, N, \emptyset, {\bf 1}} = \Bun_{G, N}(\Fq) =G(F) \backslash G(\mathbb{A}) / K_N $ and, for $\ov{x}$ a geometric point of $\on{Spec} \Fq$, we also have $H_{G, N, \emptyset, {\bf 1}, \ov{x}}^0=C_c(G(F) \backslash G(\mathbb{A}) / K_N \Xi, \Ql)$. 

\sssec{}   \label{subsection-Hecke-G-action}
Let $u$ be a place in $X \sm N$. Let $l = \on{deg}(u)$. The Hecke algebra $\ms H_{G, u}$ acts on $H_{G, N, I, W, \ov{x}}^j$ by Hecke correspondences (see \cite[\S~2.20 and~4.4]{vincent} and \cite[\S~6.1-6.2]{cusp-coho}). This action does not preserve the subspaces $H_{G, N, I, W, \ov{x}}^{j, \, \leq\mu}$.

Concretely, for any dominant coweight $\beta$ of $G$, let $h^G_{\beta} \in \ms H_{G, u} $ be the characteristic function of $G(\mc O_u) \varpi^{\beta} G(\mc O_u)$. We have the Hecke correspondence associated to $h^G_{\beta}$: 
$$ \Cht_{G, N, I, W}/\Xi   \xleftarrow{\pr_1}    \Gamma( G(\mc O_u)\varpi^{ \beta }G(\mc O_u))   \xrightarrow{\pr_2}  \Cht_{G, N, I, W} / \Xi ,$$
where $\pr_1$ and $\pr_2$ are finite étale morphisms.
For any $\mu \in \frac{1}{r} \wh R_{G^{\mr{ad}}}^+$,
the projection $\pr_1$ sends $\pr_2^{-1}(  \Cht_{G, N, I, W}^{\leq \mu} / \Xi  )$ to $\Cht_{G, N, I, W}^{\leq \mu + \beta l }/\Xi $.
It induces a morphism
$$h^G_{\beta}:  H_{G, N, I, W, \ov{x}}^{j, \, \leq \mu  } \rightarrow H_{G, N, I, W, \ov{x}}^{j, \, \leq \mu + \beta l}.$$

\sssec{}   \label{subsection-defi-geo-point-x-bar}
For $i \in I$, let $\on{pr}_i: X^I \rightarrow X$ be the projection to the $i$-th factor.
From now on, let $\ov{x}$ be a geometric point of $(X \sm N)^I$ such that for every $i, j \in I$, the image of the composition $$\ov{x} \rightarrow (X \sm N)^I \xrightarrow{(\on{pr}_i, \on{pr}_j) } (X \sm N) \times (X \sm N)$$ 
is not included in the graph of any non-zero power of Frobenius morphism $\on{Frob}: X \sm N \rightarrow X \sm N$. 
(This condition will be needed in \ref{subsection-defi-CT-morphism} and Proposition \ref{prop-TC-suite-exacte}.)

In particular, when $i=j \in I$, the above condition is equivalent to the condition that the composition $\ov{x} \rightarrow (X \sm N)^I \xrightarrow{\on{pr}_i } X \sm N $ is over the generic point $\eta$ of $X \sm N$.

One example of geometric point satisfying the above condition is $\ov{\eta^I}$, a geometric point over the generic point $\eta^I$ of $X^I$. 
Another example of geometric point satisfying the above condition is $\Delta(\ov{\eta})$, where $\Delta: X \hookrightarrow X^I$ is the diagonal inclusion and $\ov{\eta}$ is a geometric point over $\eta$.

However, for any $i \in I$, let $\Frob_{\{i\}}: X^I \rightarrow X^I$ be the morphism sending $(x_j)_{j \in I}$ to $(x_j')_{j \in I}$, with $x_i' = \Frob(x_i)$ and $x_j' = x_j$ if $j \neq i$. Then $\Frob_{\{i\}} \Delta(\ov{\eta})$ does not satisfy the above condition.

\begin{nota}
In the remaining part of Section 2, we write $H_G^{j, \, \leq \mu}$ $($resp. $H_G^{j})$ instead of $H_{G, N, I, W, \ov{x}}^{j, \, \leq\mu}$ $($resp. $H_{G, N, I, W, \ov{x}}^{j})$ to simplify the notation.
\end{nota}

\subsection{Strategy of the proof of Theorem \ref{thm-coho-cht-is-Hecke-mod-type-fini-Hecke-en-v}}

\sssec{}
We denote by $\wh R_{G^{\mr{ad}}}^+$ the dominant cone of the coroot lattice of $G^{\mr{ad}}$. We have $\wh R_{G^{\mr{ad}}}^+ \subset \wh{\Lambda}_{G^{\mr{ad}}}^{+}$.
For any positive integer $r$, we have $\wh{\Lambda}_{G^{\mr{ad}}}^{+}  \subset \frac{1}{r} \wh R_{G^{\mr{ad}}}^+ \subset \wh \Lambda_{G^{\mr{ad}}}^{+, \Q}$. 
Fix $r$ large enough as in \cite[\S~5.1.1]{cusp-coho}.

\begin{defi}
We define a filtration $F_G^\bullet$ of $H_{G}^j$:
$$\mu \in \frac{1}{r} \wh R_{G^{\mr{ad}}}^+, \quad F_G^{\mu}:=\on{Im} ( H_{G}^{j, \, \leq \mu} \rightarrow H_G^j ) .$$
\end{defi}

\sssec{}  \label{subsection-F-G-mu-1-includ-F-G-mu-2}
For any $\mu_1, \mu_2 \in  \frac{1}{r} \wh R_{G^{\mr{ad}}}^+$ and $\mu_1 \leq \mu_2$, the morphism $H_G^{j, \, \leq \mu_1 } \rightarrow H_G^{j}$ factors through $H_G^{j, \, \leq \mu_2}$. Thus $F_G^{\mu_1} \subset F_G^{\mu_2}$.

\quad

Since every $F_G^{\mu}$ has finite dimension, Theorem \ref{thm-coho-cht-is-Hecke-mod-type-fini-Hecke-en-v} is a direct consequence of the following proposition:
\begin{prop}   \label{prop-H-G-is-generated-by-Hecke-alg}
There exists $\mu_0 \in \frac{1}{r} \wh R_{G^{\mr{ad}}}^+$ large enough $($a priori depending on $u$ and $\ov{x})$, such that 
$$H_G^j = \ms H_{G, u} \cdot F_G^{\mu_0} $$
\end{prop}
Proposition \ref{prop-H-G-is-generated-by-Hecke-alg} will follow from Lemma \ref{lem-coho-case-F-mu+alpha-generated-by-F-mu} below.

\quad

\sssec{}
We denote by $\Gamma_G$ the set of vertices of the Dynkin diagram of $G$. For any $i \in \Gamma_G$, we denote by $\alpha_i$ (resp. $\check{\alpha}_i$) the simple root (resp. simple coroot) of $G$ associated to $i$. We denote by $\langle \ , \ \rangle$ the natural pairing between coweights and weights. 
Let $P_i$ be the standard maximal parabolic with Levi quotient $M_i$ such that $\Gamma_G - \Gamma_{M_i} = \{ i \}$, where $\Gamma_{M_i}$ is the set of vertices of the Dynkin diagram of $M_i$.  
By a quasi-fundamental coweight of $G$ we mean the smallest positive multiple of a fundamental coweight of $G^{\mr{ad}}$, which is a coweight of $G$.
Let $\omega_i$ be the quasi-fundamental coweight of $G$ such that $\langle \omega_i, \alpha_j \rangle =0$ for any $j \in \Gamma_{M_i}$.

\begin{notation}     \label{notation-constant-C-G-I}
Let $C_G \in \Q_{\leq 0}$ be the constant in Proposition 5.1.5 (c) of \cite{cusp-coho}. Let 
$$
C(G, u):= C_G+ \on{max}_{i \in \Gamma_G} \langle \omega_i l + \frac{1}{r}\check{\alpha}_i, \alpha_i \rangle ,
$$
where $l = \on{deg} (u)$.
\end{notation}

\begin{rem}
The above definition implies that for any $\mu \in \frac{1}{r} \wh R_{G^{\mr{ad}}}^+$ and any $i \in \Gamma_G$, if $\langle \mu, \alpha_i \rangle \geq C(G, u)$, then $\langle \mu - \omega_i l -\frac{1}{r}\check{\alpha}_i , \alpha_j \rangle \geq 0$ for all $j \in \Gamma_G$, i.e. $\mu-\omega_i l - \frac{1}{r}\check{\alpha}_i$ is dominant. 
\end{rem}

\begin{lem}   \label{lem-coho-case-F-mu+alpha-generated-by-F-mu}
Let $\mu \in \frac{1}{r} \wh R_{G^{\mr{ad}}}^+$ such that $\langle \mu, \alpha_j \rangle > C(G, u)$ for all $j \in \Gamma_G$, then for  
any $i \in \Gamma_G$,
we have
$$F_G^{\mu } \subset \ms H_{G, u} \cdot F_G^{\mu-\frac{1}{r}\check{\alpha}_i }.$$
\end{lem}
The proof of Lemma \ref{lem-coho-case-F-mu+alpha-generated-by-F-mu} will be given after Lemma \ref{lem-Hecke-operator-on-filtration-G} below.

\quad

\noindent {\bf Proof of Proposition \ref{prop-H-G-is-generated-by-Hecke-alg} admitting Lemma \ref{lem-coho-case-F-mu+alpha-generated-by-F-mu}.} 
We use the same argument as in the proof of Lemma 5.3.6 in \cite{cusp-coho}.
Let $Z(C)$ be the set of $\mu \in  \frac{1}{r} \wh R_{G^{\mr{ad}}}^+ $ such that $\langle \mu , \alpha_j \rangle > C(G, u)$ for all $j \in \Gamma_G$.  
Let $\Omega(C)$ be the set of $\mu \in Z(C)$ such that $\mu - \frac{1}{r} \check{\alpha}_i \notin Z(C)$ for all $i \in \Gamma_G$.
The set $\Omega(C)$ is bounded, thus is finite.
Let $\mu_0 \in  \frac{1}{r} \wh R_{G^{\mr{ad}}}^+ $ such that $\mu_0 > \mu$ for all $\mu \in \Omega(C)$. 

For any $\lambda \in Z(C)$, there exists a (zigzag) chain $\lambda = \lambda^{(0)} > \lambda^{(1)} > \cdots > \lambda^{(m-1)}  > \lambda^{(m)} $ in $\frac{1}{r} \wh R_{G^{\mr{ad}}}^+$ for some $m \in \Z_{\geq 0}$ such that 
\begin{enumerate}[label=$(\roman*)$]
\item for all $j$, we have $\lambda^{(j)} \in Z(C)$,
\item for all $j$, we have $\lambda^{(j)} - \lambda^{(j+1)} = \frac{1}{r} \check{\alpha}_i$ for some $i \in \Gamma_G$.
\item $\lambda^{(m)} \in \Omega(C)$.
\end{enumerate}
Applying successively Lemma \ref{lem-coho-case-F-mu+alpha-generated-by-F-mu} to $\lambda^{(0)} $, $\lambda^{(1)} $, $\cdots$, until $\lambda^{(m)} $, we deduce that $$F_G^{\lambda} \subset \ms H_{G, u} \cdot F_G^{\lambda^{(m)} }.$$
Since $\lambda^{(m)} < \mu_0$, 
by \ref{subsection-F-G-mu-1-includ-F-G-mu-2}, we have $F_G^{\lambda^{(m)}} \subset F_G^{\mu_0}.$ 
Thus $F_G^{\lambda} \subset \ms H_{G, u} \cdot F_G^{ \mu_0 }.$ 
\cqfd

\quad

\begin{notation}  \label{constr-h-G-omega}
For any quasi-fundamental coweight $\omega_i$ of $G$, let $h^G_{\omega_i} \in \ms H_{G, u} $
be the characteristic function of $G(\mc O_u) \varpi^{\omega_i} G(\mc O_u)$.
\end{notation}

\sssec{}
Suppose that $\mu-\omega_i l - \frac{1}{r}\check{\alpha}_i $ is dominant. By \ref{subsection-Hecke-G-action}, $h^G_{\omega_i}$ induces morphisms
\begin{equation*}
h^G_{\omega_i}: F_G^{\mu - \omega_i l} \rightarrow F_G^{\mu } .
\end{equation*}
\begin{equation*}
h^G_{\omega_i}: F_G^{\mu- \omega_i l - \frac{1}{r}\check{\alpha}_i } \rightarrow F_G^{\mu - \frac{1}{r} \check{\alpha}_i } .
\end{equation*}
They induce a morphism on the quotient spaces
\begin{equation}    \label{equation-h-G-omega-acts-on-grading}
F_{G}^{\mu - \omega_i l} / F_{G}^{\mu- \omega_i l -\frac{1}{r}\check{\alpha}_i }  \xrightarrow{ h^G_{\omega_i} } F_{G}^{\mu } / F_{G}^{\mu - \frac{1}{r} \check{\alpha}_i}  .
\end{equation}

\begin{lem}    \label{lem-Hecke-operator-on-filtration-G}
Let $\mu \in \frac{1}{r} \wh R_{G^{\mr{ad}}}^+$ such that $\langle \mu, \alpha_j \rangle > C(G, u)$ for all $j \in \Gamma_G$, then for  
any $i \in \Gamma_G$, the morphism $($\ref{equation-h-G-omega-acts-on-grading}$)$
is an isomorphism.
\end{lem}

The proof of this lemma will be given in \ref{subsection-proof-key-lemma}. The following figure shows an illustration of Lemma \ref{lem-Hecke-operator-on-filtration-G} for $G=GL_3$, 
and $i=2$, where the set $S_{M_2}(\mu)$ is defined in \ref{subsection-S-M-mu} below.
\begin{center}
\includegraphics[width=0.4\textwidth]{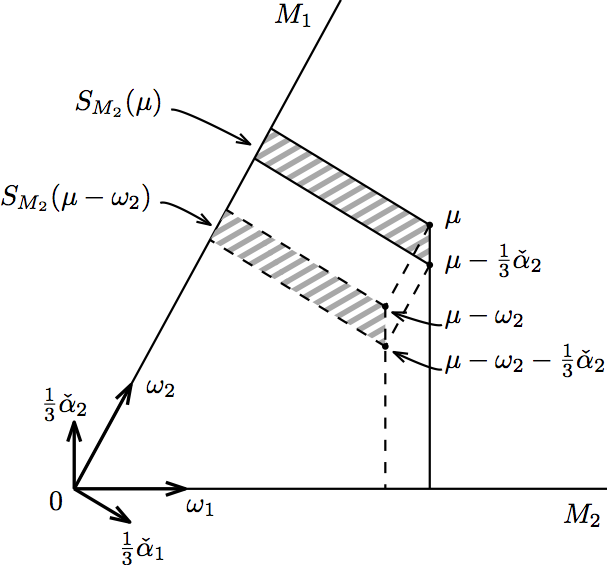} 
\end{center}

\noindent {\bf Proof of Lemma \ref{lem-coho-case-F-mu+alpha-generated-by-F-mu} admitting Lemma \ref{lem-Hecke-operator-on-filtration-G}.} 
We deduce from Lemma \ref{lem-Hecke-operator-on-filtration-G} that
\begin{equation}
F_G^{\mu } = h^G_{\omega_i} \cdot F_G^{\mu-\omega_i l} + F_G^{\mu-\frac{1}{r}\check{\alpha}_i }.
\end{equation}

For any $i \in \Gamma_G$, we have $\omega_i \geq \frac{1}{r}\check{\alpha}_i $. Indeed, since $\omega_i \in \frac{1}{r} \wh R_{G^{\mr{ad}}}^+$, we have $\omega_i = \sum_j \frac{c_{ij}}{r} \check{\alpha}_j$ for some $c_{ij} \in \Z_{\geq 0}$. Since $\langle \omega_i, \alpha_j \rangle=0$ for all $j \neq i \in \Gamma_G$, we have $\langle \omega_i, \omega_i \rangle = \langle \omega_i,  \frac{c_{ii}}{r} \alpha_i   \rangle = \frac{c_{ii}}{r} \langle \omega_i, \alpha_i \rangle$. The fact that $\langle \omega_i, \omega_i \rangle>0$ implies $c_{ii} \neq 0$. Thus $c_{ii} \geq 1$.

Taking into account that $l \geq 1$, we have $\mu-\omega_i l \leq \mu-\frac{1}{r}\check{\alpha}_i$. Thus $F_G^{\mu-\omega_i l} \subset F_G^{\mu-\frac{1}{r}\check{\alpha}_i }$. As a consequence, 
\begin{equation}
F_G^{\mu } \subset h^G_{\omega_i} \cdot F_G^{\mu-\frac{1}{r}\check{\alpha}_i} + F_G^{\mu-\frac{1}{r}\check{\alpha}_i } \subset \ms H_{G, u} \cdot F_G^{\mu-\frac{1}{r}\check{\alpha}_i} .
\end{equation}
\cqfd

\quad

The goal of the remaining part of Section 2 is to prove Lemma \ref{lem-Hecke-operator-on-filtration-G}. We need to use the cohomology group $H_{M_i}$ of the stack of $M_i$-shtukas, the constant term morphism from $H_G$ to $H_{M_i}$ and the contractibility of deep enough Harder-Narasimhan strata of $\Cht_G$. We recall these in Section 2.3. Then we study the action of $h^G_{\omega_i}$ on $H_{M_i}$ in Section 2.4. Finally in Section 2.5, we use Section 2.3 and Section 2.4 to prove Lemma \ref{lem-Hecke-operator-on-filtration-G}.

\subsection{Cohomology of stacks of $M$-shtukas and constant term morphism}    

For the convenience of the reader, we extract some results from \cite{cusp-coho}. The goal of this subsection is to state Corollary \ref{cor-CT-gr-G-equal-gr-M}.

As before, let $I$ be a finite set and $W$ be a representation of $\wh G^I$. Let $P$ be a standard parabolic subgroup of $G$. Let $M$ be the Levi quotient of $P$.

\sssec{}  
In \cite[\S~2.6.2]{cusp-coho}, we defined the classifying stack of $M$-shtukas $\Cht_{M, N, I, W}$ over $(X \sm N)^I$. It is a Deligne-Mumford stack locally of finite type. 

In \cite[\S~2.6.3]{cusp-coho}, for any $\mu \in \wh \Lambda_{G^{\mr{ad}}}^{+, \Q} $, we defined an open substack $\Cht_{M, N, I, W}^{\leq \mu}$
of $\Cht_{M, N, I, W}$. We have $\Cht_{M, N, I, W} = \bigcup_{ \mu \in \wh \Lambda_{G^{\mr{ad}}}^{+, \Q} } \Cht_{M, N, I, W}^{\leq \mu} $.

\sssec{}
Let $Z_M$ be the center of $M$, $\wh{\Lambda}_{Z_M/ Z_G}$ the coweight lattice of $Z_M / Z_G$ and
$\wh{\Lambda}_{Z_M/ Z_G}^{\Q} = \wh{\Lambda}_{Z_M/ Z_G} \otimes_{\Z} \Q.$
In \cite[Definition~1.7.2 and~2.6.3]{cusp-coho}, for any $\nu \in \wh{\Lambda}_{Z_M/ Z_G}^{\Q}$, we defined an open and closed substack $\Cht_{M, N, I, W}^{\leq \mu, \, \nu} $ of $\Cht_{M, N, I, W}^{\leq \mu}$. We have  $\Cht_{M, N, I, W}^{\leq \mu}  = \bigsqcup_{  \nu \in \wh{\Lambda}_{Z_M/ Z_G}^{\Q}   } \Cht_{M, N, I, W}^{\leq \mu, \, \nu}  $. 

Let $\Xi \subset Z_G(F) \backslash Z_G(\mb A) $ as before. We defined the quotient $\Cht_{M, N, I, W}^{\leq \mu, \, \nu} / \Xi$. By \cite[\S~2.6.3]{cusp-coho}, it is a Deligne-Mumford stack of finite type.

\begin{rem}
$\Xi$ is a lattice in $Z_G(F) \backslash Z_G(\mb A)$ but only a discrete subgroup in $Z_M(F) \backslash Z_M(\mb A)$. 
Thus the stack $\Cht_{M, N, I, W}^{\leq \mu} / \Xi$ is only locally of finite type.
\end{rem}

\sssec{}    \label{subsection-definition-of-cone-Lambda-mu}
As in \cite[\S~1.5.11]{cusp-coho}, let $$\pr_P^{ad}: \wh \Lambda_{G^{\mr{ad}}}^{\Q} \rightarrow \wh{\Lambda}_{Z_M/ Z_G}^{\Q}$$ be the projection.
In \cite[\S~1.5.13]{cusp-coho}, we defined a partial order "$\leq$" on $\wh{\Lambda}_{Z_M/ Z_G}^{\Q}$ which is induced by the partial order on $\wh \Lambda_{G^{\mr{ad}}}^{\Q}$. 
We defined a translated cone $\wh{\Lambda}_{Z_M / Z_G}^{\mu}:=\{  \nu \in \wh{\Lambda}_{Z_M / Z_G}^{\Q}, \, \nu \leq \pr_P^{ad}(\mu)    \}$ in $\wh{\Lambda}_{Z_M/ Z_G}^{\Q}$.
By \cite[Lemma~1.5.14]{cusp-coho}, if $\nu \notin \wh{\Lambda}_{Z_M / Z_G}^{\mu} \cap \pr_P^{ad}(  \frac{1}{r} \wh R_{G^{\mr{ad}}}^+)$, then the stack $\Cht_{M, N, I, W}^{\leq \mu, \, \nu}$ is empty.

\sssec{}   \label{subsection-Cht-M-prime}
In \cite[Definition~3.4.2]{cusp-coho}, we defined
$\Cht_{M, N, I, W}' :=\Cht_{M, N, I, W} \overset{P(\mc O_N)}{\times} G(\mc O_N)$. In \cite[Definition~3.4.5]{cusp-coho}, we defined $\Cht_{M, N, I, W}^{' \, \leq \mu}$ and $\Cht_{M, N, I, W}^{' \, \leq \mu, \, \nu}$ in the same way.

\sssec{}
In \cite[Definition~3.4.7]{cusp-coho} ,
we defined a Satake perverse sheaf $\mc F_{M, N, I, W}^{' \, \Xi}$ over $\Cht_{M, N, I, W}' / \Xi$ (with the perverse normalization relative to $(X \sm N)^I$).

\begin{defi} [\emph{cf.} Definitions~3.4.7 and~3.4.9 in \cite{cusp-coho}]  \label{def-H-M-N-I-W-leq-mu-nu}
Let $\mf p_M: \Cht_{M, N, I, W}' \rightarrow (X \sm N)^I$ be the morphism of paws.
For any $\mu \in \wh \Lambda_{G^{\mr{ad}}}^{+, \Q}$ and $\nu \in \wh{\Lambda}_{Z_M / Z_G}^{\mu}$, for any geometric point $\ov{x}$ of $(X \sm N)^I$, we define
$$\mc H_{M, N, I, W}^{' \, j, \, \leq\mu, \, \nu} := R^j (\mf p_M)_! ( \restr{ \mc F_{M, N, I, W}^{' \, \Xi} }{ \Cht_{M, N, I, W}^{' \, \leq \mu, \, \nu} / \Xi } ) $$
$$H_{M, N, I, W, \ov{x}}^{' \, j, \, \leq\mu, \, \nu} := \restr{ \mc H_{M, N, I, W}^{' \, j, \, \leq\mu, \, \nu} }{ \ov{x}} $$
$$H_{M, N, I, W, \ov{x}}^{' \, j, \, \leq\mu} : = \prod_{\nu \in \wh{\Lambda}_{Z_M / Z_G}^{\mu}}  H_{M, N, I, W, \ov{x}}^{' \, j, \, \leq\mu, \, \nu} $$
$$H_{M, N, I, W, \ov{x}}^{' \, j}:= \varinjlim _{\mu} H_{M, N, I, W, \ov{x} }^{' \, j, \, \leq\mu} $$
We also define 
$$H_{M, N, I, W, \ov{x}}^{' \, j, \, \nu} := \varinjlim _{\mu} H_{M, N, I, W, \ov{x}}^{' \, j, \, \leq\mu, \, \nu}. $$
\end{defi}
By definition, $H_{M, N, I, W, \ov{x}}^{' \, j, \, \leq\mu, \, \nu}$ is a finite dimensional $\Ql$-vector space and it depends on $\Xi$.

\quad

\begin{nota}
From now on, we will omit the indices $N, I, W, \ov{x}$ and $'$ to simplify the notation.
\end{nota}

\sssec{}  \label{subsection-S-M-mu}
Let $\mu \in  \wh{\Lambda}_{G^{\mr{ad}}}^{+, \Q}  $. 
In \cite[Definition~4.1.1]{cusp-coho}, 
we defined a bounded set:
$$ S_M(\mu) :  = \{ \lambda \in  \wh{\Lambda}_{G^{\mr{ad}}}^{+, \Q}  \; | \; \lambda \leq \mu \} \bigcap \{  \lambda \in \wh{\Lambda}_{G^{\mr{ad}}}^{+, \Q} \; | \;  \pr_P^{ad}(\lambda) = \pr_P^{ad}(\mu)  \} .  $$
where $\leq$ is the partial order in $\wh \Lambda_{G^{\mr{ad}}}^{\Q}$.

\sssec{}    \label{subsection-H-M-S-mu-equal-H-M-leq-mu-nu}
In \cite[Definition~4.1.10]{cusp-coho}, we defined a locally closed substack $\Cht_{G}^{S_M(\mu)}$
of $\Cht_{G}$ and an open and closed substack $\Cht_{M}^{S_M(\mu)}$ of $\Cht_M$. In \cite[Definition~4.6.2]{cusp-coho}, we defined $$H_{G}^{j, \, S_M(\mu  )}  := H_c^j ( \restr{ \Cht_{G}^{S_M(\mu  )} / \Xi }{ \ov{x} } , \mc F_{G}) ;$$
$$H_{M}^{j, \, S_M(\mu  )}  := H_c^j ( \restr{ \Cht_{M}^{S_M(\mu  )} / \Xi }{ \ov{x} } , \mc F_{M}) .$$

\sssec{}    \label{subsection-defi-CT-morphism}
In \cite[\S~3.5.8 and~4.6.3]{cusp-coho} for any $\mu \in  \frac{1}{r}\wh R_{G^{\mr{ad}}}^{+} $, we defined the constant term morphism:
$$C_G^{P, \, j, \, \leq \mu}: H_G^{j, \, \leq \mu} \rightarrow H_M^{j, \, \leq \mu} ;$$
$$C_G^{P, \, j, \, S_M( \mu )}: H_G^{j, \, S_M( \mu )} \rightarrow H_M^{j, \, S_M( \mu )} .$$
(In fact Definition \ref{def-H-M-N-I-W-leq-mu-nu} is valid for any geometric point $\ov x$. It is here that we need the conditions in \ref{subsection-defi-geo-point-x-bar}. 
See Remark \ref{rem-ov-x-ov-eta-I} for more discussion.)

\begin{notation}     \label{notation-P-i-M-i-alpha-i-omega-i}
In the following, let $i \in \Gamma_G$ and $P_i$ be a maximal parabolic subgroup with Levi quotient $M_i$ such that $\Gamma_G - \Gamma_{M_i} =  \{ i \}$. We note $P:=P_i$, $M:=M_i$, $\alpha:=\alpha_i$, $\check{\alpha}:=\check{\alpha}_i$ and $\omega:=\omega_i$.
\end{notation}

\sssec{}   \label{subsection-exact-sequences-associated-to-open-and-close}
Let $\check{\alpha}$ as in Notation \ref{notation-P-i-M-i-alpha-i-omega-i}. Note that it is a simple coroot of $G$ but not a simple coroot of $M$. Let $\mu \in  \frac{1}{r}\wh R_{G^{\mr{ad}}}^{+} $ and $\mu \geq \frac{1}{r} \check{\alpha}$. 
By the proof of Lemma 5.3.1 (1) in \cite{cusp-coho}, we have an open substack and the complementary closed substack in $\Cht_G^{\leq \mu} $:
$$\Cht_G^{\leq \mu - \frac{1}{r} \check{\alpha}  } \rightarrow \Cht_G^{\leq \mu} \leftarrow \Cht_G^{S_M(\mu)} $$
We have a long exact sequence of compact support cohomology groups:
\begin{equation}    \label{equation-coho-case-long-exact-sequence-H-G}
\cdots \rightarrow   H_{G}^{j, \, \leq \mu - \frac{1}{r} \check{\alpha} } \xrightarrow{f_j} H_{G}^{j, \, \leq \mu  } \rightarrow H_{G}^{j, \, S_M(\mu )} \rightarrow H_{G}^{j+1, \, \leq \mu - \frac{1}{r} \check{\alpha} }  \xrightarrow{f_{j+1}} H_{G}^{j+1, \, \leq \mu }  \rightarrow \cdots
\end{equation}

\sssec{}    \label{subsection-H-M-short-exact-sequence}
By the proof of Lemma 5.3.1 (1) of \cite{cusp-coho}, we also have an open substack and the complementary closed substack in $\Cht_M^{\leq \mu} $:
$$\Cht_M^{\leq \mu - \frac{1}{r} \check{\alpha}  } \rightarrow \Cht_M^{\leq \mu} \leftarrow \Cht_M^{S_M(\mu)} $$
Moreover, by \cite[\S~4.1.7]{cusp-coho}, $\Cht_M^{S_M(\mu)} = \Cht_M^{\leq \mu, \, \pr_P^{ad}(\mu)}$. So it is open and closed in $\Cht_M^{\leq \mu}$. Thus we have a (split) short exact sequence of compact support cohomology groups for every degree $j$:
$$0 \rightarrow H_{M}^{j, \, \leq \mu - \frac{1}{r} \check{\alpha}   }  \rightarrow  H_{M}^{j, \, \leq \mu }   \rightarrow H_{M}^{j, \, S_M(\mu )}  \rightarrow 0$$
And we have $H_{M}^{j, \, S_M(\mu )} = H_{M}^{j, \, \leq \mu, \, \pr_P^{ad}(\mu)} $.

\begin{prop} [\emph{cf.} Lemma~A.0.8, Proposition~4.6.4 and Proposition~5.1.5 in \cite{cusp-coho}] \label{prop-TC-suite-exacte}
Let $\mu \in  \frac{1}{r}\wh R_{G^{\mr{ad}}}^{+} $ and $\mu \geq \frac{1}{r} \check{\alpha}$. Let $j \in \Z$. Let $C_G$ be the constant in \cite{cusp-coho} Proposition 5.1.5.
\begin{enumerate}[label=$(\roman*)$]
\item We have a commutative diagram:
\begin{equation}    \label{equation-coho-case-exact-sequence-and-CT}
\xymatrix{
H_{G}^{j, \, \leq \mu - \frac{1}{r} \check{\alpha} } \ar[r]^{f_j}   \ar[d]^{C_G^{P, \, j, \, \leq \mu - \frac{1}{r} \check{\alpha}}}
& H_{G}^{j, \, \leq \mu   }   \ar[r]   \ar[d]^{C_G^{P, \, j, \, \leq \mu  }}
& H_{G}^{j, \, S_M(\mu )}  \ar[d]^{C_G^{P, j, \, \, S_M(\mu)}} \\
H_{M}^{j, \, \leq \mu - \frac{1}{r} \check{\alpha} }   \ar[r]   
& H_{M}^{j, \, \leq \mu  }   \ar[r] 
& H_{M}^{j, \, S_M(\mu)}.  
}
\end{equation}
\item If $\langle  \mu, \alpha \rangle > C_G$, then the morphism $C_G^{P, \, j, \, S_M(\mu)}$ is an isomorphism.
\item If moreover $\langle \mu, \gamma \rangle > C_G$ for every simple root $\gamma$ of $G$, then for every $j \in \Z$, the morphism $H_{G}^{j, \, \leq \mu} \rightarrow H_{G}^j$ is injective. In particular, $f_j$ and $f_{j+1}$ in (\ref{equation-coho-case-long-exact-sequence-H-G}) are injective. So we have a short exact sequence:
$$0 \rightarrow   H_{G}^{j, \, \leq \mu - \frac{1}{r} \check{\alpha}} \xrightarrow{f_j} H_{G}^{j, \, \leq \mu} \rightarrow H_{G}^{j, \, S_M(\mu)} \rightarrow 0.$$
\end{enumerate}
\end{prop}

\begin{rem}   \label{rem-ov-x-ov-eta-I}
In \cite{cusp-coho} the constant term morphism $C_G^{P, \, j, \, \leq \mu}$ and Proposition 5.1.5 are stated for $\ov{x} = \ov{\eta^I}$, a geometric generic point of $X^I$. But in fact, in \emph{loc.~cit.},
we can replace $\ov{\eta^I}$ by any geometric point $\ov{x}$ satisfying the conditions in \ref{subsection-defi-geo-point-x-bar} (of this paper). The same arguments go through.
\end{rem}

\begin{cor}    \label{cor-CT-gr-G-equal-gr-M}
Under the hypothesis of Proposition \ref{prop-TC-suite-exacte} $(iii)$, we deduce from $($\ref{equation-coho-case-exact-sequence-and-CT}$)$ a commutative diagram
\begin{equation}   \label{equation-CT-gr-G-equal-gr-M}
\xymatrix{
H_G^{j, \,  \leq \mu }  /  H_G^{j, \, \leq \mu - \frac{1}{r} \check{\alpha}  }  \ar[r]^{\quad \quad =}    \ar[d]^{ \ov{ C_G^{P, \, j, \, \leq \mu } } }
& H_{G}^{j, \, S_M(\mu  )} \ar[d]_{\simeq}^{C_G^{P, \, j, \, S_M(\mu)}} \\
H_M^{j, \,  \leq \mu  }  /  H_M^{j, \, \leq \mu - \frac{1}{r} \check{\alpha}  }   \ar[r]^{\quad \quad =}
& H_{M}^{j, \, S_M(\mu  )} .
}
\end{equation}
Thus the left vertical map $\ov{ C_G^{P, \, j, \, \leq \mu } }$
is an isomorphism.
\end{cor}

\subsection{Action of Hecke operators}    \label{subsection-action-of-Hecke}

Let $M$ as in Notation \ref{notation-P-i-M-i-alpha-i-omega-i}. The goal of this section is to prove Lemma \ref{lem-Hecke-operator-on-filtration-M}.

\sssec{}
Let $\ms H_{M, u}:=C_{c}(M(\mc O_u)\backslash M(F_u)/M(\mc O_u), \Ql)$ be the Hecke algebra of $M$ at the place $u$. 
The cohomology group $H_{M}^{j}$ is equipped with an action of $\ms H_{M, u}$ by Hecke correspondences (see \cite[\S~2.20 and~4.4]{vincent} and \cite[\S~6.1-6.2]{cusp-coho}) and equipped with an action of $\ms H_{G, u}$ via the Satake transformation $\ms H_{G, u} \hookrightarrow \ms H_{M, u}$.

\begin{defi}
We define a filtration $F_M^\bullet$ of $H_M^j$:
$$\mu \in \frac{1}{r} \wh R_{G^{\mr{ad}}}^+, \quad F_M^{\mu}:= \on{Im} ( H_{M}^{j, \, \leq \mu} \rightarrow H_M^j ) $$
\end{defi}

\begin{lem}    \label{lem-Hecke-operator-on-filtration-M}
Let $\omega$ as in Notation \ref{notation-P-i-M-i-alpha-i-omega-i} and 
$h^G_{\omega}$ as in Notation \ref{constr-h-G-omega}.
\begin{enumerate}[label=$(\roman*)$]
\item For any $\lambda  \in \frac{1}{r} \wh R_{G^{\mr{ad}}}^+$, the Hecke operator $h^G_{\omega}$ induces a morphism 
\begin{equation}   \label{equation-F-M-lambda-to-F-M-lambda+omega}
F_{M}^{\lambda} \xrightarrow{h^G_{\omega}} F_{M}^{\lambda + \omega l}  .
\end{equation}
\item If $\lambda - \frac{1}{r} \check{\alpha}$ is still in $\frac{1}{r} \wh R_{G^{\mr{ad}}}^+$, the induced morphism on the quotients
\begin{equation}    \label{equation-gr-M-lambda-to-gr-M-lambda+alpha}
F_{M}^{\lambda} / F_{M}^{\lambda -\frac{1}{r}\check{\alpha} }  \xrightarrow{ h^G_{\omega} } F_{M}^{\lambda + \omega l } / F_{M}^{\lambda + \omega l - \frac{1}{r} \check{\alpha}} 
\end{equation}
is an isomorphism.
\end{enumerate}
\end{lem}

To prove this lemma we need some preparations.

\sssec{}
We define $$\Omega:=\{   \theta \text{ coweight of } G \text{ conjugate to } \omega , \; \theta \text{ is dominant for } M \}.$$
Since $\omega$ is dominant for $G$, for any $\theta \in \Omega$, we have $\theta \leq \omega$ for the order in $\wh \Lambda_{G}$.

\sssec{}
For any dominant coweight $\theta$ of $M$,
let $h^M_{\theta} \in \ms H_{M, u}$
be the characteristic function of $M(\mc O_u) \varpi^\theta M(\mc O_u)$.

We can view $h^G_{\omega}$ as an element in $\ms H_{M, u}$ by the Satake transformation $\ms H_{G, u} \hookrightarrow \ms H_{M, u}$. We have the following equality (up to multiplication by a power of $q$) in $\ms H_{M, u}$:
\begin{equation}   \label{equation-h-G-omega=sum-h-M-theta}
h^G_{\omega} =  h^M_{\omega} + \sum_{\theta \in \Omega, \, \theta \neq \omega} c_{\theta} h^M_{\theta}, \quad c_{\theta} \in \Q_{\geq 0}.
\end{equation}

\sssec{}    \label{subsection-omega-theta-geq-alpha}
For any $\theta \in \Omega, \, \theta \neq \omega$, we have $\omega - \theta \geq \frac{1}{r}\check{\alpha} $ for the order in $\wh \Lambda_{G}^{\Q}$ (hence for the order in $\wh \Lambda_{G^{\mr{ad}}}^{\Q}$). 
Indeed, we have $\theta<\omega$, i.e. $\omega - \theta = \sum_j \frac{c_{j}}{r} \check{\alpha}_j$ for some $c_{j} \in \Z_{\geq 0}$. Recall that $\check{\alpha} = \check{\alpha}_i$. Then the fact that $\langle \omega, \alpha_j \rangle=0$ for all $j \neq i \in \Gamma_G$ implies $\langle \omega , \omega- \theta \rangle = \langle \omega ,  \frac{c_{i}}{r} \alpha_i   \rangle = \frac{c_{i}}{r} \langle \omega , \alpha_i \rangle$. Since $\| \theta \|^2 = \|\omega \|^2$, we cannot have $\langle \omega , \omega- \theta \rangle = 0$.
This implies $c_{i} \neq 0$. Thus $c_{i} \geq 1$. 

Moreover, since $l \geq 1$, we have $\omega l - \theta l \geq \frac{1}{r}\check{\alpha} $.

\quad

\noindent {\bf Proof of Lemma \ref{lem-Hecke-operator-on-filtration-M}.} 
Recall that in \ref{subsection-definition-of-cone-Lambda-mu} we defined a projection $\pr_P^{ad}: \wh \Lambda_{G^{\mr{ad}}}^{\Q} \rightarrow \wh{\Lambda}_{Z_M/ Z_G}^{\Q}$. We will write $\pr:=\pr_P^{ad}$ to simplify the notations.

(i) Let $\lambda \in \frac{1}{r} \wh R_{G^{\mr{ad}}}^+$ and $\nu \in \wh \Lambda_{Z_M / Z_G}^{\Q}$.
Note that $\omega$ is a coweight of the center of $M$.
Thus the Hecke operator $h^M_{\omega}$ induces an isomorphism of stacks (\emph{cf.} \cite[\S~ 1.3.1]{cusp-coho} applied to the reductive group $M$ and (2) of the proof of Lemma 6.3.3):
$$\Cht_{M}^{\leq \lambda, \, \nu} / \Xi \isom \Cht_{M}^{\leq \lambda+\omega l , \, \nu+\pr(\omega) l} / \Xi $$
This induces an isomorphism of cohomology groups:
\begin{equation}   \label{equation-coho-case-H-M-nu-isom-H-M-nu+omega}
h^M_{\omega}: H_M^{j, \, \leq \lambda, \, \nu} \isom H_M^{j, \, \leq \lambda+\omega l, \, \nu+\pr(\omega) l} 
\end{equation}

Now consider the Hecke operator $h^M_{\theta}$ for $\theta \in \Omega$ and $\theta \neq \omega$. As in \cite[\S~2.20 and~4.4]{vincent} and \cite[\S~6.1-6.2]{cusp-coho}, we have the Hecke correspondence associated to $\theta$:
$$\Cht_M / \Xi \xleftarrow{\pr_1} \Gamma(M(\mc O_u) \varpi^{\theta} M(\mc O_u) ) \xrightarrow{\pr_2} \Cht_M / \Xi$$
where $\pr_1$ and $\pr_2$ are finite étale morphisms.
The morphism $\pr_1$ sends $\pr_2^{-1}( \Cht_M^{\leq \lambda, \, \nu} / \Xi )$ to $\Cht_M^{\leq \lambda+\theta l, \, \nu+ \pr(\theta) l} / \Xi $.
This induces 
a morphism of cohomology groups:
\begin{equation}    \label{equation-coho-case-H-M-nu-to-H-M-nu+theta}
h^M_{\theta}: H_M^{j, \, \leq \lambda, \, \nu} \rightarrow H_M^{j, \, \leq \lambda+\theta l ,  \, \nu+\pr(\theta) l}
\end{equation}
Taking into account (\ref{equation-h-G-omega=sum-h-M-theta}), we deduce from (\ref{equation-coho-case-H-M-nu-isom-H-M-nu+omega}) and (\ref{equation-coho-case-H-M-nu-to-H-M-nu+theta}) a morphism
\begin{equation}  \label{equation-coho-case-H-M-lambda-to-H-M-lambda+theta-all-theta}
H_M^{j, \, \leq \lambda, \, \nu} \xrightarrow{  h^G_{\omega} }  \left( \prod_{\theta \in \Omega, \, \theta \neq \omega} H_M^{j, \, \leq \lambda+\theta l, \, \nu+\pr(\theta) l}  \right) \times H_M^{j, \, \leq \lambda+\omega l, \, \nu+\pr(\omega) l} .
\end{equation}
Consider the RHS of (\ref{equation-coho-case-H-M-lambda-to-H-M-lambda+theta-all-theta}). For any $\theta \in \Omega$ and $\theta \neq \omega$, by \ref{subsection-omega-theta-geq-alpha} we have $\theta l \leq \omega l -\frac{1}{r}\check{\alpha}$. Thus %(for the component indexed by $\nu+pr(\theta) l$), we have 
\begin{equation}    \label{equation-H-M-theta-smaller-than-H-M-omega-alpha}
H_M^{j, \, \leq \lambda+\theta l, \, \nu+\pr(\theta) l} \subset H_M^{j, \, \leq \lambda+\omega l - \frac{1}{r}\check{\alpha}, \, \nu+\pr(\theta) l} .
\end{equation}
By \cite[Proposition~3.1.2]{schieder}, $\pr$ is order-preserving. Thus $\pr(\theta) l \leq \pr(\omega) l - \pr(\frac{1}{r} \check{\alpha}) $. We deduce that 
%(the LHS is a factor of the RHS)
\begin{equation}  \label{equation-H-M-omega-alpha-smaller-than-prod}
H_M^{j, \, \leq \lambda+\omega l - \frac{1}{r}\check{\alpha}, \, \nu+\pr(\theta) l}  \subset \prod_{\nu' \leq \nu+\pr(\omega) l - \frac{1}{r}\pr(\check{\alpha})} H_M^{j, \, \leq \lambda+\omega l - \frac{1}{r}\check{\alpha}, \, \nu'} .
\end{equation}
%H_M^{j, \, \leq \lambda+\theta l, \, \nu+pr(\theta) l} \subset H_M^{j, \, \leq \lambda+\omega l, \, \nu+pr(\theta) l}  \subset \prod_{\nu' < \nu+pr(\omega) l} H_M^{j, \, \leq \lambda+\omega l, \, \nu'} .
Moreover, by \emph{loc.~cit.} $\pr(\check{\alpha}) > 0$, thus 
\begin{equation}   \label{equation-prod-H-M-omega-alpha-smaller-than-prod-omega}
\prod_{\nu' \leq \nu+\pr(\omega) l - \frac{1}{r}\pr(\check{\alpha})} H_M^{j, \, \leq \lambda+\omega l - \frac{1}{r}\check{\alpha}, \, \nu'}  \subset  \prod_{\nu' < \nu+\pr(\omega) l } H_M^{j, \, \leq \lambda+\omega l, \, \nu'}.
\end{equation}
Combining (\ref{equation-H-M-theta-smaller-than-H-M-omega-alpha}), (\ref{equation-H-M-omega-alpha-smaller-than-prod}) and (\ref{equation-prod-H-M-omega-alpha-smaller-than-prod-omega}), we deduce that the RHS of (\ref{equation-coho-case-H-M-lambda-to-H-M-lambda+theta-all-theta}) can be sent in 
$$\left( \prod_{\nu' < \nu+\pr(\omega) l } H_M^{j, \, \leq \lambda+\omega l, \, \nu'}  \right) \times H_M^{j, \, \leq \lambda+\omega l, \, \nu+\pr(\omega) l} ,$$
which is nothing but $\prod_{\nu' \leq \nu+\pr(\omega) l } H_M^{j, \, \leq \lambda+\omega l, \, \nu'}.$
So (\ref{equation-coho-case-H-M-lambda-to-H-M-lambda+theta-all-theta}) induces a morphism
$$h^G_{\omega}  : H_M^{j, \, \leq \lambda, \, \nu} \rightarrow \prod_{\nu' \leq \nu+\pr(\omega) l } H_M^{j, \, \leq \lambda+\omega l, \, \nu'} .$$
Now fixing $\lambda$ and taking product over all $\nu' \in \wh \Lambda_{Z_M / Z_G}^{\Q}$ such that $\nu' \leq \pr(\lambda)$, we obtain:
\begin{equation}  \label{equation-prod-H-M-leq-lambda-to-H-M-leq-lambda+omega} 
h^G_{\omega}:  \prod_{\nu' \leq \pr(\lambda)} H_M^{j, \,  \leq \lambda, \, \nu'} \rightarrow \prod_{\nu' \leq \pr(\lambda)+\pr(\omega) l } H_M^{j, \,  \leq \lambda+\omega l, \, \nu'} . 
\end{equation}
By Definition \ref{def-H-M-N-I-W-leq-mu-nu}, $$\prod_{\nu' \leq \pr(\lambda)} H_M^{j, \,  \leq \lambda, \, \nu'} = H_M^{j, \,  \leq \lambda}  \; ; \quad \prod_{\nu' \leq \pr(\lambda)+\pr(\omega) l } H_M^{j, \,  \leq \lambda+\omega l, \, \nu'} = H_M^{j, \,  \leq \lambda+\omega l} .$$ We deduce from (\ref{equation-prod-H-M-leq-lambda-to-H-M-leq-lambda+omega}) the morphism (\ref{equation-F-M-lambda-to-F-M-lambda+omega}).

\quad

(ii) 
By \ref{subsection-H-M-short-exact-sequence}, we have $$H_{M}^{j, \, \leq \lambda} / H_{M}^{j, \, \leq \lambda- \frac{1}{r}\check{\alpha}} =  H_M^{j, \, S_M(\lambda)}  =  H_M^{j, \, \leq \lambda, \, \pr(\lambda)  } ; $$ 
$$H_{M}^{j, \, \leq \lambda + \omega l } / H_{M}^{j,\,\leq \lambda - \frac{1}{r} \check{\alpha} + \omega l} = H_M^{j, \, S_M(\lambda + \omega l ) } = H_M^{j, \, \leq \lambda + \omega l, \, \pr(\lambda + \omega l) } .$$
For any $\theta \in \Omega$ and $\theta \neq \omega$, by (\ref{equation-coho-case-H-M-nu-to-H-M-nu+theta}), (\ref{equation-H-M-theta-smaller-than-H-M-omega-alpha}) and (\ref{equation-H-M-omega-alpha-smaller-than-prod}), $h^M_{\theta}$ induces a morphism
$$H_M^{j, \, \leq \lambda, \, \pr(\lambda) }  \rightarrow \prod_{\nu' \leq \pr(\lambda)+\pr(\omega) l - \frac{1}{r}\pr(\check{\alpha})} H_M^{j, \, \leq \lambda+\omega l - \frac{1}{r}\check{\alpha}, \, \nu'} = H_M^{j, \, \leq \lambda+\omega l - \frac{1}{r}\check{\alpha} } . $$
Thus the morphism
$$H_{M}^{j, \, \leq \lambda} / H_{M}^{j, \, \leq \lambda- \frac{1}{r}\check{\alpha}}    \xrightarrow{h^M_{\theta}}     H_{M}^{j, \, \leq \lambda + \omega l } / H_{M}^{j,\,\leq \lambda - \frac{1}{r} \check{\alpha} + \omega l} $$
is the zero morphism.

By (\ref{equation-h-G-omega=sum-h-M-theta}) and (\ref{equation-coho-case-H-M-nu-isom-H-M-nu+omega}), we deduce that $$H_{M}^{j, \, \leq \lambda} / H_{M}^{j, \, \leq \lambda- \frac{1}{r}\check{\alpha}}    \xrightarrow{h^G_{\omega}}     H_{M}^{j, \, \leq \lambda + \omega l } / H_{M}^{j,\,\leq \lambda - \frac{1}{r} \check{\alpha} + \omega l} $$ is an isomorphism.
Finally, 
we deduce that
(\ref{equation-gr-M-lambda-to-gr-M-lambda+alpha}) is an isomorphism.
\cqfd

\subsection{Proof of Lemma \ref{lem-Hecke-operator-on-filtration-G}  }    \label{subsection-proof-key-lemma}

By \cite[Lemma~6.2.12]{cusp-coho}, the action of $\ms H_{G, u}$ commutes with $C_G^P$. 
The following diagram is commutative, where the horizontal lines are the exact sequences in \ref{subsection-exact-sequences-associated-to-open-and-close} and \ref{subsection-H-M-short-exact-sequence}, the vertical morphisms are constant term morphisms:

{\resizebox{12cm}{!}{ $$
\xymatrix{
& F_{G}^{\mu- \omega l - \frac{1}{r}\check{\alpha}}  \ar[rr]   \ar[dd] \ar[ld]_{h^G_\omega}   \ar@{}[lddd]|{(a)}
& & F_{G}^{\mu - \omega l}   \ar[rr]   \ar[dd]    \ar[ld]_{h^G_\omega}    \ar@{}[lddd]|{(b)}
& & F_{G}^{\mu - \omega l} / F_{G}^{\mu- \omega - \frac{1}{r}\check{\alpha}}     \ar[dd] \ar[ld]_{h^G_\omega}    \ar@{}[lddd]|{(c)}     \\
 F_{G}^{\mu - \frac{1}{r} \check{\alpha}}    \ar[rr]   \ar[dd]
& & F_{G}^{\mu }   \ar[rr]   \ar[dd]
& & F_{G}^{\mu } / F_{G}^{\mu - \frac{1}{r} \check{\alpha}}   \ar[dd]  \\
& F_{M}^{\mu- \omega l - \frac{1}{r}\check{\alpha}}     \ar[rr]    \ar[ld]_{h^G_\omega}
& & F_{M}^{\mu - \omega l}  \ar[rr]    \ar[ld]_{h^G_\omega}
& & F_{M}^{\mu- \omega l} / F_{M}^{\mu- \omega l - \frac{1}{r}\check{\alpha} }  \ar[ld]_{h^G_\omega} \\
 F_M^{\mu - \frac{1}{r} \check{\alpha}}    \ar[rr]
& & F_{M}^{\mu }   \ar[rr]
& & F_{M}^{\mu } / F_{M}^{\mu - \frac{1}{r} \check{\alpha}} 
}
$$
} }

Consider the right face (c):
$$
\xymatrix{
F_{G}^{\mu- \omega l} / F_{G}^{\mu- \omega l - \frac{1}{r}\check{\alpha} }  \ar[r]^{\quad h^G_{\omega}}   \ar[d]_{\ov{C_G^P}}
& F_{G}^{\mu} / F_{G}^{\mu - \frac{1}{r} \check{\alpha} }   \ar[d]_{\ov{C_G^P}} \\
F_{M}^{\mu- \omega l} / F_{M}^{\mu- \omega - \frac{1}{r}\check{\alpha} l}   \ar[r]^{\quad h^G_{\omega}}_{ \quad \simeq }
& F_{M}^{\mu } / F_{M}^{\mu - \frac{1}{r} \check{\alpha}}  
}
$$
where the isomorphism of the lower line follows from Lemma \ref{lem-Hecke-operator-on-filtration-M}.
Moreover, by Corollary \ref{cor-CT-gr-G-equal-gr-M}, the vertical morphisms are isomorphisms. We deduce that the upper line is an isomorphism.
\cqfd

\begin{rem}
To prove Lemma \ref{lem-coho-case-F-mu+alpha-generated-by-F-mu}, we only need the surjectivity of the morphism in Lemma \ref{lem-Hecke-operator-on-filtration-G}. But our proof gives the isomorphism as well.
\end{rem}

\begin{rem}
In Section 3, we only need to use Proposition \ref{prop-H-G-is-generated-by-Hecke-alg} for $\ov x = \ov{\eta^I}$. But we prefer to state Proposition \ref{prop-H-G-is-generated-by-Hecke-alg} for general $\ov x$.
\end{rem}

\begin{rem}
\leavevmode
\begin{enumerate}[label=$(\roman*)$]
\item We do not know if the statement of Theorem \ref{thm-coho-cht-is-Hecke-mod-type-fini-Hecke-en-v} is still true for $u \in |N|$ because our argument uses the Satake transform for the spherical local Hecke algebra.
\item Let $\ms H_G:=C_{c}(K_N\backslash G(\mb A)/K_N, \Ql)$ be the global Hecke algebra. We have $\ms H_G = \otimes_{u \in |X|} ' \ms H_{G, u}$. Theorem~\ref{thm-coho-cht-is-Hecke-mod-type-fini-Hecke-en-v} implies that $H_{G, N, I, W, \ov{x}}^j$ is a $\ms H_{G}$-module of finite type.
\end{enumerate}
\end{rem}

\sssec{}
Proposition \ref{prop-H-G-is-generated-by-Hecke-alg} is proved for $\ov x$ satisfying the condition in \ref{subsection-defi-geo-point-x-bar}. Now we extend it to more general geometric points.
For $(n_i)_{i \in I} \in \N^I$, the partial Frobenius morphism 
$$\prod F_{\{i\}}^{n_i}: \left(\prod \Frob_{\{i\}}^{n_i}\right)^* \mc H_{G, N, I, W}^j \isom \mc H_{G, N, I, W}^j$$
defined in \cite[Sections~3 and~4.3]{vincent} induces a commutative diagram
$$
\xymatrix{
\restr{\mc H_{G, N, I, W}^{j, \, \leq \mu_0}}{\ov x} \ar[r] \ar[d]
& \restr{\mc H_{G, N, I, W}^{j, \, \leq \mu_0+\kappa}}{\prod \Frob_{\{i\}}^{-n_i}(\ov x)} \ar[d] \\
\restr{\mc H_{G, N, I, W}^j}{\ov x} \ar[r]^{\simeq \quad \quad}
& \restr{\mc H_{G, N, I, W}^j}{\prod \Frob_{\{i\}}^{-n_i}(\ov x)}.
}
$$
Since the action of $\ms H_{G, u}$ commutes with the action of partial Frobenius morphisms, we deduce from $\restr{\mc H_{G, N, I, W}^j}{\ov x} = \ms H_{G, u} \cdot \restr{\mc H_{G, N, I, W}^{j, \, \leq \mu_0}}{\ov x} $ that 
\begin{equation}  \label{equation-H-Frob-ov-x}
\restr{\mc H_{G, N, I, W}^j}{\prod \Frob_{\{i\}}^{-n_i}(\ov x)} = \ms H_{G, u} \cdot \restr{\mc H_{G, N, I, W}^{j, \, \leq \mu_0+\kappa}}{\prod \Frob_{\{i\}}^{-n_i}(\ov x)}
\end{equation}

Let $\ov y$ be a geometric point of $X^I$ such that for every $i$, the composition $\ov{y} \rightarrow (X \sm N)^I \xrightarrow{\on{pr}_i } X \sm N $ is over the generic point $\eta$ of $X \sm N$. Then there exists $\ov x$ satisfying the condition in \ref{subsection-defi-geo-point-x-bar} and $(n_i)_{i \in I} \in \N^I$ such that $\ov y = \prod \Frob_{\{i\}}^{-n_i}(\ov x)$. (For example, for $\ov y = \Frob_{\{i\}} \Delta(\ov{\eta})$ we have $\ov y = \prod_{j \neq i} \Frob_{\{j\}}^{-1} (\Frob(\Delta(\ov{\eta})))$, where $\Frob$ is the total Frobenius so $\Frob(\Delta(\ov{\eta})) = \Delta(\ov{\eta})$, which satisfies \ref{subsection-defi-geo-point-x-bar}.) We deduce from (\ref{equation-H-Frob-ov-x}) that Proposition \ref{prop-H-G-is-generated-by-Hecke-alg} is also true for $\ov y$.

\quad

\section{Global excursion operators}

Recall that in Definition \ref{def-H-G-I-W-sheaf}, we defined an inductive limit of $\Ql$-constructible sheaves $\mc H_{G, N, I, W}^j$ over $(X \sm N)^I$.
In this section, we write $\mc H_{I, W}^j := \mc H_{G, N, I, W}^j$ to simplify the notations.

The goals of this section are Construction \ref{constr-excursion-operator-on-C-c}, Theorem \ref{thm-decomposition-of-C-c-quotient-I-by-param-Langlands} and Construction \ref{construction-excursion-operator-on-cohomology}. For this, we need a specialization morphism constructed in Section 3.1. We also need some variant of Drinfeld's lemma to get an action of $\on{Weil}(\ov F / F)^I$ on $\restr{\mc H_{I, W}^j}{\ov{\eta^I}}$ in Secton 3.2.

\subsection{Specialization morphism}     

\sssec{}  \label{subsection-fix-specialization-from-eta-I-bar-to-delta-eta-bar}
Let $I$ be a finite set. We denote by $\Delta: X \rightarrow X^I$ the diagonal morphism. We denote by $F^I$ the function field of $X^I$ and $\eta^I = \on{Spec}(F^I)$ the generic point of $X^I$. We fix an algebraic closure $\ov{F^I}$ of $F^I$ and denote by $\ov{\eta^I}=\on{Spec}(\ov{F^I})$ the geometric point over $\eta^I$. Moreover, we fix a specialization map 
\begin{equation}  \label{equation-specialisation-ov-eta-I-to-Delta-ov-eta}
\mf{sp}:\ov{\eta^I} \rightarrow \Delta(\ov{\eta}).
\end{equation}
It induces the homomorphism of specialization:
\begin{equation}    \label{equation-specialisation-H-G}
\mf{sp}^*: \restr{ \mc H_{I, W}^j }{\Delta(\ov{\eta})} \rightarrow  \restr{ \mc H_{I, W}^j  }{\ov{\eta^I}} 
\end{equation}

By \cite[Proposition~8.32]{vincent}, this morphism is injective.

\begin{prop}   \label{prop-specialisation-is-bij}
The morphism $($\ref{equation-specialisation-H-G}$)$ is a bijection.
\end{prop}

\begin{notation}   \label{notation-Frob-partial-on-X-I}
For any $i \in I$, let $\Frob_{\{i\}}: X^I \rightarrow X^I$ be the morphism sending $(x_j)_{j \in I}$ to $(x_j')_{j \in I}$, with $x_i' = \Frob(x_i)$ and $x_j' = x_j$ if $j \neq i$. 
\end{notation}

\noindent {\bf Proof of Proposition \ref{prop-specialisation-is-bij}.}  
We only need to prove the surjectivity. The proof consists of 3 steps.

\noindent\textbf{Step~1.} Let $\mu_0 \in \wh \Lambda_{G^{\mr{ad}}}^{+, \Q}$ as in Proposition \ref{prop-H-G-is-generated-by-Hecke-alg} such that 
\begin{equation}    \label{equation-H-I-W-is-Hecke-H-I-W-leq-mu}
\restr{  \mc H_{I, W}^j  }{\ov{\eta^I}} = \restr{  \ms H_{G, u} \cdot \mc H_{I, W}^{j, \, \leq \mu_0}   }{\ov{\eta^I}} ,
\end{equation}
where $\mc H_{I, W}^{j, \, \leq \mu_0} $ is the $\Ql$-constructible sheaf over $(X \sm N)^I$ defined in Definition \ref{def-H-I-W-leq-mu-sheaf}.
Let $\Omega_0$ be an open dense subscheme of $(X \sm N)^I$ such that $\restr{  \mc H_{I, W}^{j, \, \leq \mu_0}  }{\Omega_0} $ is smooth. 
By \cite[Lemma~9.2.1]{eike-lau} (the argument is recalled in the proof of Lemme 8.12 of \cite{vincent}), the set $\{ (\prod_{i \in I} \Frob_{ \{i \} }^{m_i}  )(\Delta(\eta)) , (m_i)_{i \in I} \in \N^I \}$ is Zariski dense in $X^I$. Thus there exists $(n_i)_{i \in I} \in \N^I$ such that $(\prod_{i \in I} \Frob_{ \{i \} }^{n_i}  )(\Delta(\eta)) \in \Omega_0$. 

\noindent\textbf{Step~2.} The image by $\prod_{i \in I} \Frob_{\{i \}}^{n_i} $ of the specialization map (\ref{equation-specialisation-ov-eta-I-to-Delta-ov-eta}) gives a specialization map 
$$\wt{\mf{sp}}:  \left(\prod_{i \in I} \Frob_{ \{i \} }^{n_i}  \right)(\ov{\eta^{I}})   \rightarrow \left(\prod_{i \in I} \Frob_{\{i \}}^{n_i}  \right)(\Delta(\ov{\eta})) .$$
We have a commutative diagram:
\begin{equation}     \label{diagram-specialisation-H-G-leq-mu-0-H-G}
\xymatrix{
\restr{  \mc H_{I, W}^{j, \, \leq \mu_0}  }{ (\prod \Frob  )(\Delta(\ov{\eta}))} \ar[r]^{\wt{\mf{sp}}^*}_{(a)} \ar[d]
& \restr{  \mc H_{I, W}^{j, \, \leq \mu_0}  }{ (\prod \Frob )(\ov{\eta^{I}}) }  \ar[d]  \\
\restr{  \mc H_{I, W}^j  }{ (\prod \Frob  )(\Delta(\ov{\eta}))} \ar[r]^{\wt{\mf{sp}}^*}_{(b)}
& \restr{  \mc H_{I, W}^j }{ (\prod \Frob )(\ov{\eta^{I}}) } 
}
\end{equation}
where $\prod \Frob $ is a shortcut of $\prod_{i \in I} \Frob_{\{i \}}^{n_i} $. The horizontal maps are homomorphism of specialization induced by $\wt{\mf{sp}}$. The vertical maps come from the morphism of sheaves $\mc H_{I, W}^{j, \, \leq \mu_0}  \rightarrow \mc H_{I, W}^j$.

Since $\restr{  \mc H_{I, W}^{j, \, \leq \mu_0}  }{\Omega_0} $ is smooth and $(\prod_{i \in I} \Frob_{ \{i \} }^{n_i}  )(\Delta(\eta)) \in \Omega_0$, the morphism (a) in the diagram (\ref{diagram-specialisation-H-G-leq-mu-0-H-G}) is an isomorphism.

Now we show that the morphism (b) in the diagram (\ref{diagram-specialisation-H-G-leq-mu-0-H-G}) is surjective. Since $(\prod \Frob )(\ov{\eta^{I}}) \simeq \ov{\eta^I}$, by (\ref{equation-H-I-W-is-Hecke-H-I-W-leq-mu}) we have 
\begin{equation}
\restr{  \mc H_{I, W}^j  }{(\prod \Frob )(\ov{\eta^{I}})} = \restr{  \ms H_{G, u} \cdot \mc H_{I, W}^{j, \, \leq \mu_0}   }{(\prod \Frob )(\ov{\eta^{I}})} .
\end{equation}

Let $x \in \restr{  \mc H_{I, W}^j  }{(\prod \Frob )(\ov{\eta^{I}})}$. Then $x = \sum h_k a_k$ for $h_k \in \ms H_{G, u}$ and $a_k \in  \restr{  \mc H_{I, W}^{j, \, \leq \mu_0}   }{(\prod \Frob )(\ov{\eta^{I}})} $. Since (a) is an isomorphism, there exists $b_k \in \restr{  \mc H_{I, W}^{j, \, \leq \mu_0}  }{ (\prod \Frob  )(\Delta(\ov{\eta}))}$ such that $a_k = \wt{\mf{sp}}^*(b_k)$. Moreover, the diagram (\ref{diagram-specialisation-H-G-leq-mu-0-H-G}) is compatible with the action of $\ms H_{G, u}$. Thus $x$ is the image of $\sum h_k b_k \in  \restr{  \mc H_{I, W}^j  }{ (\prod \Frob  )(\Delta(\ov{\eta}))} $. This proves the surjectivity of (b).

\noindent\textbf{Step~3.} As in the proof of Proposition 8.31 in \cite{vincent}, we have a commutative diagram
$$
\xymatrix{
  \restr{  \mc H_{I, W}^j  }{ (\prod \Frob  )(\Delta(\ov{\eta}))}  \ar[d]^{\prod_{i\in I}F_{\{i\}}^{n_{i}} }_{\simeq}   \ar[r]^{(b)}
 &   \restr{  \mc H_{I, W}^j  }{(\prod \Frob )(\ov{\eta^{I}})} \ar[d]^{\prod_{i\in I}F_{\{i\}}^{n_{i}} }_{\simeq}    \\  
 \restr{  \mc H_{I, W}^j  }{\Delta(\ov{\eta})}   \ar[r]^{(\ref{equation-specialisation-H-G})}
& \restr{ \mc H_{I, W}^j  } {\ov{\eta^{I}}}  
}
$$
where the vertical morphisms are the partial Frobenius morphisms defined in \cite[\S~4.3]{vincent}. 
The surjectivity of (b) implies the surjectivity of (\ref{equation-specialisation-H-G}).
\cqfd

\subsection{Drinfeld's lemma}

Drinfeld's lemma is first proved in \cite[Proposition~6.1]{drinfeld-compact}, then generalized in \cite[Theorem~8.1.4]{eike-lau} and \cite[lem.~8.2]{vincent}.

\sssec{}   \label{subsection-notation-fix-sp-eta-I-bar-to-Delta-eta-bar}
We fixed an algebraic closure $\ov{F^I}$ of $F^I$ and $\mf{sp}$ in \ref{subsection-fix-specialization-from-eta-I-bar-to-delta-eta-bar}.
As in \cite[rem.~8.18]{vincent}, we define 
$$\on{FWeil(\eta^I, \ov{\eta^I})} := \left\{ \delta \in \on{Aut}_{\ov{\Fq}}(\ov{F^I}) \, | \, \exists (n_i)_{i \in I} \in \Z^I, \on{ s.t.} \restr{\delta}{(F^I)^{\on{perf}}} = \prod_{i \in I} (\on{Frob}_{  \{ i \}  })^{-n_i}  \right\} ,$$
where $(F^I)^{\on{perf}}$ is the perfection of $F^I$, $\on{Frob}_{\{i\}}$ is the partial Frobenius defined in Notation \ref{notation-Frob-partial-on-X-I}.

\sssec{}
Note that we have $\pi_1(\eta, \ov{\eta}) = \on{Gal}(\ov{F} / F)$.
Let $\on{Weil}(\eta, \ov{\eta})$ be the Weil group in $\pi_1(\eta, \ov{\eta})$ (also denoted by $\on{Weil}(\ov{F} / F)$). 

We have $\pi_1(\eta^I, \ov{\eta^I}) = \on{Gal}(\ov{F^I} / F^I)$.
Let $\on{Weil}(\eta^I, \ov{\eta^I})$ be the Weil group in $\pi_1(\eta^I, \ov{\eta^I})$ (also denoted by $\on{Weil}(\ov{F^I} / F^I)$). That is to say, 
$$\on{Weil(\eta^I, \ov{\eta^I})} = \left\{ \gamma \in \on{Aut}_{F^I}(\ov{F^I}) \, | \, \exists n \in \Z, \on{ s.t.} \restr{\gamma}{\Fqbar} = \Frob^n  \right\} .$$

As in \cite[rem.~8.18]{vincent}, the specialization map $\mf{sp}$ induces an inclusion $$\ov{F} \otimes_{\ov{\Fq}} \cdots \otimes_{\ov{\Fq}} \ov{F} \subset \ov{F^I}.$$
%By restriction of automorphisms, 
We denote $\ov{F_i}:=\ov{\Fq} \otimes_{\ov{\Fq}} \cdots \otimes_{\ov{\Fq}} \ov F \otimes_{\ov{\Fq}} \cdots \otimes_{\ov{\Fq}} \ov{\Fq}$ where $\ov F$ is the $i$-th term. 
We have morphisms
$$
\pi_1(\eta^I, \ov{\eta^I})  \rightarrow \pi_1(\eta, \ov{\eta})^I  , \quad \gamma  \mapsto (\restr{\gamma}{\ov F_i})_{i \in I}
$$
$$
\on{Weil}(\eta^I, \ov{\eta^I})  \rightarrow \on{Weil}(\eta, \ov{\eta})^I  , \quad \gamma  \mapsto (\restr{\gamma}{\ov F_i})_{i \in I}
$$
(depending on the choice of $\mf{sp}$) and
$$
\xymatrix{
0 \ar[r]
& \pi_1^{\on{geom}}(\eta^I, \ov{\eta^I}) \ar[r]  \ar@{->>}[d]
& \pi_1(\eta^I, \ov{\eta^I}) \ar[r]  \ar[d]
& \wh{\Z} \ar[r]  \ar@{^{(}->}[d]
& 0 \\
0 \ar[r]
& \pi_1^{\on{geom}}(\eta, \ov{\eta})^I \ar[r]
& \pi_1(\eta, \ov{\eta})^I \ar[r]
& \wh{\Z}^I \ar[r]
& 0
}
$$
$$
\xymatrix{
0 \ar[r]
& \pi_1^{\on{geom}}(\eta^I, \ov{\eta^I}) \ar[r]  \ar@{->>}[d]
& \on{Weil}(\eta^I, \ov{\eta^I}) \ar[r]  \ar[d]
& \Z \ar[r]   \ar@{^{(}->}[d]
& 0 \\
0 \ar[r]
& \pi_1^{\on{geom}}(\eta, \ov{\eta})^I \ar[r]
& \on{Weil}(\eta, \ov{\eta})^I \ar[r]
& \Z^I \ar[r]
& 0
}
$$
where $\Z^I$ has the discrete topology.

\sssec{}
We denote by $\on{Frob}: \ov{F^I} \rightarrow \ov{F^I}$ the absolute Frobenius morphism over $\Fq$. 
We have a morphism
\begin{equation}  \label{equation-Weil-to-FWeil}
\begin{aligned}
\on{Weil}(\eta^I, \ov{\eta^I}) & \rightarrow \on{FWeil}(\eta^I, \ov{\eta^I}) \\
\gamma & \mapsto \Frob^{-n} \cdot \gamma
\end{aligned}
\end{equation}
We have a surjective morphism
\begin{equation}   \label{equation-FWeil-to-Weil-I}
\begin{aligned}
\Psi: \on{FWeil}(\eta^I, \ov{\eta^I}) & \twoheadrightarrow \on{Weil}(\eta, \ov{\eta})^I \\
\delta & \mapsto \big( \Frob_{\{i\}}^{n_i} \cdot \restr{\delta}{\ov{F_i}} \big)_{i \in I}
\end{aligned}
\end{equation}
(depending on the choice of $\mf{sp}$). 
We have
\begin{equation}  \label{equation-exact-pi1geo-FWeil}
\xymatrix{
0 \ar[r]
& \pi_1^{\on{geom}}(\eta^I, \ov{\eta^I}) \ar[r]  \ar@{->>}[d]
& \on{FWeil}(\eta^I, \ov{\eta^I}) \ar[r]  \ar@{->>}[d]^{\Psi}
& \Z^I \ar[r]  \ar[d]^{\simeq}
& 0 \\
0 \ar[r]
& \pi_1^{\on{geom}}(\eta, \ov{\eta})^I \ar[r]
& \on{Weil}(\eta, \ov{\eta})^I \ar[r]
& \Z^I \ar[r]
& 0
}
\end{equation}

Let 
$Q$ be the kernel of $\Psi$ (not depending on the choice of $\mf{sp}$).
Note that $Q$ is equal to the kernel of $\pi_1^{\on{geom}}(\eta^I, \ov{\eta^I}) \rightarrow \pi_1^{\on{geom}}(\eta, \ov{\eta})^I $.

%\sssec{}   \label{subsection-FWeil-singleton}
When $I$ is a singleton, (\ref{equation-Weil-to-FWeil}) is an isomorphism.
Its inverse is the morphism (\ref{equation-FWeil-to-Weil-I}).

\begin{prop} [\emph{cf.} Proposition 6.1 in \cite{drinfeld-compact}] \label{prop-Drinfeld-kernel-is-open-subgp-and-normal}
The kernel $Q$ is equal to the intersection of all open subgroups of 
$\pi_1^{\on{geom}}(\eta^I, \ov{\eta^I})$, 
which are normal in $\on{FWeil(\eta^I, \ov{\eta^I})}$.

\cqfd
\end{prop}

Proposition \ref{prop-Drinfeld-kernel-is-open-subgp-and-normal} has the following consequence.

\begin{lem} \label{lem-Drinfeld-partial-Frob-to-pi-I-finite-set}
A continuous action of $\on{FWeil(\eta^I, \ov{\eta^I})}$ on a finite set factors through $\on{Weil}(\eta, \ov{\eta})^I$.
\end{lem}
\dem
We denote by $S$ a finite set and by $\rho: \on{FWeil}(\eta^I, \ov{\eta^I}) \rightarrow \Aut(S)$ a continuous homomorphism. Then:
\begin{itemize}
\item $\on{Ker}(\restr{\rho}{\pi_1^{\on{geom}}(\eta^I, \ov{\eta^I})}) = \pi_1^{\on{geom}}(\eta^I, \ov{\eta^I}) \cap \on{Ker}(\rho)$ is normal in $\on{FWeil(\eta^I, \ov{\eta^I})}$.
\item $\on{Ker}(\restr{\rho}{\pi_1^{\on{geom}}(\eta^I, \ov{\eta^I})})$ is an open subgroup of $\pi_1^{\on{geom}}(\eta^I, \ov{\eta^I})$, because $1$ is open in $\Aut(S)$ which is a finite group. 
\end{itemize} 

Proposition \ref{prop-Drinfeld-kernel-is-open-subgp-and-normal} implies that $\on{Ker}(\restr{\rho}{\pi_1^{\on{geom}}(\eta^I, \ov{\eta^I})})$ contains $Q$. 
\cqfd

\begin{rem}
In fact, Lemma \ref{lem-Drinfeld-partial-Frob-to-pi-I-finite-set} is equivalent to Proposition \ref{prop-Drinfeld-kernel-is-open-subgp-and-normal}. Here is how Lemma \ref{lem-Drinfeld-partial-Frob-to-pi-I-finite-set} implies Proposition \ref{prop-Drinfeld-kernel-is-open-subgp-and-normal}:

By Lemma \ref{lem-Drinfeld-partial-Frob-to-pi-I-finite-set}, the profinite completion $\on{FWeil}(\eta^I, \ov{\eta^I})^{\wedge}$ of $\on{FWeil}(\eta^I, \ov{\eta^I})$ is equal to the profinite completion $(  \on{Weil}(\eta, \ov{\eta})^I   )^{\wedge}$ of $\on{Weil}(\eta, \ov{\eta})^I$. Thus 
$$
\begin{aligned}
Q & = \Ker(\on{FWeil}(\eta^I, \ov{\eta^I}) \rightarrow (  \on{Weil}(\eta, \ov{\eta})^I   )^{\wedge}) \\
& = \Ker(\on{FWeil}(\eta^I, \ov{\eta^I}) \rightarrow \on{FWeil}(\eta^I, \ov{\eta^I})^{\wedge}) \\
& = \bigcap_{\mc V \text{ open subgroup of finite index in } \on{FWeil}(\eta^I, \ov{\eta^I})} \mc V
\end{aligned}
$$
Moreover, on the one hand, it is evident that
$$\bigcap_{\mc V \text{ open subgroup of finite index in } \on{FWeil}(\eta^I, \ov{\eta^I})} \mc V \supset \bigcap_{\mc U \text{ open in } \pi_1^{\on{geom}}(\eta^I, \ov{\eta^I}) \text{ normal in } \on{FWeil}(\eta^I, \ov{\eta^I})} \mc U$$
On the other hand, for any open subgroup $\mc U$ of $\pi_1^{\on{geom}}(\eta^I, \ov{\eta^I})$ which is normal in $\on{FWeil(\eta^I, \ov{\eta^I})}$, the quotient group $\on{FWeil(\eta^I, \ov{\eta^I})} / \mc U$ is an extension of $\Z^I$ by $\pi_1^{\on{geom}}(\eta^I, \ov{\eta^I}) / \mc U$. Since an extension of $\Z^I$ by a finite group injects in its profinite completion, we have
$$\bigcap_{\mc V \text{ open subgroup of finite index in } \on{FWeil}(\eta^I, \ov{\eta^I}), \mc V \supset \mc U} \mc V / \mc U=1.$$
We deduce that $$\bigcap_{\mc V \text{ open subgroup of finite index in } \on{FWeil}(\eta^I, \ov{\eta^I})} \mc V \subset \bigcap_{\mc U \text{ open in } \pi_1^{\on{geom}}(\eta^I, \ov{\eta^I}) \text{ normal in } \on{FWeil}(\eta^I, \ov{\eta^I})} \mc U.$$
\end{rem}

\begin{rem}
The general form of Drinfeld's lemma is proved in \cite[lem.~8.11]{vincent}. 
Note that the equivalence of categories 
$$\{ \text{ étale coverings of } \eta^I \} \isom \{ \text{ finite sets equipped with a continuous action of } \pi_1(\eta^I, \ov{\eta^I})  \}$$
induces an equivalence of categories
$$
\left\{\begin{array}{c}
\text{ étale coverings of } \eta^I \text{ equipped with}\\
\text{the action of partial Frobenius morphisms}
\end{array}\right\}\isom 
\left\{\begin{array}{c}
\text{ finite sets equipped with}\\
\text{a continuous action of } \on{FWeil}(\eta^I, \ov{\eta^I})
\end{array}\right\}.
$$
\cite[lem.~8.11]{vincent} immediately implies Lemma \ref{lem-Drinfeld-partial-Frob-to-pi-I-finite-set}, thus Proposition \ref{prop-Drinfeld-kernel-is-open-subgp-and-normal}.
\end{rem}

\begin{lem} [Drinfeld]  \label{lem-Drinfeld-partial-Frob-to-pi-I-coef-integral}
A continuous action of $\on{FWeil(\eta^I, \ov{\eta^I})}$ on a $\mc O_E$-module of finite type factors through $\on{Weil}(\eta, \ov{\eta})^I$.
\end{lem}

\dem
Let $\rho: \on{FWeil}(\eta^I, \ov{\eta^I}) \rightarrow GL_r(\mc O_E)$ be a continuous morphism. 
Let $\lambda_E$ be a uniformizer of $\mc O_E$.
For any $n >1$, the induced morphism $\rho_n: \on{FWeil}(\eta^I, \ov{\eta^I}) \rightarrow GL_r(\mc O_E / {\lambda_E}^n \mc O_E)$ is continuous. 
By Lemma \ref{lem-Drinfeld-partial-Frob-to-pi-I-finite-set}, $\on{Ker}(\restr{\rho_n}{\pi_1^{\on{geom}}(\eta^I, \ov{\eta^I})})$ contains $Q$. 
Taking the limit on $n$ we get the lemma.
\cqfd

\quad

\begin{rem}
In the following, Lemma \ref{lem-Drinfeld-partial-Frob-to-pi-I-coef-rational}, Lemma \ref{lem-Drinfeld-kernel-is-closed-subgp-topo-f-t-and normal} and their proofs are due to Drinfeld and were communicated to the author by V. Lafforgue. 
\end{rem}

\begin{lem} [Drinfeld \cite{drinfeld-lemma}]  \label{lem-Drinfeld-partial-Frob-to-pi-I-coef-rational}
A continuous action of $\on{FWeil(\eta^I, \ov{\eta^I})}$ on an $E$-vector space of finite dimension factors through $\on{Weil}(\eta, \ov{\eta})^I$.
\end{lem}
\dem
Let $\rho: \on{FWeil}(\eta^I, \ov{\eta^I}) \rightarrow GL_r(\Ql)$ be a continuous morphism. In this case 
\begin{itemize}
\item $\on{Ker}(\restr{\rho}{\pi_1^{\on{geom}}(\eta^I, \ov{\eta^I})})$ is normal in $\on{FWeil(\eta^I, \ov{\eta^I})}$.
\item $\on{Ker}(\restr{\rho}{\pi_1^{\on{geom}}(\eta^I, \ov{\eta^I})})$ is closed in $\pi_1^{\on{geom}}(\eta^I, \ov{\eta^I})$ (we do not know if it is open). 
\item $\rho( \pi_1^{\on{geom}}(\eta^I, \ov{\eta^I}) )$ is topologically finitely generated (i.e. there exists a dense finitely generated subgroup).\footnote{Also called "of finite type" in some literature.}
\end{itemize} 

To prove the last statement, note that $\pi_1^{\on{geom}}(\eta^I, \ov{\eta^I}) $ is compact and the morphism $\rho$ is continuous. So $\rho( \pi_1^{\on{geom}}(\eta^I, \ov{\eta^I}) ) \subset GL_r(E) $ is a closed subgroup. 
By the theory of $\ell$-adic Lie groups \footnote{Also called "analytical $\ell$-adic groups"} in \cite{lazard} (recalled in \cite{serre}) or \cite{schneider}, $GL_r(E)$ is an $\ell$-adic Lie group.
By \cite[\S~1.3]{serre}, any closed subgroup of an $\ell$-adic Lie group is still an $\ell$-adic Lie group.
By \cite[th.~2]{serre}, any $\ell$-adic Lie group is topologically finitely generated.
We deduce that $\rho( \pi_1^{\on{geom}}(\eta^I, \ov{\eta^I}) )$ is topologically finitely generated. Then apply Lemma \ref{lem-Drinfeld-kernel-is-closed-subgp-topo-f-t-and normal} below to $H=\on{Ker}(\restr{\rho}{\pi_1^{\on{geom}}(\eta^I, \ov{\eta^I})})$. We deduce that $\on{Ker}(\restr{\rho}{\pi_1^{\on{geom}}(\eta^I, \ov{\eta^I})})$ contains $Q$.
\cqfd

\begin{lem} [Drinfeld \cite{drinfeld-lemma}]  \label{lem-Drinfeld-kernel-is-closed-subgp-topo-f-t-and normal}
Let $H$ be a closed subgroup of $\pi_1^{\on{geom}}(\eta^I, \ov{\eta^I})$ such that
\begin{itemize}
\item $H$ is normal in $\on{FWeil}(\eta^I, \ov{\eta^I}) $;
\item $\pi_1^{\on{geom}}(\eta^I, \ov{\eta^I}) / H$ is topologically finitely generated.
\end{itemize} 
Then $H$ contains $Q$.
\end{lem}
\dem
For any finite group $L$, we define $$U_L:= \bigcap_{f \text{ trivial on } H} \on{Ker}\left(  \pi_1^{\on{geom}}(\eta^I, \ov{\eta^I})  \xrightarrow{f} L  \right) .$$ Since $\pi_1^{\on{geom}}(\eta^I, \ov{\eta^I}) / H$ is topologically finitely generated, there is only a finite number of morphisms $f$ in the right hand side. Thus $U_L$ is an open subgroup of $\pi_1^{\on{geom}}(\eta^I, \ov{\eta^I}) $. 
Moreover, by definition, $U_L$ is preserved by any automorphism of $\pi_1^{\on{geom}}(\eta^I, \ov{\eta^I}) $ which preserves $H$. By hypothesis $H$ is normal in $\on{FWeil}(\eta^I, \ov{\eta^I}) $. We deduce that $U_L$ is normal in $\on{FWeil}(\eta^I, \ov{\eta^I}) $. Then by Proposition \ref{prop-Drinfeld-kernel-is-open-subgp-and-normal}, $U_L$ contains $Q$.

Since $\pi_1^{\on{geom}}(\eta^I, \ov{\eta^I})$ is profinite and $H$ is a closed and normal subgroup, the quotient group $\pi_1^{\on{geom}}(\eta^I, \ov{\eta^I}) / H$ is profinite. So $\pi_1^{\on{geom}}(\eta^I, \ov{\eta^I}) / H = \underset{\longleftarrow}{\on{lim}} \; \pi_1^{\on{geom}}(\eta^I, \ov{\eta^I}) / U_L$. Thus $H=\bigcap_{L \text{ finite group }} U_L$, which contains $Q$.
\cqfd

\sssec{}  \label{subsection-A-and-M}
Let $A$ be a finitely generated commutative $E$-algebra. 
%with trivial Jacobson radical. 
Let $M$ be an $A$-module of finite type. 

By (\ref{equation-exact-pi1geo-FWeil}), we have $0 \rightarrow \pi_1^{\on{geom}}(\eta^I, \ov{\eta^I}) \rightarrow \on{FWeil(\eta^I, \ov{\eta^I})} \rightarrow \Z^I \rightarrow 0$.
An action of $\on{FWeil(\eta^I, \ov{\eta^I})}$ on $M$ is said to be continuous if $M$ is a union of finite dimensional $E$-vector subspaces which are stable under $\pi_1^{\on{geom}}(\eta^I, \ov{\eta^I})$ and on which the action of $\pi_1^{\on{geom}}(\eta^I, \ov{\eta^I})$ is continuous.

This implies that for any finite codimensional ideal $I$ of $A$, the action of $\on{FWeil(\eta^I, \ov{\eta^I})}$ on $M / IM$ is continuous.

\begin{lem}   \label{lem-Drinfeld-partial-Frob-to-pi-I-module-of-finite-type}
Let $A$ and $M$ as in \ref{subsection-A-and-M}.
A continuous $A$-linear action of $\on{FWeil(\eta^I, \ov{\eta^I})}$ on $M$ factors through $\on{Weil}(\eta, \ov{\eta})^I$.
\end{lem}

\dem
For any maximal ideal $\mf m$ of $A$, since $A$ is finitely generated over $E$ (in particular $A$ is Noetherian), %by Hilbert Nullstellensatz $A / \mf m$ is of finite dimension over $E$. Thus 
for any $n \in \N$, the quotient $A / \mf m^n$ is of finite dimension over $E$.
%${\mf m}^n$ is an ideal of $A$ of finite codimension. 
Since $M$ is an $A$-module of finite type, $M / \mf m^n M$ is an $A / \mf m^n$-module of finite type. Thus $M / \mf m^n M$ is an $E$-vector space of finite dimension.

Applying Lemma \ref{lem-Drinfeld-partial-Frob-to-pi-I-coef-rational} to $M / \mf m^n M$, we deduce that the action of $Q$ on $M / \mf m^n M$ is trivial. 
Since $A$ is Noetherian, for any $q \in Q$ and $x \in M$, we have 
$$q \cdot x - x \in \bigcap_{\mf m \text{ max ideal }} \, \bigcap_{n=1}^{\infty} \mf m^n M \, \overset{(a)}\subset \bigcap_{\mf m \text{ max ideal }} \on{Ker} (M \rightarrow M_{\mf m}) \overset{(b)}= \{ 0 \} .$$ 
where $M_{\mf m}$ is the localization of $M$ on $A-\mf m$. (a) follows from \cite[Theorem~8.9]{matsumura}  and (b) follows from \cite[Theorem~4.6]{matsumura}. 
We deduce that $q \cdot x=x$. Thus the action of $Q$ on $M$ is trivial. 
\cqfd

\quad

\sssec{}     \label{subsection-action-of-FWeil-eta-I-eta-I-bar}
By the discussion after remarque 8.18 of \cite{vincent}, we have a continuous action of $\on{FWeil(\eta^I, \ov{\eta^I})}$ on $\restr{ \mc H_{I, W}^j }{ \ov{\eta^I} }$ (depending on the choice of $\ov{\eta^I}$ and $\mf{sp}$) which combines the action of $\pi_1(\eta^I, \ov{\eta^I})$ and the action of the partial Frobenius morphisms. 

Concretely, let $\theta \in \on{FWeil(\eta^I, \ov{\eta^I})}$ such that $\restr{\theta}{(F^I)^{\on{perf}}} = \prod_{i \in I} \on{Frob}_{\{i\}}^{-d_i}$. It induces a specialization map (which is in fact a morphism): 
$$\mf{sp}_{\theta}: \prod_{i \in I} \on{Frob}_{\{i\}}^{-d_i} (  \ov{\eta^I}  ) \rightarrow  \ov{\eta^I}.$$
The action of $\theta$ on $\restr{ \mc H_{I, W}^j }{ \ov{\eta^I} }$ is defined to be the composition:
$$\restr{ \mc H_{I, W}^j }{ \ov{\eta^I} } \xrightarrow{ \mf{sp}_{\theta}^* } \restr{ \mc H_{I, W}^j }{ \prod_{i \in I} \on{Frob}_{\{i\}}^{-d_i} (  \ov{\eta^I}  ) } \xrightarrow{ \prod_{i \in I} \on{F}_{\{i\}}^{-d_i} } \restr{ \mc H_{I, W}^j }{ \ov{\eta^I} },$$
where the second map is the partial Frobenius morphism on $\mc H_{I, W}^j$ defined in \cite[\S~4.3]{vincent}.

To see that this action is continuous, note that $\restr{\mc H_{I, W}^j}{\ov{\eta^I}} = \varinjlim _{\mu} \restr{\mc H_{I, W}^{j, \, \leq\mu} }{\ov{\eta^I}} $. Each $\restr{ \mc H_{I, W}^{j, \, \leq \mu} }{\ov{\eta^I}}$ is finite dimensional and stable by $\pi_1^{\on{geom}}(\eta^I, \ov{\eta^I})$; and the action of $\pi_1^{\on{geom}}(\eta^I, \ov{\eta^I})$ is continuous on $\restr{ \mc H_{I, W}^{j, \, \leq \mu} }{\ov{\eta^I}}$.

\begin{prop}   \label{prop-FWeil-on-H-G-factors-through-Weil}  
The action of $\on{FWeil(\eta^I, \ov{\eta^I})}$ on $\restr{\mc H_{I, W}^j}{\ov{\eta^I}} $ factors through $\on{Weil}(\eta, \ov{\eta})^I$.
\end{prop}
\dem
The action of $\on{FWeil(\eta^I, \ov{\eta^I})}$ on $\restr{ \mc H_{I, W}^j }{ \ov{\eta^I} }$ commutes with the action of the Hecke algebra $\ms H_{G, u}$. The Hecke algebra $\ms H_{G, u}$
%$\ms H_{G, u} \simeq \ms H_{T, u}^W$ 
is finitely generated over $E$ and is commutative.   
%(because $\ms H_{T, u}$ is a ring of polynomial).
%and has trivial Jacobson radical . 
By Theorem \ref{thm-coho-cht-is-Hecke-mod-type-fini-Hecke-en-v}, $\restr{\mc H_{I, W}^j}{\ov{\eta^I}} $ is a $\ms H_{G, u}$-module of finite type.
Applying Lemma \ref{lem-Drinfeld-partial-Frob-to-pi-I-module-of-finite-type} to $A=\ms H_{G, u}$ and $M = \restr{\mc H_{I, W}^j}{\ov{\eta^I}} $, we obtain the proposition.

\cqfd

\begin{rem}
Proposition \ref{prop-FWeil-on-H-G-factors-through-Weil} generalizes \cite[prop.~8.27]{vincent}.
\end{rem}

\subsection{More on Drinfeld's lemma}   \label{subsection-more-Drinfeld-lemmas}

We need the following variant of \cite[lem.~8.2]{vincent}.

\sssec{}
As in \cite[\S~1.1.10]{deligne}, let $Y$ be a connected scheme of finite type over $\Fq$. Let $\ov{\eta_Y}$ be a geometric generic point of $Y$. We define the Weil group $\on{Weil}(Y, \ov{\eta_Y})$ as the inverse image of $\on{Weil}(\Fqbar / \Fq) $ in $\pi_1(Y, \ov{\eta_Y})$.

A Weil $\Ql$-sheaf $\mc F$ on $Y$ consists a pair $(\ov{\mc F}, \phi)$, where $\ov{\mc F}$ is a constructible $\Ql$-sheaf on $\ov Y = Y \underset{\Fq}{\times} \ov{\Fq}$, equipped with an action of $\on{Weil}(\ov{\Fq} / \Fq)$: $\phi: (\Id_{Y} \times \Frob)^* \ov{\mc F} \isom \ov{\mc F} $. 

The category of smooth $E$-sheaves over $Y$ is a full subcategory of the category of smooth Weil $\Ql$-sheaves over $Y$.

We have an equivalence of categories
$$\{ \text{smooth Weil } \Ql \text{-sheaves over } Y \} \isom \{ \text{continuous finite dimensional } \Ql \text{-representations of } \on{Weil}(Y, \ov{\eta_Y})  \}.$$

\begin{lem}[rational coefficients version of lemme~8.2 in \cite{vincent}] \label{lem-Drinfeld-Weil-U-on-Ql-vector-space-E-faisceau}  
%(a variant of \cite{vincent} Lemma 8.11 for smooth $E$-sheaves)
Let $U$ be an open dense subscheme of $X$. We have an equivalence between 
%constructed in the proof below, 
\begin{enumerate}[label=(\alph*)]
\item the category of smooth Weil $\Ql$-sheaves over $U^I$, equipped with the partial Frobenius morphisms;
\item the category of continuous representations of $\on{Weil}(U, \ov{\eta})^I$ on $\Ql$-vector spaces of finite dimension. 
\end{enumerate}
which is characterized by the following two facts as in \cite[lem.~8.2]{vincent}:
\begin{itemize}
\item the composition with the restriction functor of the representations of $\on{Weil}(U, \ov{\eta})^I$ to the representations of $\on{Weil}(U, \ov{\eta})$ (diagonally) is the functor $\mf E \mapsto \restr{\mf E}{\Delta(\ov{\eta})}$;
\item if $(\mc F_{i})_{i \in I}$ is a family of smooth Weil $\Ql$-sheaves over $U$, then the image of the functor of $\boxtimes_{i \in I} \mc F_i$ is $\restr{(\boxtimes_{i \in I} \mc F_i)}{\Delta(\ov{\eta})} = \otimes_{i \in I} \restr{(\mc F_i)}{\ov{\eta}}$ equipped with the action of $\on{Weil}(U, \ov{\eta})^I$ coming from the fact that each $\restr{(\mc F_i)}{\ov{\eta}}$ is equipped with an action of $\on{Weil}(U, \ov{\eta})$.
\end{itemize} 
\end{lem}

\sssec{}
As in \cite[rem.~8.3]{vincent}, the functor (b) $\rightarrow$ (a) is explicit: the image of a continuous representation $\on{Weil}(U, \ov{\eta})^I \rightarrow GL(V)$ where $V$ is a $\Ql$-vector space of finite dimension, is the smooth Weil $\Ql$-sheaf $\mc F$ over $U^I$ associated to the continuous representation $\on{Weil}(U^I, \Delta(\ov{\eta})) \rightarrow \on{Weil}(U, \ov{\eta})^I \rightarrow GL(V)$. 

\quad

\noindent {\bf Proof of Lemma \ref{lem-Drinfeld-Weil-U-on-Ql-vector-space-E-faisceau}.} 
("hard" direction, the functor (a) to (b))
Let $\mf E$ be a smooth Weil $E$-sheaf over $U^I$ equipped with the partial Frobenius morphisms. It induces a continuous morphism $$\on{FWeil}(\eta^I, \ov{\eta^I}) \rightarrow Aut(\restr{\mf E}{\ov{\eta^I}}).$$ By Lemma \ref{lem-Drinfeld-partial-Frob-to-pi-I-coef-rational}, this morphism factors through $\on{Weil}(\eta, \ov{\eta})^I$. We deduce a continuous morphism: 
%Thus $\mc F$ induces a continuous morphism 
$$\on{Weil}(\eta, \ov{\eta})^I \rightarrow Aut(\restr{\mf E}{\ov{\eta^I}}).$$
Since $\mf E$ is unramified over $U^I$, the above morphism factors through $\on{Weil}(U, \ov{\eta})^I$. We deduce a continuous morphism: $$\on{Weil}(U, \ov{\eta})^I \rightarrow Aut(\restr{\mf E}{\ov{\eta^I}}).$$
Since $\mf E$ is smooth over $U^I$, the homomorphism of specialization $\mf{sp}^*: \restr{\mf E}{\Delta(\ov{\eta})} \rightarrow \restr{\mf E}{\ov{\eta^I}}$ is an isomorphism. We deduce a continuous morphism (not depending on $\ov{\eta^I}$ and $\mf{sp}$): $$\on{Weil}(U, \ov{\eta})^I \rightarrow Aut(\restr{\mf E}{\Delta(\ov{\eta})}).$$
\cqfd

We will also need the following lemma, whose proof is the same as in \cite[lem.~8.12]{vincent}.
\begin{lem}[rational coefficients version of lemme~8.12 in \cite{vincent}] \label{lem-Drinfeld-smooth-sheaf-extend-to-U-I}
Let $\Omega$ be an open dense subscheme of $X^I$. Let $\mf E$ be a smooth $E$-sheaf over $\Omega$ equipped with the partial Frobenius morphisms. Then there exists an open dense subscheme $U$ of $X$ such that $\mf E$ can be extended to a smooth $E$-sheaf over $U^I$. 
\cqfd
\end{lem}

\subsection{Excursion operators}    \label{subsection-excursion-operator}

The goal of this subsection is to construct the excursion operators in Construction \ref{constr-excursion-operator-on-C-c}. We begin with some preparations.

\begin{construction}     \label{const-Weil-I-on-H-Delta-eta}
Let $(\gamma_i)_{i \in I} \in \on{Weil}(\eta, \ov{\eta})^I$.
We construct an action of $(\gamma_i)_{i \in I}$ on $\restr{ \mc H_{I, W}^j   }{\Delta(\ov{\eta})}$ for any $j \in \Z$ as the composition of morphisms
\begin{equation}
\xymatrixrowsep{2pc}
\xymatrixcolsep{4pc}
\xymatrix{ 
    \restr{ \mc H_{I, W}^j   }{\Delta(\ov{\eta})}  \ar[r]^{\mf{sp}^{*}}_{\sim}
&   \restr{ \mc H_{I, W}^j  }{\ov{\eta^I}} \ar[d]^{(\gamma_i)_{i\in I}} \\  
  \restr{ \mc H_{I, W}^j  }{\Delta(\ov{\eta})} 
&       \restr{ \mc H_{I, W}^j  }{\ov{\eta^I}}  \ar[l]_{(\mf{sp}^{*})^{-1}}^{\sim}
    }
\end{equation}
where the isomorphism $\mf{sp}^*$ is defined in \ref{subsection-fix-specialization-from-eta-I-bar-to-delta-eta-bar} and Proposition \ref{prop-specialisation-is-bij}, the action of $\on{Weil}(\eta, \ov{\eta})^I$ on $\restr{\mc H_{I, W}^j}{\ov{\eta^I}}$ is defined in Proposition \ref{prop-FWeil-on-H-G-factors-through-Weil}.
\end{construction}  

\begin{lem}[\emph{cf.} lemme~9.4 in \cite{vincent} for the Hecke-finite part]  \label{lem-excursion-operator-independ-of-eta-I-bar-and-sp}  
For any $j \in \Z$, the action of $\on{Weil}(\eta, \ov{\eta})^I$ on $\restr{ \mc H_{I, W}^j   }{\Delta(\ov{\eta})}$ defined in Construction \ref{const-Weil-I-on-H-Delta-eta} is independent of the choice of $\ov{\eta^I}$ and $\mf {sp}$. 
\end{lem}

We need some preparations before the proof.

\sssec{}   \label{subsection-def-E-and-f}
Let $\ms I$ be an ideal of $\ms H_{G, u}$ of finite codimension. 
Let
%$\restr{ (\mc H_{I, W}^j / \mc I \cdot \mc H_{I, W}^j ) }{ \eta^I  } = $
$\restr{\mc H_{I, W}^j}{\eta^I} / \ms I \cdot \restr{\mc H_{I, W}^j}{\eta^I} $
be the quotient sheaf. It is stable under the action of the partial Frobenius morphisms. By Theorem \ref{thm-coho-cht-is-Hecke-mod-type-fini-Hecke-en-v}, $\restr{\mc H_{I, W}^j}{\ov{\eta^I}}$ is of finite type as $\ms H_{G, u}$-module, thus $\restr{\mc H_{I, W}^j}{\eta^I} / \ms I \cdot \restr{\mc H_{I, W}^j}{\eta^I} $ is a constructible (and smooth) sheaf over $\eta^I$.
By Lemma \ref{lem-Drinfeld-smooth-sheaf-extend-to-U-I} applied to this sheaf, there exists an open subscheme $U \subset X \sm N$ and a smooth $\Ql$-sheaf $\mf E$ over $U^I$ such that 
\begin{equation}   \label{equation-E-eta-I-equal-quotient-eta-I}
\restr{\mf E}{\eta^I} = \restr{\mc H_{I, W}^j}{\eta^I} / \ms I \cdot \restr{\mc H_{I, W}^j}{\eta^I} .
\end{equation}
%It induces $\restr{ \mf E  }{ \ov{\eta^I} } = \restr{  \mc H_{I, W}^j }{ \ov{\eta^I}  }  / \ms I \cdot \restr{ \mc H_{I, W}^j  }{ \ov{\eta^I} } $.

We define a natural morphism $f$
as the following composition of morphisms:  
\begin{equation}   \label{diagram-independent-of-sp-and-ov-eta-I}
\xymatrix{
\restr{  \mc H_{I, W}^j }{ \Delta(\ov{\eta})  }  / \ms I \cdot \restr{ \mc H_{I, W}^j  }{ \Delta(\ov{\eta}) }   \ar[d]^{f}   \ar[r]^{\mf{sp}^*}_{\simeq}   
& \restr{  \mc H_{I, W}^j }{ \ov{\eta^I}  }  / \ms I \cdot \restr{ \mc H_{I, W}^j  }{ \ov{\eta^I} }  \ar[d]_{\simeq}^{(\ref{equation-E-eta-I-equal-quotient-eta-I}) } \\
\restr{  \mf E }{ \Delta(\ov{\eta}) }  
& \restr{ \mf E  }{ \ov{\eta^I} }  \ar[l]^{(\mf{sp}^*)^{-1}}_{\simeq}  
}
\end{equation}
where both 
$\mf{sp}^*$ and $(\mf{sp}^*)^{-1}$ 
%horizontal lines 
are defined by the same specialization map $\mf{sp}$ fixed in \ref{subsection-fix-specialization-from-eta-I-bar-to-delta-eta-bar} (thus $f$ is independent of the choice of $\ov{\eta^I}$ and $\mf{sp}$). The isomorphism of the upper line follows from Proposition \ref{prop-specialisation-is-bij}. The isomorphism of the lower line follows from the fact that $\mf E$ is smooth.

\sssec{}  \label{subsection-H-I-W-includ-in-H-I-W-completion-m-adic}
For any maximal ideal $\mf m$ of $\ms H_{G, u}$, we denote by 
%$(\restr{  \mc H_{I, W}^j  }{ \ov{\eta^I} } )_{\mf m}^{\wedge} $ the $\mf m$-adic completion of $\restr{  \mc H_{I, W}^j  }{ \ov{\eta^I} } $ and by 
$(\restr{  \mc H_{I, W}^j  }{ \Delta(\ov{\eta}) } )_{\mf m}^{\wedge} $ the $\mf m$-adic completion of $ \restr{  \mc H_{I, W}^j  }{ \Delta(\ov{\eta}) } $. 
By Theorem \ref{thm-coho-cht-is-Hecke-mod-type-fini-Hecke-en-v}, $\restr{  \mc H_{I, W}^j  }{ \Delta(\ov{\eta}) }$ is of finite type as $\ms H_{G, u}$-module. By \cite[Theorem~4.6]{matsumura}, the morphism 
\begin{equation}   \label{equation-H-includ-in-prod-H-m-completion}
\restr{  \mc H_{I, W}^j  }{ \Delta(\ov{\eta}) }  \rightarrow \prod_{\mf m} (\restr{  \mc H_{I, W}^j  }{ \Delta(\ov{\eta}) } )_{\mf m}^{\wedge}
\end{equation} is injective.

\noindent {\bf Proof of Lemma \ref{lem-excursion-operator-independ-of-eta-I-bar-and-sp}.}
The proof consists of 2 steps.

\noindent\textbf{Step~1.} Let $\ms I$ and $f$ defined as in \ref{subsection-def-E-and-f}.
The morphism $f$ is compatible with the action of $\on{Weil}(\eta, \ov{\eta})^I$, where 
\begin{itemize}
\item the actions on $\restr{  \mc H_{I, W}^j }{ \ov{\eta^I}  }  / \ms I \cdot \restr{ \mc H_{I, W}^j  }{ \ov{\eta^I} }$ and on $\restr{ \mf E  }{ \ov{\eta^I} } $ are given by Lemma \ref{lem-Drinfeld-partial-Frob-to-pi-I-coef-rational};
\item the action on $\restr{  \mc H_{I, W}^j }{ \Delta(\ov{\eta})  }  / \ms I \cdot \restr{ \mc H_{I, W}^j  }{ \Delta(\ov{\eta}) }$ is induced by $\mf{sp}^*$ (defined in Construction \ref{const-Weil-I-on-H-Delta-eta});
\item the action on $\restr{  \mf E }{ \Delta(\ov{\eta}) } $ is given by Lemma \ref{lem-Drinfeld-Weil-U-on-Ql-vector-space-E-faisceau}, hence is independent of the choice of $\ov{\eta^I}$ and $\mf{sp}$.
\end{itemize} 
By \ref{subsection-def-E-and-f}, $f$ is independent of the choice of $\ov{\eta^I}$ and $\mf{sp}$.
%The action of $\on{Weil}(\eta, \ov{\eta})^I$ on $\restr{  \mf E }{ \Delta(\ov{\eta}) } $ given by Lemma \ref{lem-Drinfeld-Weil-U-on-Ql-vector-space-E-faisceau} is independent of the choice of $\ov{\eta^I}$ and $\mf{sp}$. 
We deduce that the action of $\on{Weil}(\eta, \ov{\eta})^I$ on $\restr{  \mc H_{I, W}^j }{ \Delta(\ov{\eta})  }  / \mc I \cdot \restr{ \mc H_{I, W}^j  }{ \Delta(\ov{\eta}) } $ is independent of the choice of $\ov{\eta^I}$ and $\mf{sp}$. 

\quad

\noindent\textbf{Step~2.} Through the morphism (\ref{equation-H-includ-in-prod-H-m-completion}),
the action of $\on{Weil}(\eta, \ov{\eta})^I$ on $\restr{  \mc H_{I, W}^j  }{ \Delta(\ov{\eta}) }$ is compatible with the action on $\prod_{\mf m} (\restr{  \mc H_{I, W}^j  }{ \Delta(\ov{\eta}) } )_{\mf m}^{\wedge}$. 

Applying Step~1 to $\ms I = \mf m^n$, $n \in \N$, we deduce that the action of $\on{Weil}(\eta, \ov{\eta})^I$ on $\prod_{\mf m} (\restr{  \mc H_{I, W}^j  }{ \Delta(\ov{\eta}) } )_{\mf m}^{\wedge} $ is independent of the choice of $\ov{\eta^I}$ and $\mf{sp}$. So is the action of $\on{Weil}(\eta, \ov{\eta})^I$ on $\restr{  \mc H_{I, W}^j  }{ \Delta(\ov{\eta}) } $.
\cqfd

\sssec{}   \label{subsection-creation-annilation}
Let $I$ be a finite set and $W$ be a representation of $\wh G^I$. Let $x \in W$ and $\xi \in W^*$ be invariant under the diagonal action of $\wh G$.
We denote by $E_{(X \sm N)}$ the constant sheaf over $X \sm N$.
In \cite[d\'efi.~5.1-5.2]{vincent} and the beginning of Section 9 of \emph{loc. cit.}, V. Lafforgue defined the creation operator 
\begin{equation}    \label{equation-creation-operator-over-X}
\mc C^{\sharp}_{x}: C_c(\Bun_{G, N}(\Fq) / \Xi, E) \boxtimes E_{(X \sm N)}  \rightarrow \restr{\mc H_{I, W}^0 }{\Delta(X \sm N)} ,
\end{equation}
and 
the annihilation operator 
\begin{equation}    \label{equation-annihilation-operator-over-X}
\mc C^{\flat}_{\xi}: \restr{\mc H_{I, W}^0 }{\Delta(X \sm N)} \rightarrow C_c(\Bun_{G, N}(\Fq) / \Xi, E) \boxtimes E_{(X \sm N)}  .
\end{equation}

We denote by $\restr{ \mc C^{\sharp}_{x}  }{ \ov{\eta} }$ (resp. $\restr{ \mc C^{\flat}_{\xi}  }{ \ov{\eta} }$) the restriction of (\ref{equation-creation-operator-over-X}) (resp. (\ref{equation-annihilation-operator-over-X})) on $\ov{\eta}$.

\begin{construction}    \label{constr-excursion-operator-on-C-c}
Let $I$, $W$, $x$ and $\xi$ as in \ref{subsection-creation-annilation}. Let $(\gamma_i)_{i \in I} \in \on{Weil}(\eta, \ov{\eta})^I$. We construct an excursion operator $S_{I, W, x, \xi, (\gamma_i)_{i \in I}} $ on $C_c(\Bun_{G, N}(\Fq) / \Xi, E)$ as the composition of morphisms:
$$
\xymatrix{ 
    C_c(\Bun_{G, N}(\Fq) / \Xi, E)   \ar[r]^{ \quad \quad \restr{ \mc C^{\sharp}_{x}  }{ \ov{\eta} } } 
&    \restr{ \mc H_{I, W}^{0}   }{\Delta(\ov{\eta})}   \ar[d]^{(\gamma_i)_{i\in I}} \\  
C_c(\Bun_{G, N}(\Fq) / \Xi, E)  
&  \restr{ \mc H_{I, W}^{0}  }{\Delta(\ov{\eta})}  \ar[l]_{ \quad \quad \restr{ \mc C^{\flat}_{\xi}  }{ \ov{\eta} }}
    }
$$
where the action of $(\gamma_i)_{i\in I}$ is defined in Construction \ref{const-Weil-I-on-H-Delta-eta}.

The action of $S_{I, W, x, \xi, (\gamma_i)_{i \in I}} $ commutes with the action of the global Hecke algebra $\ms H_G = C_{c}(K_N\backslash G(\mb A)/K_N, E)$. Thus 
$$S_{I, W, x, \xi, (\gamma_i)_{i \in I}} \in \on{End}_{\ms H_G} (C_c(\Bun_{G, N}(\Fq) / \Xi, E)) .$$
\end{construction}

\begin{rem}     \label{rem-excursion-op-restrict-to-cusp}
Let $C_c^{\on{cusp}}(\Bun_{G, N}(\Fq) / \Xi, E)$ be the vector subspace of cuspidal automorphic forms. It is well-known that it is of finite dimension.
The restriction of $S_{I, W, x, \xi, (\gamma_i)_{i \in I}}$ on $C_c^{\on{cusp}}(\Bun_{G, N}(\Fq) / \Xi, E)$ coincides with the excursion operator defined in \cite[d\'efi.-prop.~9.1]{vincent}.
\end{rem}

\begin{rem}
In \cite[d\'efi.~9.3]{vincent}, the notation $H_{I, W}$ is used for the Hecke-finite part of the cohomology group in degree $0$, which we denote in this paper by %$(   \varinjlim _{\mu} \restr{  \mc  H_{I, W}^{0, \, \leq\mu} }{ \Delta(\ov{\eta}) }    )^{\on{Hf}}$. 
$( \restr{  \mc  H_{I, W}^{0} }{ \Delta(\ov{\eta}) }    )^{\on{Hf}}$. 
In this paper, the notation $H_{I, W}^j$ is used for the cohomology group in degree $j$, which is also denoted by $ \restr{  \mc  H_{I, W}^{j} }{ \Delta(\ov{\eta}) } $.
\end{rem}

\quad

Now we extend the properties satisfied by the excursion operators acting on $C_c^{\on{cusp}}(\Bun_{G, N}(\Fq) / \Xi, E)$ (proven in \cite[Sections~9 and~10]{vincent}) to the excursion operators acting on $C_c(\Bun_{G, N}(\Fq) / \Xi, E)$.

%(for the cuspidal part, see \cite{vincent} Lemme 10.1)
\begin{lem}  \label{lem-S-I-f-gamma-commutes}
The properties in \cite[lem.~10.1]{vincent} for the excursion operators acting on $C_c^{\on{cusp}}(\Bun_{G, N}(\Fq) / \Xi, E)$ extend to the same properties for the excursion operators acting on $C_c(\Bun_{G, N}(\Fq) / \Xi, E)$ with the same proof. 
In particular, the excursion operators commute with each other.
\end{lem}

\begin{lem} [\emph{cf.} lemme~10.6  in \cite{vincent} for the cuspidal part]  \label{lem-excursion-op-W-x-xi-equal-f}
The excursion operator $S_{I, W, x, \xi, (\gamma_i)_{i \in I}} $ depends only on $I, f$ and $(\gamma_i)_{i \in I}$, where $f \in \mc O(\wh G \backslash \wh G^I / \wh G )$ is the function given by
$f: (g_i)_{i \in I} \mapsto \langle \xi , (g_i)_{i \in I} \cdot x \rangle .$
\end{lem}
The proof is the same as in \cite[lemme~10.6]{vincent}.

\quad

\sssec{}    \label{subsection-fix-sp-v}
For any place $v$ of $X$, fix an algebraic closure $\ov{F_v}$ of $F_v$ and fix the embeddings such that the following diagram commutes: 
\begin{equation}
\xymatrix{
\ov{F} \ar@{^{(}->}[r]
& \ov{F_v} \\
F  \ar@{^{(}->}[r]   \ar@{^{(}->}[u]
& F_v .  \ar@{^{(}->}[u]
}
\end{equation}

Let $k_v$ be the residue field of $F_v$ and let $\ov{k_v}$ be the residue field of the maximal unramified extension of $F_v$ in $\ov{F_v}$. Let $\ov{v} = \on{Spec}(\ov{k_v})$ be the associated geometric point over $v = \on{Spec} (k_v)$. 
Let 
\begin{equation}   \label{equation-sp-v}
\mf{sp}_v: \ov{\eta} \rightarrow \ov{v}
\end{equation}
be the specialization map associated to $\ov{F} \subset \ov{F_v}$. We denote still by $\mf{sp}_v$ the image by $\Delta$ of the above specialization map (\ref{equation-sp-v}) 
\begin{equation}   \label{equation-sp-v-Delta}
\mf{sp}_v: \Delta(\ov{\eta}) \rightarrow \Delta(\ov{v}).
\end{equation}

The inclusion $\ov{F} \subset \ov{F_v}$ induces $\on{Weil}(\ov{F_v} / F_v) \subset \on{Weil}(\ov F / F)$.

\begin{lem} [\emph{cf.} lemme~10.4 in \cite{vincent} for the cuspidal part]   \label{lem-Frob-partiel-commute-with-gamma-excursion-op}
Let $v$ be a place in $X \sm N$.
Let us consider $(\gamma_i)_{i \in I} \in  \on{Weil}(\ov{F_v} / F_v)^I \subset \on{Weil}(\ov{F} / F)^I$. Let $d_i = \on{deg} ( \gamma_i )$.
We have a commutative diagram
%where the commutativity of the middle square is the lemma \ref{lem-Frob-partiel-commute-with-gamma}.
$$
\xymatrixrowsep{3pc}
\xymatrixcolsep{3pc}
\xymatrix{
C_c(\Bun_{G, N}(\Fq) / \Xi, E)   \ar[d]_{\restr{ \mc C^{\sharp}_{x}  }{ \ov{v} }}  \ar[dr]^{\restr{ \mc C^{\sharp}_{x}  }{ \ov{\eta} }}  \\
 \restr{ \mc H_{I, W}^0 }{ \Delta(\ov{v}) }    \ar[r]^{\mf{sp}_v^*}   \ar[d]_{ \prod_{i \in I} \on{F}_{\{i\}}^{\on{deg}(v)d_i} } 
& \restr{ \mc H_{I, W}^0 }{ \Delta(\ov{\eta}) }    \ar[d]^{(\gamma_i)_{i \in I}}  \\
 \restr{ \mc H_{I, W}^0 }{ \Delta(\ov{v}) }    \ar[r]^{\mf{sp}_v^*}  \ar[d]_{\restr{ \mc C^{\flat}_{\xi}  }{ \ov{v} }}
& \restr{ \mc H_{I, W}^0 }{ \Delta(\ov{\eta}) }    \ar[dl]^{\restr{ \mc C^{\flat}_{\xi}  }{ \ov{\eta} }}    \\
C_c(\Bun_{G, N}(\Fq) / \Xi, E) 
}
$$
where $\prod_{i \in I} \on{F}_{\{i\}}^{\on{deg}(v)d_i} $ is the partial Frobenius morphism, $\restr{ \mc C^{\sharp}_{x}  }{ \ov{v} }$ (resp. $\restr{ \mc C^{\flat}_{\xi}  }{ \ov{v} }$) is the restriction of (\ref{equation-creation-operator-over-X}) (resp. (\ref{equation-annihilation-operator-over-X})) on $\ov{v}$.
\end{lem}

The proof will be given in the next subsection.

\sssec{}
Let $V$ be a representation of $\wh G$ and $V^*$ be the dual of $V$. Let $\delta_V: {\bf 1} \rightarrow V \otimes V^*$ and $ev_V: V \otimes V^* \rightarrow \bf 1$ be the canonical morphisms. 

Let $v \in | X \sm N |$. Let $h_{V, v} \in \ms H_{G, v}$ be the spherical function associated to $V$ by the Satake isomorphism. Let $T(h_{V, v})$ be the Hecke operator on $C_c(\Bun_{G, N}(\Fq) / \Xi, E)$ associated to $h_{V, v}$.

\begin{lem}[\emph{cf.} lemme~10.2 in \cite{vincent} for the cuspidal part] \label{lem-S-V-v-equal-T-V-v}
Let $d \in \N$ and $\gamma \in \on{Weil}(\ov{F_v} / F_v) \subset \on{Weil}(\ov F / F)$ such that $\on{deg}(\gamma) = d$. Then $S_{\{1, 2\}, V\boxtimes V^*, \delta_V, ev_V, (\gamma, 1)}$ depends only on $d$, and if $d=1$ it equals to $T(h_{V, v})$.
\end{lem}
\dem
This follows from Lemma \ref{lem-Frob-partiel-commute-with-gamma-excursion-op} above and \cite[prop.~6.2]{vincent} (where the statement is already for the whole cohomology, not only for the cuspidal part).

\cqfd

\begin{prop}[\emph{cf.} lemme~10.10 in \cite{vincent} for the cuspidal part]   \label{prop-factorise-by-weil-X-N} 
For any $I$ and $f$, $S_{I, f, (\gamma_i)_{i\in I}}$ depends only on the image of $(\gamma_i)_{i \in I}$ in $\on{Weil}(X \sm N, \ov{\eta})^I$.
\end{prop}
\dem
(The argument is the same as in \cite[lem.~10.10]{vincent}.)\\
Let $(\delta_i)_{i \in I} \in (I_v)^I \subset \on{Weil}(\ov{F_v} / F_v)^I \subset \on{Weil}(\ov{F} / F)^I$. Applying Lemma \ref{lem-Frob-partiel-commute-with-gamma-excursion-op} to $d_i=0$, we deduce that the image of $C_c(\Bun_{G, N}(\Fq) / \Xi, E)$ in $\restr{ \mc H_{I, W}^0 }{ \ov{\eta^I} }$ is invariant by $(I_v)^I$.
Thus for $(\gamma_i)_{i \in I} \in \on{Weil}(\ov F / F)^I$ and $(\delta_i)_{i \in I} \in (I_v)^I$, we have $S_{I, f, (\gamma_i)} = S_{I, f, (\delta_i \gamma_i)}$. We have this for any choice of inclusion $\ov F \subset \ov{F_v}$. 

Since $\on{Weil}(X \sm N, \ov{\eta})$ is the quotient of $\on{Weil}(\eta, \ov{\eta})$ by the subgroup generated by the $I_v$ for $v \in X \sm N$ and their conjugates, we deduce that $S_{I, f, (\gamma_i)}$ depends only on the image of $(\gamma_i)$ by $\on{Weil}(\eta, \ov{\eta})^I \rightarrow \on{Weil}(X \sm N, \ov{\eta})^I$. \cqfd

\begin{rem}
The statement \cite[prop.~8.10]{vincent} can be generalized to the excursion operators acting on $C_c(\Bun_{G, N}(\Fq) / \Xi, E)$ and the proof is the same. So we will not state them here.
\end{rem}

\subsection{Proof of Lemma \ref{lem-Frob-partiel-commute-with-gamma-excursion-op} }

Lemma \ref{lem-Frob-partiel-commute-with-gamma-excursion-op} will follow from Lemma \ref{lem-Frob-partiel-commute-with-gamma} below. 

\begin{lem}[rational coefficients version of lemme~8.15 in \cite{vincent}] \label{lem-j-*-delta-v-to-delta-eta}
Let $U$ be an open subscheme of $X$. We denote by $j^I: U^I \hookrightarrow X^I$ the inclusion.
Let $\mf E$ be a smooth $E$-sheaf over $U^I$, equipped with the partial Frobenius morphisms. 
Let $v$ be a place in $X$ as in \ref{subsection-fix-sp-v}.
Let $(\gamma_i)_{i \in I}$ and $(d_i)_{i \in I}$ as in Lemma \ref{lem-Frob-partiel-commute-with-gamma-excursion-op}.
Then the following diagram is commutative
\begin{equation}    \label{diagram-j-*-E-sp-v-to-E}
\xymatrixrowsep{2pc}
\xymatrixcolsep{3pc}
\xymatrix{
\restr{  (j^I)_* \mf E  }{ \Delta(\ov{v}) }      \ar[d]_{ \prod_{i \in I} \on{F}_{\{i\}}^{\on{deg}(v)d_i} }    \ar[r]^{  (\mf{sp}_v )^* }
& \restr{ \mf E }{ \Delta(\ov{\eta}) }  \ar[d]^{ (\gamma_i)_{i \in I}  } \\
\restr{   (j^I)_* \mf E  }{ \Delta(\ov{v}) }   \ar[r]^{ ( \mf{sp}_v  )^* }
& \restr{ \mf E }{ \Delta(\ov{\eta}) }  
}
\end{equation}
where the vertical map on the right is the action of $\on{Weil}(U, \ov{\eta})^I$ on $\restr{ \mf E }{ \Delta(\ov{\eta}) }$ given by Lemma \ref{lem-Drinfeld-Weil-U-on-Ql-vector-space-E-faisceau}. 
\end{lem}
\dem
(The argument is the same as in \cite[lem.~8.15]{vincent}.) 
It is enough to prove the lemma with $\mf E$ of the form $\boxtimes_{i \in I} \mf E_i$ (as in Lemma \ref{lem-Drinfeld-Weil-U-on-Ql-vector-space-E-faisceau}). Noting $j: U \hookrightarrow X$ the inclusion, we have $\restr{((j^I)_* \mf E)}{\Delta(\ov{v})} = \otimes_{i \in I}(\restr{j_* \mf E_i}{\ov v})$. Thus it is enough to prove the lemma in the case where $I$ is a singleton. In this case, the commutativity follows from the definition of $\on{deg}: \on{Weil}(\ov{F_v} / F_v) \rightarrow \on{Weil}(\ov{k(v)} / k(v))$ by restriction of the action on the maximal unramified extension of $F_v$ on its residue field.
\cqfd

\sssec{}
By \cite[\S~4.3]{vincent}, for any $(d_i)_{i \in I} \in \Z^I$, there exists $\kappa \in \wh\Lambda_{G^{\mr{ad}}}^+$, such that for any $\mu \in \wh\Lambda_{G^{\mr{ad}}}^+$, we have morphisms of partial Frobenius
$$\prod_{i \in I} \on{F}_{\{i\}}^{\on{deg}(v)d_i}: \restr{  \mc H_{I, W}^{j, \, \leq \mu}  }{ \Delta(\ov{v}) }  \rightarrow \restr{  \mc H_{I, W}^{j, \, \leq \mu + \kappa}  }{ \Delta(\ov{v}) }.$$

\begin{lem}   \label{lem-H-Frob-to-H-quotient-gamma-commute}
Let $\ms I$ be an ideal of $\ms H_{G, u}$ of finite codimension as in \ref{subsection-def-E-and-f}. Let $(\gamma_i)_{i \in I}$ and $(d_i)_{i \in I}$ as in Lemma \ref{lem-Frob-partiel-commute-with-gamma-excursion-op}.
For any $\mu \in \wh\Lambda_{G^{\mr{ad}}}^+$, the following diagram is commutative
$$
\xymatrixrowsep{2pc}
\xymatrixcolsep{5pc}
\xymatrix{
\restr{  \mc H_{I, W}^{j, \, \leq \mu}  }{ \Delta(\ov{v}) }      \ar[d]_{ \prod_{i \in I} \on{F}_{\{i\}}^{\on{deg}(v)d_i} }    \ar[r]^{ ( \mf{sp}_v  )^* \quad \quad }
& \restr{ \mc H_{I, W}^j }{ \Delta(\ov{\eta}) } / \ms I \cdot \restr{ \mc H_{I, W}^j }{ \Delta(\ov{\eta}) }   \ar[d]^{ (\gamma_i)_{i \in I}  } \\
\restr{  \mc H_{I, W}^{j, \, \leq \mu + \kappa}  }{ \Delta(\ov{v}) }   \ar[r]^{ ( \mf{sp}_v  )^* \quad \quad }
& \restr{ \mc H_{I, W}^j }{ \Delta(\ov{\eta}) }  / \ms I \cdot \restr{ \mc H_{I, W}^j }{ \Delta(\ov{\eta}) } 
.}
$$
\end{lem}
\dem
Let $\mf E$ and $U$ as defined in \ref{subsection-def-E-and-f}. There exists an open subscheme $\Omega$ of $U^I$ such that $\mc H_{I, W}^{j, \, \leq \mu}$ and $\mc H_{I, W}^{j, \, \leq \mu + \kappa}$ are smooth over $\Omega$. 
Let $\iota: \Omega \hookrightarrow U^I$ and $j_{\Omega}: \Omega \hookrightarrow X^I$ be the inclusions. We have $j_{\Omega} = j^I \circ \iota$.
Since both $(j^I)^* \mc H_{I, W}^{j, \, \leq \mu} $ and $\mf E$ are smooth over $\Omega$, the morphism 
\begin{equation}
\restr{  \mc H_{I, W}^{j, \, \leq \mu}  }{ \eta^I }  \rightarrow \restr{  \mc H_{I, W}^{j}  }{ \eta^I }  \rightarrow \restr{  \mc H_{I, W}^{j}  }{ \eta^I } / \ms I \cdot \restr{  \mc H_{I, W}^{j}  }{ \eta^I } = \restr{  \mf E  }{ \eta^I } 
\end{equation}
extends to a morphism
\begin{equation}  % \label{equation-iota-H-I-W-to-iota-mf-E}
\iota^* (j^I)^* \mc H_{I, W}^{j, \, \leq \mu} \rightarrow \iota^* \mf E .
\end{equation} 
By adjunction, this morphism induces a morphism (the fonctors are not derived)
\begin{equation}
\mc H_{I, W}^{j, \, \leq \mu} \rightarrow (j^I)_* \iota_* \iota^* \mf E \simeq (j^I)_* \mf E
\end{equation}
where the isomorphism follows from the fact that $\mf E$ is smooth over $U^I$ so the adjunction morphism $\mf E \rightarrow \iota_* \iota^* \mf E$ is an isomorphism. Similarly, we have a morphism
\begin{equation}
\mc H_{I, W}^{j, \, \leq \mu + \kappa} \rightarrow (j^I)_* \iota_* \iota^* \mf E \simeq (j^I)_* \mf E
\end{equation}
The following diagram is commutative
\begin{equation}    \label{diagram-H-I-W-to-j-*-E}
\xymatrix{
\restr{  \mc H_{I, W}^{j, \, \leq \mu}  }{ \Delta(\ov{v}) }      \ar[d]_{ \prod_{i \in I} \on{F}_{\{i\}}^{\on{deg}(v)d_i} }    \ar[r]
& \restr{ (j^I)_* \mf E }{ \Delta(\ov{v}) }   \ar[d]^{ \prod_{i \in I} \on{F}_{\{i\}}^{\on{deg}(v)d_i}   } \\
\restr{  \mc H_{I, W}^{j, \, \leq \mu + \kappa}  }{ \Delta(\ov{v}) }   \ar[r]
& \restr{ (j^I)_* \mf E }{ \Delta(\ov{v}) } .
}
\end{equation}

Applying Lemma \ref{lem-j-*-delta-v-to-delta-eta} to $\mf E$, we deduce that the following diagram (which is the composition of (\ref{diagram-H-I-W-to-j-*-E}) and (\ref{diagram-j-*-E-sp-v-to-E})) is commutative
\begin{equation}
\xymatrixrowsep{2pc}
\xymatrixcolsep{5pc}
\xymatrix{
\restr{  \mc H_{I, W}^{j, \, \leq \mu}  }{ \Delta(\ov{v}) }      \ar[d]_{ \prod_{i \in I} \on{F}_{\{i\}}^{\on{deg}(v)d_i} }    \ar[r]^{ ( \mf{sp}_v  )^*  }
& \restr{ \mf E }{ \Delta(\ov{\eta}) }   \ar[d]^{ (\gamma_i)_{i \in I}  } \\
\restr{  \mc H_{I, W}^{j, \, \leq \mu + \kappa}  }{ \Delta(\ov{v}) }   \ar[r]^{ ( \mf{sp}_v  )^*  }
& \restr{ \mf E }{ \Delta(\ov{\eta}) } .
}
\end{equation}
Taking into account the isomorphism $$f: \restr{  \mc H_{I, W}^j }{ \Delta(\ov{\eta})  }  / \ms I \cdot \restr{ \mc H_{I, W}^j  }{ \Delta(\ov{\eta}) } \isom \restr{  \mf E }{ \Delta(\ov{\eta}) }  $$
defined in \ref{subsection-def-E-and-f}, we deduce the lemma.

(By construction, it is easy to verify that the compositions of morphisms
$$\restr{  \mc H_{I, W}^{j, \, \leq \mu}  }{ \Delta(\ov{v}) } \xrightarrow{\mf{sp}_v^*} \restr{ \mf E }{ \Delta(\ov{\eta}) } \xrightarrow{f^{-1}} \restr{  \mc H_{I, W}^j }{ \Delta(\ov{\eta})  }  / \ms I \cdot \restr{ \mc H_{I, W}^j  }{ \Delta(\ov{\eta}) } $$
and
$$\restr{  \mc H_{I, W}^{j, \, \leq \mu}  }{ \Delta(\ov{v}) } \xrightarrow{\mf{sp}_v^*} \restr{  \mc H_{I, W}^j }{ \Delta(\ov{\eta})  } \rightarrow \restr{  \mc H_{I, W}^j }{ \Delta(\ov{\eta})  }  / \ms I \cdot \restr{ \mc H_{I, W}^j  }{ \Delta(\ov{\eta}) } $$
are the same.)
\cqfd

\begin{lem} [\emph{cf.}  lemme~10.4 in \cite{vincent} for the cuspidal part]    \label{lem-Frob-partiel-commute-with-gamma}
Let $v \in |X \sm N|$. Let $(\gamma_i)_{i \in I} \in \on{Weil}(\ov{F_v} / F_v)^I$ with degree $(d_i)_{i \in I} \in \Z^I$.
Then for any $j \in \Z$, the following diagram is commutative:
$$
\xymatrixrowsep{2pc}
\xymatrixcolsep{3pc}
\xymatrix{
\restr{  \mc H_{I, W}^j  }{ \Delta(\ov{v}) }      \ar[d]_{ \prod_{i \in I} \on{F}_{\{i\}}^{\on{deg}(v)d_i} }    \ar[r]^{ ( \mf{sp}_v )^* }
& \restr{ \mc H_{I, W}^j }{ \Delta(\ov{\eta}) }  \ar[d]^{ (\gamma_i)_{i \in I}  } \\
\restr{  \mc H_{I, W}^j  }{ \Delta(\ov{v}) }   \ar[r]^{ ( \mf{sp}_v )^* }
& \restr{ \mc H_{I, W}^j }{ \Delta(\ov{\eta}) }  .
}
$$
\end{lem}
\dem
Applying Lemma \ref{lem-H-Frob-to-H-quotient-gamma-commute} to $\ms I = \mf m^n$, $n \in \N$, we deduce that the following diagram is commutative
\begin{equation}
\xymatrixrowsep{2pc}
\xymatrixcolsep{5pc}
\xymatrix{
\restr{  \mc H_{I, W}^{j, \, \leq \mu}  }{ \Delta(\ov{v}) }      \ar[d]_{ \prod_{i \in I} \on{F}_{\{i\}}^{\on{deg}(v)d_i} }    \ar[r]^{ ( \mf{sp}_v  )^* }
&  \prod_{\mf m} (\restr{  \mc H_{I, W}^j  }{ \Delta(\ov{\eta}) } )_{\mf m}^{\wedge}     \ar[d]^{ (\gamma_i)_{i \in I}  } \\
\restr{  \mc H_{I, W}^{j, \, \leq \mu + \kappa}  }{ \Delta(\ov{v}) }   \ar[r]^{ ( \mf{sp}_v  )^* }
&  \prod_{\mf m} (\restr{  \mc H_{I, W}^j  }{ \Delta(\ov{\eta}) } )_{\mf m}^{\wedge} .
}
\end{equation}
By \ref{subsection-H-I-W-includ-in-H-I-W-completion-m-adic}, we deduce that the following diagram is commutative
\begin{equation}
\xymatrixrowsep{2pc}
\xymatrixcolsep{5pc}
\xymatrix{
\restr{  \mc H_{I, W}^{j, \, \leq \mu}  }{ \Delta(\ov{v}) }      \ar[d]_{ \prod_{i \in I} \on{F}_{\{i\}}^{\on{deg}(v)d_i} }    \ar[r]^{ ( \mf{sp}_v  )^* }
&  \restr{  \mc H_{I, W}^j  }{ \Delta(\ov{\eta}) }     \ar[d]^{ (\gamma_i)_{i \in I}  } \\
\restr{  \mc H_{I, W}^{j, \, \leq \mu + \kappa}  }{ \Delta(\ov{v}) }   \ar[r]^{ ( \mf{sp}_v  )^* }
&   \restr{  \mc H_{I, W}^j  }{ \Delta(\ov{\eta}) }  .
}
\end{equation}
Taking the limit on $\mu$, we deduce Lemma \ref{lem-Frob-partiel-commute-with-gamma}.
\cqfd

\begin{rem}
For another proof of Lemma \ref{lem-Frob-partiel-commute-with-gamma}, see the first version of this paper on the arXiv.
\end{rem}

\subsection{Langlands parametrizations}   \label{subsection-Langlands-parametrization}

\begin{defi}
We denote by ${\mc B}^{\sim}$ the sub-$E$-algebra of $$\on{End}_{\ms H_G} \big( C_c(\Bun_{G, N}(\Fq) / \Xi, E) \big)$$ generated by all the excursion operators $S_{I, f, (\gamma_i)_{i \in I}}$. It may be infinite dimensional.
By Lemma \ref{lem-S-I-f-gamma-commutes}, ${\mc B}^{\sim}$ is commutative. 
\end{defi}

By Lemma \ref{lem-S-V-v-equal-T-V-v}, ${\mc B}^{\sim}$ contains the Hecke algebras at all the places of $X \sm N$.

\sssec{}   
%Remark that the space $\on{End}_{C_c(K_N \backslash G(\mb A) / K_N, E)} (C_c(\Bun_{G, N}(\Fq) / \Xi, E))$ may have infinite dimension.
Let $u \in |X \sm N|$.
Let $\ms I$ be an ideal of $\ms H_{G, u}$ of finite codimension.
By Proposition \ref{prop-space-of-auto-form-is-Hecke-mod-type-fini}, the quotient $E$-vector space 
\begin{equation}    \label{equation-C-c-quotient-by-I}
C_c(\Bun_{G, N}(\Fq) / \Xi, E) / \ms I \cdot C_c(\Bun_{G, N}(\Fq) / \Xi, E)
\end{equation}
is of finite dimension.

\begin{rem}
In the case $\on{dim} \ms H_{G, u} / \ms I=1$, we have $\ms I = \on{Ker} \chi$ for some character $\chi: \ms H_{G, u} \rightarrow \Ql$. In this case (\ref{equation-C-c-quotient-by-I}) is the largest quotient of $C_c(\Bun_{G, N}(\Fq) / \Xi, \Qlbar) $ on which $\ms H_{G, u}$ acts by $\chi$.
\end{rem}

\sssec{}
Replacing $\restr{\mc H_{I, W}^0}{\Delta(\ov{\eta})}$ in Construction \ref{constr-excursion-operator-on-C-c} by the quotient $\restr{\mc H_{I, W}^0}{\Delta(\ov{\eta})} / \ms I \cdot \restr{\mc H_{I, W}^0}{\Delta(\ov{\eta})}$, we define the excursion operators $S_{I, f, (\gamma_i)}$ acting on (\ref{equation-C-c-quotient-by-I}). All the properties in \S~\ref{subsection-excursion-operator} are still true for $S_{I, f, (\gamma_i)}$ acting on this quotient vector space.

\begin{defi}    \label{def-algebra-excursion-B-I}
We denote by ${\mc B}_{\ms I}$ the sub-$E$-algebra of $$\on{End}_{\ms H_G} \big( C_c(\Bun_{G, N}(\Fq) / \Xi, E) / \ms I \cdot C_c(\Bun_{G, N}(\Fq) / \Xi, E) \big)$$ generated by all the excursion operators $S_{I, f, (\gamma_i)_{i \in I}}$. It is finite dimensional. By Lemma \ref{lem-S-I-f-gamma-commutes}, ${\mc B}_{\ms I}$ is commutative. 
\end{defi}

By Lemma \ref{lem-S-V-v-equal-T-V-v}, ${\mc B}_{\ms I}$ contains the Hecke algebras at all the places of $X \sm N$.

\begin{lem} [\emph{cf.} proposition~10.10 in \cite{vincent} for the cuspidal part]
The morphism $$\on{Weil}(X \sm N, \ov{\eta})^I \rightarrow {\mc B}_{\ms I}, \quad (\gamma_i)_{i \in I} \mapsto S_{I, f, (\gamma_i)}$$ is continuous, where ${\mc B}_{\ms I}$ is endowed with the $E$-adic topology.
\end{lem}
\dem
The proof is the same as \cite[prop.~10.10]{vincent}, except that we use Lemma \ref{lem-Frob-partiel-commute-with-gamma} (of this paper) instead of \cite[lem.~10.4]{vincent}. 
\cqfd

Then we use the same arguments as in \cite{vincent} Section 11, except that we replace $C_c^{\on{cusp}}(\Bun_{G, N}(\Fq) / \Xi, \Qlbar) $ by $C_c(\Bun_{G, N}(\Fq) / \Xi, \Qlbar) / \ms I \cdot C_c(\Bun_{G, N}(\Fq) / \Xi, \Qlbar) $ and replace $\pi_1(X \sm N, \ov{\eta})$ by $\on{Weil}(X \sm N, \ov{\eta})$. We obtain:
\begin{thm}    \label{thm-decomposition-of-C-c-quotient-I-by-param-Langlands}
We have a canonical decomposition of $C_c(K_N \backslash G(\mb A) / K_N, \Qlbar)$-modules
$$C_c(\Bun_{G, N}(\Fq) / \Xi, \Qlbar) / \ms I \cdot C_c(\Bun_{G, N}(\Fq) / \Xi, \Qlbar) = \oplus_{\sigma} \mf H_{\sigma}$$
where the direct sum is indexed by $\wh G(\Qlbar)$-conjugacy classes of morphisms $\sigma: \on{Weil}(\ov F / F) \rightarrow \wh G(\Qlbar)$ defined over a finite extension of $\mathbb{Q}_{\ell}$, continuous, semisimple and unramified outside $N$.

This decomposition is characterized by the following property: $\mf H_{\sigma}$ is equal to the generalized eigenspace $\mf H_{\nu}$ associated to the character $\nu$ of $\mc B_{\ms I}$ defined by $\nu(S_{I, f, (\gamma_i)_{i \in I}}) = f(  (\sigma(\gamma_i) )_{i \in I} ).$

It is compatible with the Satake isomorphism at every place $v$ of $X \sm N$: for any irreducible representation $V$ of $\wh G$, we denote by $T(h_{V, v})$ the Hecke operator at $v$ associated to $V$. 
Then $\mf H_{\sigma}$ is included in the generalized eigenspace of $T(h_{V, v})$ for the eigenvalue $\nu (T(h_{V, v})) =\chi_V(\sigma(\on{Frob}_v))$, where $\chi_V$ is the character of $V$ and $\on{Frob}_v$ is an arbitary lifting of the Frobenius element on $v$.
\end{thm}

\begin{rem}
Contrary to the case of the cuspidal part $C_c^{\on{cusp}}(\Bun_{G, N}(\Fq) / \Xi, \Qlbar)$ in \cite{vincent}, here the actions of the Hecke operators on $C_c(\Bun_{G, N}(\Fq) / \Xi, \Qlbar) / \ms I \cdot C_c(\Bun_{G, N}(\Fq) / \Xi, \Qlbar) $ are not always diagonalizable. 
\end{rem}

\subsection{Excursion operators on cohomology groups}   \label{subsection-excursion-operator-on-cohomology}

\sssec{}   \label{subsection-notation-I-W-J-V}
Let $J$ be a finite set and $V$ be a representation of $\wh G^J$. Let $j \in \Z$. Applying Definition \ref{def-H-G-I-W-sheaf} to $J$ and $V$, we define $\mc  H_{J, V}^j $, which is an inductive limit of constructible sheaves on $(X \sm N)^J$. 

Let $I$ be a finite set and $W$ be a representation of $\wh G^I$. Let $x \in W$ and $\xi \in W^*$ be invariant under the diagonal action of $\wh G$. Let $(\gamma_i)_{i \in I} \in \on{Weil}(\eta, \ov{\eta})^I$. 
We will construct an excursion operator $S_{I, W, x, \xi, (\gamma_i)_{i \in I}}$ on $H_{J, V}^j = \restr{\mc  H_{J, V}^j  }{\Delta^J(\ov{\eta})}$, where $\Delta^J: X \sm N \hookrightarrow (X \sm N)^J$ is the diagonal morphism.

\begin{prop} [\emph{cf.}  proposition~4.12 in \cite{vincent}]  \label{proposition-fusion}
Let $I_1$, $I_2$ be two finite sets and $\zeta: I_1 \rightarrow I_2$ be a map. Let $\Delta_{\zeta}: X^{I_2} \rightarrow X^{I_1}$, $(x_j)_{j \in I_2} \mapsto (x_{\zeta(i)})_{i \in I_1}$ be the morphism associated to $\zeta$. Let $W$ be a representation of $\wh G^{I_1}$. We denote by $W^{\zeta}$ the representation of $\wh G^{I_2}$ which is the composition of $W$ with the morphism $\wh G^{I_2} \rightarrow \wh G^{I_1}$, $(g_j)_{j \in I_2} \mapsto (g_{\zeta(i)})_{i \in I_1}$. Then there is a canonical isomorphism of sheaves over $(X \sm N)^{I_2}$: 
$$\chi_{\zeta}: \Delta_{\zeta}^*(\mc H_{I_1, W}^j) \isom \mc H_{I_2, W^{\zeta}}^j .$$
This is called the fusion.
\end{prop}

\sssec{}
In the following the unions $J \cup I$ and $J \cup \{0\}$ always mean the disjoint unions.
Applying the above proposition to the map $J \hookrightarrow J \cup \{0\}$, we deduce 
$$\mc H_{J, V}^j \boxtimes E_{(X \sm N)} \isom \mc H_{J \cup \{0\}, V \boxtimes \bf 1}^j .$$
Let $\zeta: I \twoheadrightarrow \{0\}$.
Applying the above proposition to the map $(\on{Id}_J, \zeta): J \cup I \twoheadrightarrow J \cup \{0\}$, and denoting by $\Delta^I: X \hookrightarrow X^I$ the diagonal morphism, we deduce
$$\restr{\mc H_{J \cup I, V \boxtimes W}^j }{(X \sm N)^J \times \Delta^I(X \sm N)} \isom \mc H_{J \cup \{0\}, V \boxtimes W^{\zeta}}^j .$$

\begin{defi} [see d\'efinition~5.1 in \cite{vincent} where the notations are different]
The creation morphism $\mc C_x^{\sharp}$ is defined to be the composition
\begin{equation}   \label{equation-creation-operator-over-X-J-0}
\mc H_{J, V}^j \boxtimes E_{(X \sm N)} \isom \mc H_{J \cup \{0\}, V \boxtimes \bf 1}^j \xrightarrow{\mc H(\on{Id}_V \boxtimes x)} \mc H_{J \cup \{0\}, V \boxtimes W^{\zeta}}^j \isom \restr{\mc H_{J \cup I, V \boxtimes W}^j }{(X \sm N)^J \times \Delta^I(X \sm N)}, 
\end{equation}
where the middle morphism is the functoriality of the cohomology for $(\on{Id}_V, x): V \boxtimes {\bf 1} \rightarrow V \boxtimes W^{\zeta}$ (\cite[notation~4.9]{vincent}).

The annihilation morphism $\mc C_{\xi}^{\flat}$ is defined to be the composition
\begin{equation}  \label{equation-annihilation-operator-over-X-J-0}
\restr{\mc H_{J \cup I, V \boxtimes W}^j }{(X \sm N)^J \times \Delta^I(X \sm N)} \isom \mc H_{J \cup \{0\}, V \boxtimes W^{\zeta}}^j \xrightarrow{\mc H(\on{Id}_V \boxtimes \xi)} \mc H_{J \cup \{0\}, V \boxtimes \bf 1 }^j  \isom \mc H_{J, V}^j \boxtimes E_{(X \sm N)} .
\end{equation}
where the middle morphism is the functoriality of the cohomology for $(\on{Id}_V, \xi): V \boxtimes W^{\zeta} \rightarrow V \boxtimes \bf 1$.
\end{defi}

\sssec{}
We restrict (\ref{equation-creation-operator-over-X-J-0}) to $\Delta^{J \cup \{0\}}(\ov{\eta})$ and obtain 
$$\restr{\mc C_{  x}^{\sharp}}{\ov{\eta}}: \restr{ \mc  H_{J, V}^j }{\Delta^J(\ov{\eta}) }  \rightarrow  \restr{ \mc H_{J \cup I, V \boxtimes W }^{j}   }{  \Delta^{J \cup I}(\ov{\eta})  }  .$$
We restrict (\ref{equation-annihilation-operator-over-X-J-0}) to $\Delta^{J \cup \{0\}}(\ov{\eta})$ and obtain 
$$\restr{\mc C_{ \xi}^{\flat } }{\ov{\eta}}: \restr{ \mc H_{J \cup I, V \boxtimes W }^{j}   }{ \Delta^{J \cup I}(\ov{\eta}) }  \rightarrow \restr{ \mc  H_{J, V}^j }{\Delta^J(\ov{\eta}) } .$$

\begin{construction}     \label{construction-excursion-operator-on-cohomology}
With the notations in \ref{subsection-notation-I-W-J-V}, we construct an excursion operator $S_{I, W, x, \xi, (\gamma_i)_{i \in I}} $ acting on $H_{J, V}^j = \restr{\mc  H_{J, V}^j  }{\Delta^J(\ov{\eta})}$ as the composition of morphisms:
$$
\xymatrixrowsep{2pc}
\xymatrixcolsep{4pc}
\xymatrix{ 
 \restr{ \mc  H_{J, V}^j }{\Delta^J(\ov{\eta}) }  \ar[r]^{  \restr{\mc C_{  x}^{\sharp}}{\ov{\eta}} \quad \quad } 
&    \restr{ \mc H_{J \cup I, V \boxtimes W }^{j}   }{  \Delta^{J \cup I}(\ov{\eta})  }  \ar[d]^{ (\gamma_i)_{i\in I} } \\  
  \restr{ \mc  H_{J, V}^j }{\Delta^J(\ov{\eta}) }  
&  \restr{ \mc H_{J \cup I, V \boxtimes W }^{j}   }{ \Delta^{J \cup I}(\ov{\eta}) }  \ar[l]_{  \restr{\mc C_{ \xi}^{\flat } }{\ov{\eta}}  \quad \quad }
    }
$$
where the action of $\on{Weil}(\eta, \ov{\eta})^I$ on $ \restr{ \mc H_{J \cup I, V \boxtimes W}^{j}  }{ \Delta^{J \cup I}(\ov{\eta}) } $ is given by $(1, \on{Id}_I): \on{Weil}(\eta, \ov{\eta})^I \hookrightarrow \on{Weil}(\eta, \ov{\eta})^{J \cup I}$ and by Construction \ref{const-Weil-I-on-H-Delta-eta} applied to $J \cup I$.
\end{construction}

\begin{rem}
When $J = \emptyset$ and $V = \bf 1$, this construction coincides with Construction \ref{constr-excursion-operator-on-C-c}.
\end{rem}

The vector space $H_{J, V}^j$ is equipped with an action of the global Hecke algebra $\ms H_G$ and an action of $\on{Weil}(\eta, \ov{\eta})^J$ (Construction \ref{const-Weil-I-on-H-Delta-eta} applied to $J$).

\begin{lem}   \label{lem-S-commute-with-Hecke-on-H-J-V}
The action of $S_{I, W, x, \xi, (\gamma_i)_{i \in I}}$ on $H_{J, V}^j$ commutes with the action of $\ms H_G$ and the action of $\on{Weil}(\eta, \ov{\eta})^J$.
\end{lem}
\dem
For the Hecke algebra: by \cite{vincent}, the action of Hecke algebra commutes with the functoriality, the fusion, the action of the partial Frobenius morphisms and the action of Galois groups.

For $\on{Weil}(\eta, \ov{\eta})^J$: by \cite{vincent}, the action of $\on{Weil}(\eta, \ov{\eta})^J$ commutes with the functoriality and the fusion. Moreover, by the above construction, the action of $\on{Weil}(\eta, \ov{\eta})^J$ on $\restr{ \mc H_{J \cup I, V \boxtimes W}^{j}  }{ \Delta^{J \cup I}(\ov{\eta}) }$ is given by 
$$(\on{Id}_J, 1): \on{Weil}(\eta, \ov{\eta})^J \rightarrow \on{Weil}(\eta, \ov{\eta})^J \times \on{Weil}(\eta, \ov{\eta})^I $$ and the action of 
$\on{Weil}(\eta, \ov{\eta})^I$ on $\restr{ \mc H_{J \cup I, V \boxtimes W}^{j}  }{ \Delta^{J \cup I}(\ov{\eta}) }$ is given by 
$$(1, \on{Id}_I): \on{Weil}(\eta, \ov{\eta})^I \rightarrow \on{Weil}(\eta, \ov{\eta})^J \times \on{Weil}(\eta, \ov{\eta})^I .$$ Thus the action of $\on{Weil}(\eta, \ov{\eta})^J$ commutes with the action of $\on{Weil}(\eta, \ov{\eta})^I$.
 \cqfd

\begin{rem}
By Lemma \ref{lem-S-commute-with-Hecke-on-H-J-V}, $S_{I, W, x, \xi, (\gamma_i)_{i \in I}}$ sends the Hecke finite part $H_{J, V}^{j, \, \on{Hf}}$ (\cite[d\'efi.~8.19]{vincent}) to $H_{J, V}^{j, \, \on{Hf}}$.
\end{rem}

\begin{prop}
The excursion operators constructed in Construction \ref{construction-excursion-operator-on-cohomology} above have the same properties as in Section \ref{subsection-excursion-operator}. Moreover, for any $u \in |X \sm N|$ and any ideal $\ms I$ of $\ms H_{G, u}$ of finite codimension, the quotient $E$-vector space $H_{J, V}^j / \ms I \cdot H_{J, V}^j$ has finite dimension and admits a similar decomposition as in Theorem \ref{thm-decomposition-of-C-c-quotient-I-by-param-Langlands}.
\end{prop}

\quad

\section{Compatibility with the constant term morphisms}

\subsection{Commutativity of the excursion operators and the constant term morphisms}

\sssec{}
Let $P$ be a parabolic subgroup of $G$ and $M$ its Levi quotient.
In \cite[\S~1.5]{cusp-coho}, for any $\nu \in \wh \Lambda_{Z_M / Z_G}^{\Q}$, we have defined an open and closed substack $\Bun_{M, N}^{\nu}$ of $\Bun_{M, N}$ and we have $\Bun_{M, N} = \bigsqcup_{\nu \in \wh \Lambda_{Z_M / Z_G}^{\Q}} \Bun_{M, N}^{\nu}$. As in \cite[\S~3.4]{cusp-coho}, let $\Bun_{M, N}^{' \, \nu}(\Fq)=\Bun_{M, N}^{\nu}(\Fq) \overset{P(\mc O_N)}{\times} G(\mc O_N)$. We have defined the constant term morphism $C_G^{P, \, \nu}: C_c(\Bun_{G, N}(\Fq) / \Xi, E) \rightarrow C_c(\Bun_{M, N}^{' \, \nu}(\Fq) / \Xi, E) $, which coincides with the classical constant term morphism (\cite[Remark~3.5.11 and Example~3.5.15]{cusp-coho}).

\sssec{}
Let $I$ be a finite set and $W$ be a representation of $\wh G^I$. Let $x \in W$ and $\xi \in W^*$ be invariant by the diagonal action of $\wh G$. Let $(\gamma_i)_{i \in I} \in \on{Weil}(\eta, \ov{\eta})^I$. In Construction \ref{const-Weil-I-on-H-Delta-eta} and Construction \ref{constr-excursion-operator-on-C-c}, we defined the excursion operator $S_{I, W, x, \xi, (\gamma_i)_{i \in I}}^G $ acting on $C_c(\Bun_{G, N}(\Fq) / \Xi, E)$. 
Similarly, we define an action of $\on{Weil}(\eta, \ov{\eta})^I$ on $\mc H_{M, N, I, W}^{' \, 0, \, \nu}$ (defined in Definition \ref{def-H-M-N-I-W-leq-mu-nu}, where we view $W$ as a representation of $\wh M^I$ via $\wh M^I \hookrightarrow \wh G^I$) to be the composition of morphisms $$\restr{ \mc H_{M, N, I, W}^{' \, 0, \, \nu}   }{\Delta(\ov{\eta})} \xrightarrow{\mf{sp}^*} \restr{ \mc H_{M, N, I, W}^{' \, 0, \, \nu}   }{\ov{\eta^I}} \xrightarrow{ (\gamma_i)_{i\in I}  } \restr{ \mc H_{M, N, I, W}^{' \, 0, \, \nu}   }{\ov{\eta^I}} \xrightarrow{(\mf{sp}^*)^{-1}} \restr{ \mc H_{M, N, I, W}^{' \, 0, \, \nu}   }{\Delta(\ov{\eta})} $$
and we define the excursion operator $S_{I, W, x, \xi, (\gamma_i)_{i \in I}}^M $ acting on $C_c(\Bun_{M, N}^{' \, \nu}(\Fq) / \Xi, E)$ to be the composition of morphisms:
$$
\xymatrixrowsep{2pc}
\xymatrixcolsep{3pc}
\xymatrix{ 
    C_c(\Bun_{M, N}^{' \, \nu}(\Fq) / \Xi, E)   \ar[r]^{ \quad \quad \mc C_{  x}^{\sharp}} 
&    \restr{ \mc H_{M, N, I, W}^{' \, 0, \, \nu}   }{\Delta(\ov{\eta})}  \ar[d]^{(\gamma_i)_{i\in I}} \\  
C_c(\Bun_{M, N}^{' \, \nu}(\Fq) / \Xi, E)  
&  \restr{ \mc H_{M, N, I, W}^{' \, 0, \, \nu}  }{\Delta(\ov{\eta})}  \ar[l]_{ \quad \quad \mc C_{ \xi}^{\flat }}.
    }
$$

\begin{prop}   \label{prop-CT-commutes-with-excursion}
For any $\nu \in \wh \Lambda_{Z_M / Z_G}^{\Q}$, the following diagram is commutative:
$$
\xymatrixrowsep{2pc}
\xymatrixcolsep{6pc}
\xymatrix{
C_c(\Bun_{G, N}(\Fq) / \Xi, E)   \ar[r]^{\quad S^G_{I, W, x, \xi, (\gamma_i)_{i \in I}} \quad}  \ar[d]^{C_{G}^{P, \, \nu}}
& C_c(\Bun_{G, N}(\Fq) / \Xi, E) \ar[d]^{C_{G}^{P, \, \nu}} \\
C_c(\Bun_{M, N}^{' \, \nu}(\Fq) / \Xi, E)  \ar[r]^{\quad S^M_{I, W, x, \xi, (\gamma_i)_{i \in I}} \quad}
& C_c(\Bun_{M, N}^{' \, \nu}(\Fq) / \Xi, E) .
}
$$
\end{prop}

\dem
The proof consists of 4 parts.

\noindent\textbf{Step~1.} The constant term morphisms commute with the creation operators: the following diagram is commutative
\begin{equation}
\xymatrix{
C_c(\Bun_{G, N}(\Fq) / \Xi, E)  = \restr{\mc H_{G, N, \{ 0 \}, \bf{1} }^0}{\ov \eta}  \ar[r]^{\quad \quad \quad \quad \mc H(x)}  \ar[d]^{C_{G}^{P, \, \nu}}    \ar@{}[rd]|{(a)}
& \restr{\mc H_{G, N, \{ 0 \}, W^{\zeta_I} }^0}{\ov \eta}  \ar[r]^{\chi_{\zeta_I}^{-1}}_{\sim}   \ar[d]^{\mc C_{G}^{P, \, \nu}}      \ar@{}[rd]|{(b)}
& \restr{\mc H_{G, N, I, W}^0}{\Delta(\ov \eta)} \ar[d]^{\mc C_{G}^{P, \, \nu}} \\
C_c(\Bun_{M, N}^{' \, \nu}(\Fq) / \Xi, E)= \restr{\mc H_{M, N, \{ 0 \}, \bf{1} }^{' \, 0, \, \nu} }{\ov \eta}   \ar[r]^{\quad \quad \quad \quad \mc H(x)}
& \restr{\mc H_{M, N, \{ 0 \}, W^{\zeta_I} }^{' \, 0, \, \nu} }{\ov \eta}   \ar[r]^{\chi_{\zeta_I}^{-1}}_{\sim} 
& \restr{\mc H_{M, N, I, W}^{' \, 0, \, \nu}}{\Delta(\ov \eta)}
}
\end{equation}
where $W^{\zeta_I}$ is the representation of $\wh G$ via the diagonal inclusion $\wh G \hookrightarrow \wh G^I$, $\chi_{\zeta_I}^{-1}$ and $\mc H(x)$ are defined in \cite[d\'efi.~5.1]{vincent} and $\mc C_x^{\sharp}=\chi_{\zeta_I}^{-1} \circ \mc H(x)$.

Indeed, the commutativity of (a) comes from the fact that $\mc C_{G}^{P, \, \nu}$ is functorial on $W$. The commutativity of (b) follows from Lemma \ref{lem-CT-commutes-with-fusion} below applied to $J = \{ 0 \}$.

\noindent\textbf{Step~2.} The constant term morphisms commute with the specialization morphisms: the following diagram is commutative:
\begin{equation}    \label{equation-TC-commutes-with-specialization}
\xymatrix{
\restr{  \mc H_{G, N, I, W}^j  }{ \Delta(\ov{\eta})}   \ar[r]^{\mf{sp}^*}_{\simeq}   \ar[d]^{\mc C_{G}^{P, \, \nu}}
& \restr{  \mc H_{G, N, I, W}^j  }{ \ov{\eta^{I}} }   \ar[d]^{\mc C_{G}^{P, \, \nu}}   \\
\restr{  \mc H_{M, N, I, W}^{' \, j, \, \nu}  }{ \Delta(\ov{\eta}) }   \ar[r]^{\mf{sp}^*}_{\simeq}
& \restr{  \mc H_{M, N, I, W}^{' \, j, \, \nu} }{ \ov{\eta^{I}} } .
}
\end{equation}

Indeed, by \cite[\S~3.5.6]{cusp-coho}, for any $\mu \in \wh \Lambda^{+, \Q}_{G^{\mr{ad}}}$ and $\nu \in \wh \Lambda_{Z_M / Z_G}^{\Q}$, the morphism $\mc C_{G}^{P, \, \leq \mu, \, \nu}$ is defined over $\Omega^{\leq \mu, \, \nu}$, which is an open subscheme of $X^I$ containing $\Delta(\eta)$. 
We have a morphism of sheaves
$$\mc C_{G}^{P, \, \leq \mu, \, \nu}: \restr{  \mc H_{G, N, I, W}^{j, \, \leq \mu}  }{\Omega^{\leq \mu, \, \nu}} \rightarrow \restr{  \mc H_{M, N, I, W}^{' \, j, \, \leq \mu, \, \nu}  }{\Omega^{\leq \mu, \, \nu}} .$$
Since the homomorphism of specialization is functorial, for any $\mu \in \wh \Lambda^{+, \Q}_{G^{\mr{ad}}}$, it induces a commutative diagram:
\begin{equation}    \label{equation-CT-commutes-with-specialization-leq-mu}
\xymatrix{
\restr{  \mc H_{G, N, I, W}^{j, \, \leq \mu}  }{ \Delta(\ov{\eta})}   \ar[r]^{\mf{sp}^*}   \ar[d]^{\mc C_{G}^{P, \, \leq \mu, \, \nu}}
& \restr{  \mc H_{G, N, I, W}^{j, \, \leq \mu}  }{ \ov{\eta^{I}} }   \ar[d]^{\mc C_{G}^{P, \, \leq \mu, \, \nu}}   \\
\restr{  \mc H_{M, N, I, W}^{' \, j, \, \leq \mu, \, \nu}  }{ \Delta(\ov{\eta}) }   \ar[r]^{\mf{sp}^*}
& \restr{  \mc H_{M, N, I, W}^{' \, j, \, \leq \mu, \,  \nu} }{ \ov{\eta^{I}} } .
}
\end{equation}

Moreover, by \cite[Remark~3.5.4]{cusp-coho}, for $\mu_1 \leq \mu_2$, we have $\Omega^{\leq \mu_2, \, \nu} \subset \Omega^{\leq \mu_1, \, \nu}$. By \cite[\S~3.5.9]{cusp-coho}, we have a commutative diagram:
$$
\xymatrix{
\restr{  \mc H_{G, N, I, W}^{j, \, \leq \mu_1}  }{U^{\leq \mu_2, \, \nu}} \ar[r]    \ar[d]^{\mc C_{G}^{P, \, \leq \mu_1, \, \nu}}
& \restr{  \mc H_{G, N, I, W}^{j, \, \leq \mu_2}  }{U^{\leq \mu_2, \, \nu}}  \ar[d]^{\mc C_{G}^{P, \, \leq \mu_2, \, \nu}}  \\
\restr{  \mc H_{M, N, I, W}^{' \, j, \, \leq \mu_1, \, \nu}  }{U^{\leq \mu_2, \, \nu}} \ar[r]
& \restr{  \mc H_{M, N, I, W}^{' \, j, \, \leq \mu_2, \, \nu}  }{U^{\leq \mu_2, \, \nu}}  
}
$$
The diagram (\ref{equation-CT-commutes-with-specialization-leq-mu}) is compatible with $\mu_1 \leq \mu_2$. Taking the inductive limit on $\mu$ in (\ref{equation-CT-commutes-with-specialization-leq-mu}), we deduce that the diagram (\ref{equation-TC-commutes-with-specialization}) is commutative.

\noindent\textbf{Step~3.} The following diagram is commutative:
\begin{equation}
\xymatrixrowsep{2pc}
\xymatrixcolsep{3pc}
\xymatrix{
\restr{  \mc H_{G, N, I, W}^j  }{ \ov{\eta^{I}}  }   \ar[r]^{(\gamma_i)_{i\in I}}   \ar[d]^{\mc C_{G}^{P, \, \nu}}
& \restr{  \mc H_{G, N, I, W}^j  }{ \ov{\eta^{I}} }   \ar[d]^{\mc C_{G}^{P, \, \nu}}   \\
\restr{  \mc H_{M, N, I, W}^{' \, j, \, \nu}  }{ \ov{\eta^{I}}  }   \ar[r]^{(\gamma_i)_{i\in I}}
& \restr{  \mc H_{M, N, I, W}^{' \, j, \, \nu} }{ \ov{\eta^{I}} } .
}
\end{equation}
Indeed, by Lemma \ref{lem-CT-commutes-with-partial-Frob} below, $\mc C_{G}^{P, \, \nu}$ commutes with the partial Frobenius morphisms. And $\mc C_{G}^{P, \, \nu}$ commutes with the action of $\pi_1(\eta^I, \ov{\eta^I})$ (since $\mc C_{G}^{P, \, \nu}$ is defined over $\eta^I$). So the constant term morphism commutes with the action of $\on{FWeil}(\eta^I, \ov{\eta^I})$ and thus with the action of $\on{Weil}(\ov F / F)^I$.

\noindent\textbf{Step~4.} The constant term morphisms commute with the annihilation operators: in the same way as in Step~1.
\cqfd

\sssec{}
Let $I$, $J$ be two finite sets and $\zeta: I \rightarrow J$ be a map. Let $\Delta_{\zeta}: X^J \rightarrow X^I$, $(x_j)_{j \in J} \mapsto (x_{\zeta(i)})_{i \in I}$ be the morphism associated to $\zeta$. Let $W$ be a representation of $\wh G^I$ (resp. $\wh M^I$). We denote by $W^{\zeta}$ the representation of $\wh G^J$ (resp. $\wh M^J$) which is the composition of $W$ with the morphism $\wh G^J \rightarrow \wh G^I$, $(g_j)_{j \in J} \mapsto (g_{\zeta(i)})_{i \in I}$ (resp. $\wh M^J \rightarrow \wh M^I$, $(m_j)_{j \in J} \mapsto (m_{\zeta(i)})_{i \in I}$). Let $\eta^J$ be the generic point of $X^J$. Fix a geometric point $\ov{\eta^J}$ over $\eta^J$.

\begin{lem}[The constant term morphisms commute with the fusion]\label{lem-CT-commutes-with-fusion}
The following diagram is commutative:
$$
\xymatrix{
\restr{\mc H_{G, N, I, W}^j}{\Delta_{\zeta}(\ov{\eta^J})} \ar[d]^{\mc C_{G}^{P, \, \nu}}  \ar[r]^{\chi_{\zeta}}_{\sim} 
&  \restr{\mc H_{G, N, J, W^{\zeta} }^j}{\ov{\eta^J}}   \ar[d]^{\mc C_{G}^{P, \, \nu}}  \\
\restr{\mc H_{M, N, I, W}^{' \, j, \, \nu}}{\Delta_{\zeta}(\ov{\eta^J})}   \ar[r]^{\chi_{\zeta}}_{\sim} 
&  \restr{\mc H_{M, N, J, W^{\zeta} }^{' \, j, \, \nu} }{\ov{\eta^J}}.
}
$$
where $\chi_{\zeta}$ are defined in Proposition \ref{proposition-fusion} applied to the reductive group $G$ and to the reductive group $M$ respectively.
\end{lem}
\dem The proof consists of 3 steps.

\noindent\textbf{Step~1.} By \cite[Section~5.3]{bd} or \cite[Section~15]{BR}, we have a commutative diagram of Beilinson-Drinfeld grassmannians (see \cite[Definition~1.1.11 and~3.1.2]{cusp-coho} for a reminder of definitions and vertical morphisms):
$$
\xymatrix{
\Gr_{G, J, W^{\zeta} }   \ar[r]_{\simeq \quad}
& \restr{ \Gr_{G, I, W} }{\Delta_{\zeta}(X^J)}  \ar[r]^{\quad \Delta_{\zeta}^G}
& \Gr_{G, I, W}  \\
\Gr_{P, J, W^{\zeta} }   \ar[r]_{\simeq \quad}   \ar[u]^{i_J^0}   \ar[d]_{\pi_J^0} 
& \restr{ \Gr_{P, I, W} }{\Delta_{\zeta}(X^J)}  \ar[u]   \ar[d]   \ar[r]^{\quad \Delta_{\zeta}^P}
& \Gr_{P, I, W} \ar[u]^{i_I^0}   \ar[d]_{\pi_I^0}   \\
 \Gr_{M, J, W^{\zeta} }   \ar[r]_{\simeq \quad}
& \restr{ \Gr_{M, I, W} }{\Delta_{\zeta}(X^J)}   \ar[r]^{\quad \Delta_{\zeta}^M}
& \Gr_{M, I, W}
}
$$
where all squares are cartesian.
Moreover, we have a commutative diagram of categories, where the back and front faces are recalled in \cite[Theorem~3.2.6]{cusp-coho} applied to $I$ and $J$ respectively. 

{\resizebox{14cm}{!}{ $$
\xymatrix{
& \on{Perv}_{G_{I, \infty}}(\Gr_{G, I}, \Ql)^{\on{MV}}  \ar[rr]^{(\pi_I^0)_!(i_I^0)^*}   \ar[ld]_{(\Delta_{\zeta}^G)^*}
& & \on{Perv}_{M_{I, \infty}}(\Gr_{M, I}, \Ql)^{\on{MV}}      \ar[ld]_{(\Delta_{\zeta}^M)^*}  \\
\on{Perv}_{G_{J, \infty}}(\Gr_{G, J}, \Ql)^{\on{MV}}     \ar[rr]^{\quad \quad \quad \quad \quad \quad (\pi_J^0)_!(i_J^0)^*}   
& & \on{Perv}_{M_{J, \infty}}(\Gr_{M, J}, \Ql)^{\on{MV}}   \\
& \on{Rep}_E(\wh G^I)    \ar[rr]^{ \on{Res} \quad \quad \quad \quad \quad }    \ar[ld]     \ar[uu] 
& & \on{Rep}_E(\wh M^I)      \ar[ld]  \ar[uu]  \\
\on{Rep}_E(\wh G^J)     \ar[rr]^{\on{Res}}     \ar[uu]
& &  \on{Rep}_E(\wh M^J)   \ar[uu]  
}
$$
}
}

By \emph{loc.~cit.} Theorem 3.2.6,
%the compatibility of the geometric Satake equivalence with the constant term functor (recalled in $loc$. $cit$. Theorem 3.2.6). 
there exists canonical isomorphisms
\begin{equation}    \label{equation-CT-S-G-equal-S-M}
 (\pi^0_I)_! (i^0_I)^* \mc S_{G, I, W} \isom  \mc S_{M, I, W}, \quad  (\pi^0_J)_! (i^0_J)^* \mc S_{G, J, W^{\zeta} }  \isom   \mc S_{M, J, W^{\zeta} } . 
\end{equation}
The commutativity of the left and right face follows from the compatibility of the geometric Satake equivalence with the fusion. In other words, there exists canonical isomorphisms 
\begin{equation}     \label{equation-Delta-S-I-equal-S-J}
(\Delta_{\zeta}^G)^* \mc S_{G, I, W}   \simeq   \mc S_{G, J, W^{\zeta} }, \quad (\Delta_{\zeta}^M)^* \mc S_{M, I, W}  \simeq  \mc S_{M, J, W^{\zeta} }.
\end{equation}
The commutativity of the bottom face is evident.

We deduce 
a commutative diagram:
\begin{equation}     \label{equation-CT-commutes-with-fusion-for-S}
\xymatrix{
(\pi^0_J)_! (i^0_J)^* (\Delta_{\zeta}^G)^* \mc S_{G, I, W}  \ar[d]^{\simeq}_{\text{base change}}   \ar[r]^{\simeq} 
&  (\pi^0_J)_! (i^0_J)^* \mc S_{G, J, W^{\zeta} }  \ar[dd]^{\simeq}  \\
(\Delta_{\zeta}^M)^* (\pi^0_I)_! (i^0_I)^* \mc S_{G, I, W}  \ar[d]^{\simeq}   \\
(\Delta_{\zeta}^M)^* \mc S_{M, I, W}    \ar[r]^{\simeq} 
&  \mc S_{M, J, W^{\zeta} } .
}
\end{equation}
%where morphism (*) comes from the base change.

\noindent\textbf{Step~2.} Similarly, we have a commutative diagram of stacks of shtukas, where all squares are cartesian:
$$
\xymatrix{
\Cht_{G, N, J, W^{\zeta} }   \ar[r]_{\simeq \quad \quad}
& \restr{ \Cht_{G, N, I, W} }{\Delta_{\zeta}(X^J)}  \ar[r]^{\quad \Delta_{\zeta}^G}
& \Cht_{G, N, I, W}  \\
\Cht_{P, N, J, W^{\zeta} }'   \ar[r]_{\simeq \quad \quad }   \ar[u]^{i_J}   \ar[d]_{\pi_J} 
& \restr{ \Cht_{P, N, I, W}' }{\Delta_{\zeta}(X^J)}  \ar[u]   \ar[d]   \ar[r]^{\quad \Delta_{\zeta}^P}
& \Cht_{P, N, I, W}' \ar[u]^{i_I}   \ar[d]_{\pi_I}   \\
\Cht_{M, N, J, W^{\zeta} }'  \ar[r]_{\simeq \quad \quad }
& \restr{ \Cht_{M, N, I, W}' }{\Delta_{\zeta}(X^J)}   \ar[r]^{\quad \Delta_{\zeta}^M}
& \Cht_{M, N, I, W}'
}
$$

By \emph{loc.~cit.} Definition 2.4.5, $\mc F_{G, N, I, W}$ is defined as the inverse image of $\mc S_{G, I, W}$ by the morphism $\Cht_{G, N, I, W} \rightarrow [G_{I, \infty} \backslash \Gr_{G, I, W}]$. Thus (\ref{equation-Delta-S-I-equal-S-J}) implies 
\begin{equation}     \label{equation-Delta-F-G-I-isom-F-G-J}    
(\Delta_{\zeta}^G)^* \mc F_{G, N, I, W}   \simeq   \mc F_{G, N, J, W^{\zeta} }, \quad (\Delta_{\zeta}^M)^* \mc F_{M, N, I, W}'  \simeq  \mc F_{M, N, J, W^{\zeta} }'.
\end{equation}

Moreover, in \emph{loc.~cit.} (3.31), we constructed a canonical morphism:
\begin{equation}     \label{equation-pi-i-F-G-to-F-M}
 (\pi_I)_! (i_I)^* \mc F_{G, N, I, W} \rightarrow  \mc F_{M, N, I, W}', \quad ( \text{resp. } (\pi_J)_! (i_J)^* \mc F_{G, N, J, W^{\zeta} }  \rightarrow  \mc F_{M, N, J, W^{\zeta} }' ). 
\end{equation}
The construction uses the isomorphism (\ref{equation-CT-S-G-equal-S-M}) and $\on{Tr}_{\pi_{I, d}}: (\pi_{I, d})_! (\pi_{I, d})^! \rightarrow \on{Id}$ (resp. $\on{Tr}_{\pi_{J, d}}: (\pi_{J, d})_! (\pi_{J, d})^! \rightarrow \on{Id}$), where $\pi_{I, d}$ (resp. $\pi_{J, d}$) is defined in %the right (resp. left) vertical morphism in 
the following diagram:
$$
\xymatrix{
\Cht_{P, N, J, W^{\zeta} }'   \ar[r]_{\simeq \quad}     \ar[d]_{\pi_{J, d}} 
& \restr{ \Cht_{P, N, I, W}' }{\Delta_{\zeta}(X^J)}   \ar[d]   \ar[r]^{\quad \Delta_{\zeta}^P}
& \Cht_{P, N, I, W}'    \ar[d]_{\pi_{I, d}}   \\
\wt{\Cht}_{M, N, J, W^{\zeta} }'    \ar[r]_{\simeq \quad}     \ar[d]_{ \wt{\pi_{J, d}^0} } 
& \restr{ \wt{\Cht}_{M, N, I, W}'  }{\Delta_{\zeta}(X^J)}   \ar[r]   \ar[d]
&\wt{\Cht}_{M, N, I, W}'  \ar[d]_{ \wt{\pi_{I, d}^0 } }  \\
\Cht_{M, N, J, W^{\zeta} }'  \ar[r]_{\simeq \quad}
& \restr{ \Cht_{M, N, I, W}' }{\Delta_{\zeta}(X^J)}   \ar[r]^{\quad \Delta_{\zeta}^M}
& \Cht_{M, N, I, W}'
}
$$
where $\pi_{I} = \wt{\pi_{I, d}^0} \circ \pi_{I, d}$ and $\pi_{J} = \wt{\pi_{J, d}^0} \circ \pi_{J, d}$. %The morphisms (\ref{equation-pi-i-F-G-to-F-M}) is given by the composition of some isomorphisms (containing (\ref{equation-CT-S-G-equal-S-M}))) 

The above diagram is commutative and all squares are cartesian. By \cite[XVIII, th.~2.9]{sga4}, the trace morphism is compatible with base change, thus 
\begin{equation}   \label{equation-trace-J-equal-Delta-trace-I}
\on{Tr}_{\pi_{J, d}} = (\Delta^M_{\zeta})^* \on{Tr}_{\pi_{I, d}}.
\end{equation}
We deduce from (\ref{equation-CT-commutes-with-fusion-for-S}), (\ref{equation-Delta-F-G-I-isom-F-G-J}), (\ref{equation-pi-i-F-G-to-F-M}) and (\ref{equation-trace-J-equal-Delta-trace-I}) a commutative diagram:
\begin{equation}   \label{eqaution-Delta-CT-I-to-CT-J-functor}
\xymatrix{
(\pi_J)_! (i_J)^* (\Delta_{\zeta}^G)^* \mc F_{G, N, I, W}  \ar[d]^{\simeq}_{\text{base change}}   \ar[r]^{\simeq}
&  (\pi_J)_! (i_J)^* \mc F_{G, N, J, W^{\zeta} }  \ar[dd]^{ (\ref{equation-CT-S-G-equal-S-M}) + \on{Tr}_{\pi_{J, d}}  }  \\
(\Delta_{\zeta}^M)^* (\pi_I)_! (i_I)^* \mc F_{G, N, I, W}  \ar[d]_{ (\ref{equation-CT-S-G-equal-S-M}) + \on{Tr}_{\pi_{I, d}}  }    \\
(\Delta_{\zeta}^M)^* \mc F_{M, N, I, W}    \ar[r]^{\simeq} 
&  \mc F_{M, N,  J, W^{\zeta} } .
}
\end{equation}

\noindent\textbf{Step~3.} In \emph{loc.~cit.} 3.5.6, we defined $$\mc C_{G, I}^{P, \, \nu}: \restr{\mc H_{G, N, I, W}^{j, \, \leq \mu}}{U^{\leq \mu, \, \nu}} \rightarrow \restr{ \mc H_{M, N, I, W}^{' \, j, \,  \leq \mu, \, \nu} }{{U^{\leq \mu, \, \nu}}} $$
$$(\text{resp. } \mc C_{G, J}^{P, \, \nu}: \restr{\mc H_{G, N, J, W^{\zeta} }^{j, \, \leq \mu} }{U^{\leq \mu, \, \nu}}  \rightarrow  \restr{ \mc H_{M, N,  J, W^{\zeta} }^{' \, j, \, \leq \mu, \, \nu} }{U^{\leq \mu, \, \nu}}  )$$
by using $\on{adj_{i_I}}: \on{Id} \rightarrow (i_I)_*(i_I)^*$ (resp. $\on{adj_{i_J}}: \on{Id} \rightarrow (i_J)_*(i_J)^*$) and (\ref{equation-pi-i-F-G-to-F-M}). The adjunction morphism is compatible with base change, thus
$$\on{adj}_{i_J} = (\Delta^G_{\zeta})^* \on{adj}_{i_I} .$$
This fact together with the commutativity of the diagram (\ref{eqaution-Delta-CT-I-to-CT-J-functor}) induce a commutative diagram
\begin{equation}    \label{equation-Delta-CT-I-to-CT-J-coho}
\xymatrix{
(\Delta_{\zeta}^G)^* \restr{\mc H_{G, N, I, W}^{j, \, \leq \mu}}{U^{\leq \mu, \, \nu}}  \ar[d]^{\mc C_{G, I}^{P, \, \nu}}  \ar[r]^{\simeq} 
&  \restr{\mc H_{G, N, J, W^{\zeta} }^{j, \, \leq \mu} }{U^{\leq \mu, \, \nu}} \ar[d]^{\mc C_{G, J}^{P, \, \nu}}  \\
(\Delta_{\zeta}^M)^* \restr{ \mc H_{M, N, I, W}^{' \, j, \,  \leq \mu, \, \nu} }{U^{\leq \mu, \, \nu}}    \ar[r]^{\simeq} 
&  \restr{ \mc H_{M, N,  J, W^{\zeta} }^{' \, j, \, \leq \mu, \, \nu} }{U^{\leq \mu, \, \nu}} .
}
\end{equation}

In particular, we have a commutative diagram:
$$
\xymatrix{
\restr{\mc H_{G, N, I, W}^{j, \, \leq \mu} }{\Delta_{\zeta}(\ov{\eta^J})} \ar[d]^{\mc C_{G, I}^{P, \, \nu}}  \ar[r]^{\chi_{\zeta}}_{\sim} 
& \restr{\mc H_{G, N, J, W^{\zeta} }^{j, \, \leq \mu}  }{ \ov{\eta^J} }    \ar[d]^{\mc C_{G, J}^{P, \, \nu}} \\
\restr{\mc H_{M, N, I, W}^{' \, j, \,  \leq \mu, \, \nu} }{\Delta_{\zeta}(\ov{\eta^J})}  \ar[r]^{\chi_{\zeta}}_{\sim} 
& \restr{\mc H_{M, N, J, W^{\zeta} }^{' \, j, \,  \leq \mu, \, \nu}  }{ \ov{\eta^J} } \,  .
}
$$
Taking the limit on $\mu$, we prove the lemma.
\cqfd

\sssec{}
For any partition $(I_1, \cdots, I_k)$ of $I$, let $\Cht_{G, N, I, W }^{(I_1, \cdots, I_k)} $ be the stack of iterated shtukas defined in \cite[d\'efi.~2.1]{vincent} and let $\mc F_{G, N, I, W}^{(I_1, \cdots, I_k), \,\Xi}$ be the perverse sheaf on $\Cht_{G, N, I, W }^{(I_1, \cdots, I_k)} / \Xi$ defined in \emph{loc.~cit.} d\'efinition 4.5. Note that $\Cht_{G, N, I, W}$ defined in \ref{subsection-Cht-G-N-I-W} corresponds to the partition $(I)$.
Moreover, by \emph{loc.~cit.} corollaire 2.18 and d\'efinition 4.1, there is an equality:
\begin{equation}   \label{equation-coho-independ-of-partition}
\begin{aligned}
\mc H_{G, N, I, W}^{j, \, \leq\mu} & = R^j (\mf{p}_G)_! \left( \restr{ \mc F_{G, N, I, W}^{\Xi} }{  \Cht_{G, N, I, W}^{\leq \mu} / \Xi   } \right) \\
& = R^j (\mf{p}_G)_! \left( \restr{ \mc F_{G, N, I, W}^{(I_1, \cdots, I_k), \,\Xi} }{  \Cht_{G, N, I, W}^{(I_1, \cdots, I_k), \,\leq \mu} / \Xi   } \right)
\end{aligned}
\end{equation}
In \emph{loc.~cit.} \S~3 and 4.3, $\Cht_{G, N, I, W }^{(I_1, \cdots, I_k)} $ and thus $\mc H_{G, N, I, W}^{j} $ are equipped with the partial Frobenuis morphisms.

\begin{lem}   \label{lem-CT-commutes-with-partial-Frob}
The constant term morphisms commute with the partial Frobenius morphisms.
\end{lem}
\dem
The construction of constant term morphisms in \cite{cusp-coho} works for any partition $(I_1, \cdots, I_k)$, i.e. we can use
\begin{equation}
\Cht_{G, N, I, W }^{(I_1, \cdots, I_k)} \xleftarrow{i} \Cht_{P, N, I, W }^{' \, (I_1, \cdots, I_k)} \xrightarrow{\pi} \Cht_{M, N, I, W }^{' \, (I_1, \cdots, I_k)}
\end{equation}
and the cohomological correspondence 
\begin{equation}   \label{equation-CT-functor-F-G-I-1-I-k-to-F-M-I-1-I-k}
\pi_! i^* \mc F_{G, N, I, W}^{(I_1, \cdots, I_k)} \rightarrow \mc F_{M, N, I, W}^{(I_1, \cdots, I_k) \, '} .
\end{equation}
to construct the constant term morphism
\begin{equation}  \label{equation-C-G-P-iterated-coho}
\mc C_G^{P, \, \nu}: \restr{\mc H_{G, N, I, W}^{j}}{\eta^I} \rightarrow \restr{\mc H_{M, N, I, W}^{' \, j, \, \nu}}{\eta^I} .
\end{equation}
(By (\ref{equation-coho-independ-of-partition}), the morphism (\ref{equation-C-G-P-iterated-coho}) is independent of the choice of the partition $(I_1, \cdots, I_k)$, thus coincides with the constant term morphism we used before.)

The following diagram is commutative:
$$
\xymatrixrowsep{2pc}
\xymatrixcolsep{5pc}
\xymatrix{
\restr{ \Cht_{G, N, I, W }^{(I_1, \cdots, I_k)} }{\ov{\eta^I}}  \ar[r]^{\on{Fr}_{I_1}}
& \restr{ \Cht_{G, N, I, W}^{(I_2, \cdots, I_k,  I_1)} }{\on{Frob}_{I_1} (\ov{\eta^I})}  \\
\restr{ \Cht_{P, N, I, W }^{' \, (I_1, \cdots, I_k)} }{\ov{\eta^I}}   \ar[r]^{\on{Fr}_{I_1}}   \ar[u]^{i_1}   \ar[d]_{\pi_1}
& \restr{ \Cht_{P, N, I, W}^{' \, (I_2, \cdots, I_k, I_1)} }{\on{Frob}_{I_1} (\ov{\eta^I})}   \ar[u]^{i_2}   \ar[d]_{\pi_2}  \\
\restr{ \Cht_{M, N, I, W }^{' \, (I_1, \cdots, I_k)} }{\ov{\eta^I}}   \ar[r]^{\on{Fr}_{I_1}}
& \restr{ \Cht_{M, N, I, W}^{' \, (I_2, \cdots, I_k, I_1)} }{\on{Frob}_{I_1} (\ov{\eta^I})} 
}
$$
where $\on{Fr}_{I_1}$ is defined in \cite[\S~3]{vincent}, the squares are Cartesian up to homeomorphism which is locally radical (so we have proper base change).

We have a commutative diagram, where $F_{I_1}$ is defined in \emph{loc.~cit.} proposition 3.3 and the vertical functors come from (\ref{equation-CT-functor-F-G-I-1-I-k-to-F-M-I-1-I-k}): 
%and the vertical morphisms are defined in \cite{cusp-coho} construction 3.3.5: 
$$
\xymatrix{
(\on{Fr}_{I_1})^*  (\pi_2)_! (i_2)^*  (\mc F_{G, N, I, W}^{(I_2, \cdots, I_1)})  \ar[r]^{\simeq} \ar[d]
& (\pi_1)_! (i_1)^* (\on{Fr}_{I_1})^*  (\mc F_{G, N, I, W}^{(I_2, \cdots, I_1)})   \ar[r]_{\quad \simeq}^{\quad F_{I_1}}  
& (\pi_1)_! (i_1)^* \mc F_{G, N, I, W}^{(I_1, \cdots, I_k)}  \ar[d]  \\
(\on{Fr}_{I_1})^* (\mc F_{M, N, I, W}^{' \, (I_2, \cdots, I_1)})   \ar[rr]_{\simeq}^{F_{I_1}} 
& & \mc F_{M, N, I, W}^{' \, (I_1, \cdots, I_k)}
}
$$
(Here we use the fact that the trace morphism is compatible with base change, as in (\ref{eqaution-Delta-CT-I-to-CT-J-functor}).)

We deduce that the following diagram is commutative:
$$
\xymatrix{
\restr{\mc H_{G, N, I, W}^j}{\on{Frob}_{I_1} (\ov{\eta^I})}    \ar[r]^{F_{I_1}}  \ar[d]^{\mc C_{G}^{P, \, \nu}}
& \restr{\mc H_{G, N, I, W}^j}{\ov{\eta^I}}  \ar[d]^{\mc C_{G}^{P, \, \nu} } \\
\restr{\mc H_{M, N, I, W}^{' \, j, \, \nu}}{\on{Frob}_{I_1} (\ov{\eta^I})}  \ar[r]^{F_{I_1}}
& \restr{\mc H_{M, N, I, W}^{' \, j, \, \nu}}{\ov{\eta^I}} 
}
$$ 
\cqfd

\begin{rem} 
In Section \ref{subsection-excursion-operator-on-cohomology}, for $J$ a finite set, $V$ a representation of $\wh G^J$ and $j \in \Z$, we defined the excursion operators acting on $H_{G, N, J, V}^j = \restr{\mc  H_{G, N, J, V}^j  }{\Delta^J(\ov{\eta})}$. Similarly, for any $\nu \in \wh \Lambda_{Z_M / Z_G}^{\Q}$, we can define excursion operators acting on $H_{M, N, J, V}^{' \, j, \, \nu} = \restr{\mc  H_{M, N, J, V}^{' \, j, \, \nu}  }{\Delta^J(\ov{\eta})}$. 
The same arguments as in the proof of Proposition \ref{prop-CT-commutes-with-excursion} prove that the following diagram is commutative:
$$
\xymatrixrowsep{2pc}
\xymatrixcolsep{6pc}
\xymatrix{
H_{G, N, J, V}^j  \ar[r]^{\quad S^G_{I, W, x, \xi, (\gamma_i)_{i \in I}} \quad}  \ar[d]^{C_{G}^{P, \, \nu}}
& H_{G, N, J, V}^j   \ar[d]^{C_{G}^{P, \, \nu}} \\
H_{M, N, J, V}^{' \, j, \, \nu}  \ar[r]^{\quad S^M_{I, W, x, \xi, (\gamma_i)_{i \in I}} \quad}
& H_{M, N, J, V}^{' \, j, \, \nu}
}
$$
where the constant term morphism $C_{G}^{P, \, \nu}$ is defined in \cite[Remark~3.5.11]{cusp-coho}.
\end{rem}

\subsection{Compatibility of Langlands parametrizations}

\sssec{}   \label{subsection-CT-commute-with-excursion-op}
Let $P$ be a parabolic subgroup of $G$ and $M$ its Levi quotient.
For any finite set $I$ and any function $f \in \mc O(\wh G_{\Qlbar} \backslash (\wh G_{\Qlbar} )^I / \wh G_{\Qlbar} )$, let $f^M$ be the composition $\wh M_{\Qlbar} \backslash (\wh M_{\Qlbar})^I / \wh M_{\Qlbar} \hookrightarrow \wh G_{\Qlbar} \backslash (\wh G_{\Qlbar})^I / \wh G_{\Qlbar} \xrightarrow{f} {\Qlbar}$. For any $(\gamma_i)_{i \in I} \in \on{Weil}(\eta, \ov{\eta})^I$, Proposition \ref{prop-CT-commutes-with-excursion} implies that the following diagram is commutative:
$$
\xymatrixrowsep{2pc}
\xymatrixcolsep{5pc}
\xymatrix{
C_c(\Bun_{G, N}(\Fq) / \Xi, \Qlbar)   \ar[r]^{\quad S^G_{I, f, (\gamma_i)_{i \in I}} \quad}  \ar[d]^{C_{G}^{P}}
& C_c(\Bun_{G, N}(\Fq) / \Xi, \Qlbar) \ar[d]^{C_{G}^{P}} \\
\prod_{\nu \in \wh{\Lambda}_{Z_M/ Z_G}^{\Q} } C_c(\Bun_{M, N}^{' \, \nu}(\Fq) / \Xi, \Qlbar)  \ar[r]^{\quad S^M_{I, f^M, (\gamma_i)_{i \in I}} \quad}
& \prod_{\nu \in \wh{\Lambda}_{Z_M/ Z_G}^{\Q} } C_c(\Bun_{M, N}^{' \, \nu}(\Fq) / \Xi, \Qlbar) .
}
$$

\sssec{}
Let $u \in |X \sm N|$.
Let $\ms I$ be an ideal of $\ms H_{G, u}$ of finite codimension.
We denote by $$\mc Q_{G, \ms I}:= C_c(\Bun_{G, N}(\Fq) / \Xi, \Qlbar) \big{/}  \ms I \cdot C_c(\Bun_{G, N}(\Fq) / \Xi, \Qlbar) $$ the quotient vector space.
In Definition \ref{def-algebra-excursion-B-I}, we defined ${\mc B}_{\ms I} \subset \on{End}_{\ms H_{G} } (  \mc Q_{G, \ms I} )$ the algebra of excursion operators. Here we denote it by ${\mc B}_{\ms I}^G$. We have Theorem \ref{thm-decomposition-of-C-c-quotient-I-by-param-Langlands} for $\mc Q_{G, \ms I}$.

\sssec{}
Recall that $\Xi$ is a lattice in $Z_G(F) \backslash Z_G(\mb A)$. 
Let $\Xi_M$ be a lattice in $Z_M(F) \backslash Z_M(\mb A)$ small enough. Then for any $\nu \in \wh{\Lambda}_{Z_M/ Z_G}^{\Q} $, the composition
$$\Bun_{M, N}^{\nu} / \Xi \rightarrow \Bun_{M, N} / \Xi \rightarrow \Bun_{M, N} / \Xi_M$$
is an open and closed immersion. We deduce that 
\begin{equation}
C_c(\Bun_{M, N}^{' \, \nu}(\Fq) / \Xi, \Qlbar) \hookrightarrow C_c(\Bun_{M, N}^{'}(\Fq) / \Xi_M, \Qlbar) .
\end{equation}

Applying Proposition \ref{prop-space-of-auto-form-is-Hecke-mod-type-fini} to $M$, we deduce that $C_c(\Bun_{M, N}^{'}(\Fq) / \Xi_M, \Qlbar)$ is of finite type as $\ms H_{M, u}$-module. Together with the fact that the Hecke operator $h^M_{\omega}$ acts on $\prod_{\nu} C_c(\Bun_{M, N}^{' \, \nu}(\Fq) / \Xi, \Qlbar)$ by translating the component indexed by $\nu$ (\emph{cf.} Section \ref{subsection-action-of-Hecke}), we deduce that $\prod_{\nu} C_c(\Bun_{M, N}^{' \, \nu}(\Fq) / \Xi, \Qlbar)$ is of finite type as $\ms H_{M, u}$-module. 

\sssec{}
Since $\ms I$ is of finite codimension in $\ms H_{G, u}$ and $\ms H_{G, u} \hookrightarrow \ms H_{M, u}$ is of finite type, we deduce that $\ms I \cdot \ms H_{M, u}$ is of finite codimension in $\ms H_{M, u}$. Thus the quotient vector space $$\mc Q_{M, \ms I}:= \prod_{\nu} C_c(\Bun_{M, N}^{' \, \nu}(\Fq) / \Xi, \Qlbar) \big{/}  (\ms I \cdot \ms H_{M, u}) \cdot \prod_{\nu} C_c(\Bun_{M, N}^{' \, \nu}(\Fq) / \Xi, \Qlbar) $$ is of finite dimension.

We denote by $${\mc B}_{\ms I}^{M, G} \subset \on{End}_{\ms H_{M}} \big(   \mc Q_{M, \ms I}     \big)$$ the sub-$\Qlbar$-algebra generated by all the excursion operators $S_{I, f^M, (\gamma_i)_{i \in I}}^M$, where $f \in \mc O(\wh G_{\Qlbar} \backslash (\wh G_{\Qlbar} )^I / \wh G_{\Qlbar} )$.
Similarly to Theorem \ref{thm-decomposition-of-C-c-quotient-I-by-param-Langlands}, there is a canonical decomposition of $\ms H_{M}$-modules
\begin{equation}    \label{equation-Q-M-m-sum-h-M-rho}
 \mc Q_{M, \ms I} = \bigoplus_{\rho} \mf H^{M, G}_{\rho}
\end{equation}
where the direct sum is indexed by $\wh G(\Qlbar)$-conjugacy classes of morphisms $\rho: \on{Weil}(\ov F / F) \rightarrow \wh G(\Qlbar)$ defined over a finite extension of $\mathbb{Q}_{\ell}$, continuous, semisimple and unramified outside $N$. The decomposition is characterized by the following property: $\mf H^{M, G}_{\rho}$ is equal to the generalized eigenspace associated to the character $\nu$ of ${\mc B}_{\ms I}^{M, G}$ defined by $\nu(S^M_{I, f^M, (\gamma_i)_{i \in I}}) = f( ( \rho(\gamma_i) )_{i \in I} ).$

\sssec{}
We denote by $${\mc B}_{\ms I}^M \subset \on{End}_{\ms H_{M}} \big(   \mc Q_{M, \ms I}     \big)$$ the sub-$\Qlbar$-algebra generated by all the excursion operators $S_{I, g, (\gamma_i)_{i \in I}}^M$, where $g \in \mc O(\wh M_{\Qlbar} \backslash (\wh M_{\Qlbar})^I / \wh M_{\Qlbar} )$.
Similarly to Theorem \ref{thm-decomposition-of-C-c-quotient-I-by-param-Langlands}, there is a canonical decomposition of $\ms H_{M}$-modules
\begin{equation}      \label{equation-Q-M-m-sum-h-M-rho-prime}
 \mc Q_{M, \ms I} = \bigoplus_{\rho'} \mf H^{M}_{\rho'}
\end{equation}
where the direct sum is indexed by $\wh M(\Qlbar)$-conjugacy classes of morphisms $\rho': \on{Weil}(\ov F / F) \rightarrow \wh M(\Qlbar)$ defined over a finite extension of $\mathbb{Q}_{\ell}$, continuous, semisimple and unramified outside $N$. And $\mf H^{M}_{\rho'}$ is equal to the generalized eigenspace associated to the character $\nu$ of ${\mc B}_{\ms I}^{M}$ defined by $\nu(S^M_{I, g, (\gamma_i)_{i \in I}}) = g(  (\rho'(\gamma_i))_{i \in I}  ).$

Denote by $j: \wh M(\Qlbar) \hookrightarrow \wh G(\Qlbar)$ the inclusion.
We have 
\begin{equation}
\mf H^{M, G}_{\rho} = \underset{\rho': \on{Weil}(\ov F / F) \rightarrow \wh M(\Qlbar), \;  j \circ \rho' = \rho }{\bigoplus} \mf H^{M}_{\rho'} ,
\end{equation}
where $j \circ \rho' = \rho$ is up to $\wh G(\Qlbar)$-conjugacy.

\sssec{}
We have a commutative diagram as $\ms H_{G}$-modules:
$$
\xymatrixrowsep{3pc}
\xymatrixcolsep{6pc}
\xymatrix{
\mc Q_{G, \ms I} \ar[r]^{ =  }_{ \text{Theorem } \ref{thm-decomposition-of-C-c-quotient-I-by-param-Langlands}}  \ar[d]^{\ov{C_G^P}}
& \underset{\rho}{\bigoplus} \mf H^G_{\rho}  \ar[d]^{\ov{C_G^P}}   \\
\mc Q_{M, \ms I}   \ar[r]^{ = }_{ (\ref{equation-Q-M-m-sum-h-M-rho})  }
&  \underset{\rho}{\bigoplus} \mf H^{M, G}_{\rho}  .
}
$$

By \ref{subsection-CT-commute-with-excursion-op}, for $f \in \mc O(\wh G_{\Qlbar} \backslash (\wh G_{\Qlbar})^I / \wh G_{\Qlbar} )$, the excursion operator $S_{I, f^M, (\gamma_i)}^M$ acts on $\ov{C_G^P}(\mf H^G_{\rho})$ by $f((\rho(\gamma_i))_{i \in I})$. As a consequence, if $\ov{C_G^P}(\mf H^G_{\rho}) \neq 0$, then $\rho: \on{Weil}(\ov F / F) \rightarrow \wh G(\Qlbar)$ factors through
$$\on{Weil}(\ov F / F) \xrightarrow{\rho'} \wh M(\Qlbar) \rightarrow \wh G(\Qlbar)$$
for some $\rho'$ (which may not be unique).  We conclude that in this case (as $\ms H_{M}$-modules) 
$$\ov{C_G^P}(\mf H^G_{\rho}) \subset \mf H^{M, G}_{\rho}  = \underset{\rho': \on{Weil}(\ov F / F) \rightarrow \wh M(\Qlbar), \;  j \circ \rho' = \rho }{\bigoplus} \mf H^M_{\rho'}  .$$

\quad

%%%%%%%%%%%%%%%%%%%%%
% References
%%%%%%%%%%%%%%%%%%%%%

\end{document}